\documentclass{amsart}

\usepackage{amssymb}
\usepackage{eepic}
\usepackage{color}

\allowdisplaybreaks

\numberwithin{equation}{section}

\newtheorem{theorem}{Theorem}[section]
\newtheorem{lemma}[theorem]{Lemma}
\newtheorem{proposition}[theorem]{Proposition}

\newtheorem{remark}[theorem]{Remark}

\newtheorem{Atheorem}{Theorem}

\newtheorem{Acorollary}[Atheorem]{Corollary}

\newtheorem{Aremark}[Atheorem]{Remark}

\newtheorem{TheoA}{Theorem A}
\newtheorem{TheoB}{Theorem B}

\newtheorem{localprob}{An $L_2$-localization problem}
\newtheorem{localtheo}{A pseudo-localization principle}
\newtheorem{Intcor}{Corollary}
\newtheorem{Cuculescutheo}{Cuculescu's construction}
\newtheorem{shifttheo}{Shifted \boldmath $T1$ \unboldmath theorem}
\newtheorem{cotlar}{Cotlar lemma}
\newtheorem{schur}{Schur lemma}
\newtheorem{localest}{A localization estimate}

\newcommand{\N}{\mathbb{N}}
\newcommand{\Z}{\mathbb{Z}}
\newcommand{\R}{\mathbb{R}}
\newcommand{\C}{\mathbb{C}}

\newcommand{\summ}{\sum\nolimits}

\newcommand{\Mn}{\mathcal{A}}

\def\M{\mathcal{M}}
\def\Q{\mathcal{Q}}
\def\1{\mathbf{1}}

\newcommand{\dem}{\noindent {\bf Proof. }}

\newcommand{\skeB}{\noindent {\bf Sketch of the proof of Theorem B. }}
\newcommand{\fin}{\hspace*{\fill} $\square$ \vskip0.2cm}

\def\esssup{\mathop{\mathrm{ess \ sup}}}

\addtocounter{tocdepth}{-1}

\begin{document}

\title[Noncommutative Calder{\'o}n-Zygmund
theory] {Pseudo-localization of singular integrals \\ and
noncommutative Calder{\'o}n-Zygmund theory}

\author[Javier Parcet]
{Javier Parcet}

\maketitle

\vskip-0.5cm

\null

\tableofcontents

\addtolength{\parskip}{+1ex}

\vskip-2.5cm

\footnote{Partially supported by \lq Programa Ram{\'o}n y Cajal,
2005\rq${}$ and \\ \indent also by Grants MTM2004-00678 and
CCG06-UAM/ESP-0286, Spain.}



\section*{Introduction}

After the pioneer work of Calder{\'o}n and Zygmund in the 50's, the
systematic study of singular integrals has become a corner stone
in harmonic analysis with deep implications in mathematical
physics, partial differential equations and other mathematical
disciplines. Subsequent generalizations of Calder{\'o}n-Zygmund theory
have essentially pursued two lines. We may either consider more
general domains or ranges for the functions considered. In the
first case, the Euclidean space is replaced by metric spaces
equipped with a doubling or non-doubling measure of polynomial
growth. In the second case, the real or complex fields are
replaced by a Banach space in which martingale differences are
unconditional. Historically, the study of singular integrals
acting on matrix or operator valued functions has been considered
part of the vector-valued theory. This is however a limited
approach in the noncommutative setting and we propose to regard
these functions as operators in a suitable von Neumann algebra,
generalizing so the domain and not the range of classical
functions. A far reaching aspect of our approach is the stability
of the product $fg$ and the absolute value $|f| = \sqrt{f^*f}$ for
operator-valued functions, a fundamental property not exploited in
the vector theory. In this paper we follow the original
Calder{\'o}n-Zygmund program and present a non-commutative
scalar-valued Calder{\'o}n-Zygmund theory, emancipated from the vector
theory.

Noncommutative harmonic analysis (understood in a wide sense) has
received much attention in recent years. The functional analytic
approach given by operator space theory and the new methods from
quantum/free probability have allowed to study a great variety of
topics. We find in the recent literature noncommutative analogs of
Khintchine and Rosenthal inequalities, a settled noncommutative
theory of martingale inequalities, new results on Fourier/Schur
multipliers, matrix $A_p$ weights and a sharpened Carleson
embedding theorem, see \cite{Ha,JPX,JX2,NPTV,PX1,Vo} and the
references therein. However, no essential progress has been made
in the context of singular integral operators.

Our original motivation was the weak type boundedness of
Calder{\'o}n-Zygmund operators acting on operator-valued functions, a
well-known problem which has remained open since the beginning of
the vector-valued theory in the 80's. This fits in the context of
Mei's recent paper \cite{Me}. Our main tools for its solution are
two. On one hand, the failure of some classical estimates in the
noncommutative setting forces us to have a deep understanding of
how the $L_2$-mass of a singular integral is concentrated around
the support of the function on which it acts. To that aim, we have
developed a \emph{pseudo-localization principle} for singular
integrals which is of independent interest, even in the classical
theory. This is used in conjunction with a \emph{noncommutative
form of Calder{\'o}n-Zygmund decomposition} which we have constructed
using the theory of noncommutative martingales. As a byproduct of
our weak type inequality, we obtain the sharp asymptotic behavior
of the constants for the strong $L_p$ inequalities as $p \to 1$
and $p \to \infty$, which are not known. At the end of the paper
we generalize our results to certain singular integrals including
operator-valued kernels and functions at the same time. A deep
knowledge of this kind of \emph{fully noncommutative} operators is
a central aim in noncommutative harmonic analysis. Our methods in
this paper open a door to work in the future with more general
classes of operators.

\noindent \textbf{1. Terminology.} Let us fix some notation that
will remain fixed all through out the paper. Let $\M$ be a
semifinite von Neumann algebra equipped with a normal semifinite
faithful trace $\tau$. Let us consider the algebra $\Mn_B$ of
essentially bounded $\M$-valued functions
$$\Mn_B = \Big\{ f: \R^n \to \M \, \big| \, f \ \mbox{strongly
measurable s.t.} \ \esssup_{x \in \R^n} \|f(x)\|_{\M} < \infty
\Big\},$$ equipped with the \emph{n.s.f.} trace
$$\varphi(f) = \int_{\R^n} \tau(f(x)) \, dx.$$ The weak-operator
closure $\Mn$ of $\Mn_B$ is a von Neumann algebra. If $1 \le p \le
\infty$, we write $L_p(\M)$ and $L_p(\Mn)$ for the noncommutative
$L_p$ spaces associated to the pairs $(\M,\tau)$ and
$(\Mn,\varphi)$. The lattices of projections are written $\M_\pi$
and $\Mn_{\pi}$, while $\1_\M$ and $\1_{\Mn}$ stand for the unit
elements.

The set of dyadic cubes in $\R^n$ is denoted by $\Q$. The size of
any cube $Q$ in $\R^n$ is defined as the length $\ell(Q)$ of one
of its edges. Given an integer $k \in \Z$, we use $\Q_k$ for the
subset of $\Q$ formed by cubes $Q$ of the $k$-th generation, i.e.
those of size $1/2^{k}$. If $Q$ is a dyadic cube and $f: \R^n \to
\M$ is integrable on $Q$, we set the average
$$f_Q = \frac{1}{|Q|} \int_Q f(y) \, dy.$$
Let us write $(\mathsf{E}_k)_{k \in \Z}$ for the family of
conditional expectations associated to the classical dyadic
filtration on $\R^n$. $\mathsf{E}_k$ will also stand for the
tensor product $\mathsf{E}_k \otimes id_\M$ acting on $\Mn$. If $1
\le p < \infty$ and $f \in L_p(\Mn)$
$$\mathsf{E}_k(f) = f_k = \sum_{Q \in \Q_k}^{\null} f_Q 1_Q.$$
We shall denote by $(\Mn_k)_{k \in \Z}$ the corresponding
filtration $\Mn_k = \mathsf{E}_k(\Mn)$.

If $Q \in \Q$, its dyadic father $\widehat{Q}$ is the only dyadic
cube containing $Q$ with double size. Given $\delta
> 1$, the $\delta$-concentric father of $Q$ is the only cube
$\delta \hskip1pt Q$ concentric with the cube $Q$ and such that
$\ell(\delta \hskip1pt Q) = \delta \hskip1pt \ell(Q)$. In this
paper we will mainly work with dyadic and $9$-concentric fathers.
Note that in the classical theory $2$-concentric fathers are
typically enough. We shall write just $L_p$ to refer to the
commutative $L_p$ space on $\R^n$ equipped with the Lebesgue
measure $dx$.

\noindent \textbf{2. Statement of the problem.} Just to motivate
our problem and for the sake of simplicity, the reader may think
for the moment that $(\M,\tau)$ is given by the pair
$(M_m,\mbox{tr})$ formed by the algebra of $m \times m$ square
matrices equipped with the standard trace. In this particular
case, the von Neumann algebra $\Mn = \Mn_B$ becomes the space of
essentially bounded matrix-valued functions. Let us consider a
Calder{\'o}n-Zygmund operator formally given by
$$T \! f(x) = \int_{\R^n} k(x,y) f(y) \, dy.$$ As above, let
$L_p(\M)$ be the noncommutative $L_p$ space associated to
$(\mathcal{M},\tau)$. If $\mathcal{M}$ is the algebra of $m \times
m$ matrices we recover the Schatten $p$-class over $M_m$, for a
general definition see below. The first question which arises is
whether or not the singular integral $T$ is bounded on
$L_p(\mathcal{A})$ for $1 < p < \infty$. The space
$L_p(\mathcal{A})$ is defined as the closure of $\mathcal{A}_B$
with respect to the norm $$\|f\|_p = \Big( \int_{\R^n} \tau \,
\big( |f(x)|^p \big) \, dx \Big)^{\frac1p}.$$ In other words,
$L_p(\mathcal{A})$ is isometric to the Bochner $L_p$ space with
values in $L_p(\mathcal{M})$. In particular, when dealing with the
Hilbert transform and by a well-known result of Burkholder
\cite{Bu,Bu2}, the boundedness on $L_p(\mathcal{A})$ reduces to
the fact that $L_p(\M)$ is a UMD Banach space for $1 < p <
\infty$, see also \cite{Bo,Bo3}. After some partial results of
Bourgain \cite{Bo2}, it was finally Figiel \cite{Fi} who showed in
1989 (using an ingenious martingale approach) that the UMD
property implies the $L_p$ boundedness of the corresponding
vector-valued singular integrals associated to generalized
kernels.

The second natural question has to do with a suitable weak type
inequality for $p=1$. Namely, such inequality is typically
combined in the classical theory with the real interpolation
method to produce extrapolation results on the $L_p$ boundedness
of Calder{\'o}n-Zygmund and other related operators. The problem of
finding the right weak type inequality is subtler since arguments
from the vector-valued theory are no longer at our disposal.
Indeed, in terms of Bochner spaces we may generalize the previous
situation by considering the mapping $T$ from $L_1(\R^n;
\mathrm{X})$ to $L_{1,\infty}(\R^n; \mathrm{X})$ with $\mathrm{X}
= L_1(\M)$. However, $L_1(\M)$ is not UMD and the resulting
operator is not bounded. On the contrary, using operators rather
than vectors (i.e. working directly on the algebra $\mathcal{A}$)
we may consider the operator $T: L_1(\mathcal{A}) \to
L_{1,\infty}(\mathcal{A})$ where $L_{1,\infty}(\mathcal{A})$
denotes the corresponding noncommutative Lorentz space, to be
defined below. The only result on this line is the weak type
$(1,1)$ boundedness of the Hilbert transform for operator-valued
functions, proved by Randrianantoanina in \cite{R1}. He followed
Kolmogorov's approach, exploiting the conjugation nature of the
Hilbert transform (defined in a very wide setting via Arveson's
\cite{Ar} maximal subdiagonal algebras) and applying complex
variable methods. As is well-known this is no longer valid for
other Calder{\'o}n-Zygmund operators and new real variable methods are
needed. In the classical case, these methods live around the
celebrated Calder{\'o}n-Zygmund decomposition. One of the main
purposes of this paper is to supply the right real variable
methods in the noncommutative context. As we will see, there are
significant differences.

Using real interpolation, our main result gives an extrapolation
method which produces the $L_p$ boundedness results discussed in
the paragraph above and provides the sharp asymptotic behavior of
the constants, for which the UMD approach is inefficient.
Moreover, when working with operator-valued kernels we obtain new
strong $L_p$ inequalities. We should warn the reader not to
confuse this setting with that of Rubio de Francia, Ruiz and
Torrea \cite{RRT}, Hyt\"onen \cite{Hy} and Hyt\"onen/Weis
\cite{HW1,HW2}, where the mentioned limitations of the
vector-valued theory appear.

\noindent \textbf{3. Calder{\'o}n-Zygmund decomposition.} Let us
recall the formulation of the classical decomposition for
scalar-valued integrable functions. If $f \in L_1$ is positive and
$\lambda \in \R_+$, we consider the level set
$$\mathrm{E}_\lambda = \Big\{ x \in \R^n \, \big|
\, M_d f(x) > \lambda \Big\},$$ where the dyadic Hardy-Littlewood
maximal function $M_df$ is greater than $\lambda$. If we write
$\mathrm{E}_\lambda = \bigcup_j Q_j$ as a disjoint union of
maximal dyadic cubes, we may decompose $f = g+b$ where the good
and bad parts are given by $$g = f 1_{\mathrm{E}_{\lambda}^c} +
\summ_j^{\null} f_{Q_j} 1_{Q_j} \quad \mbox{and} \quad b =
\summ_j^{\null} (f - f_{Q_j}) 1_{Q_j}$$ Letting $b_j = (f -
f_{Q_j}) 1_{Q_j}$, we have
\begin{itemize}
\item[i)] $\|g\|_1 \le \|f\|_1$ and $\|g\|_\infty \le 2^n
\lambda$.

\item[ii)] $\mbox{supp} \, b_j \subset Q_j$, $\int_{Q_j} b_j = 0$
and $\sum_j \|b_j\|_1 \le 2 \|f\|_1$.
\end{itemize}
These properties are crucial for the analysis of singular integral
operators.

In this paper we use the so-called Cuculescu's construction
\cite{Cu} to produce a sequence $(p_k)_{k \in \Z}$ of disjoint
projections in $\mathcal{A}$ which constitute the noncommutative
counterpart of the characteristic functions supported by the sets
$$\mathrm{E}_\lambda(k) = \bigcup_{\begin{subarray}{c} Q_j
\subset \mathrm{E}_\lambda \\ \ell(Q_j) = 1/2^k
\end{subarray}} Q_j.$$ Cuculescu's construction will be
properly introduced in the text. It has proved to be the right
tool from the theory of noncommutative martingales to deal with
inequalities of weak type. Indeed, Cuculescu proved in \cite{Cu}
the noncommutative Doob's maximal weak type inequality. Moreover,
these techniques were used by Randrianantoanina to prove several
weak type inequalities for noncommutative martingales
\cite{R2,R3,R4} and by Junge and Xu in their remarkable paper
\cite{JX2}. In fact, a strong motivation for this paper relies on
\cite{PR}, where similar methods were applied to obtain Gundy's
decomposition for noncommutative martingales. It is well-known
that the probabilistic analog of Calder{\'o}n-Zygmund decomposition is
precisely Gundy's decomposition. However, in contrast to the
classical theory, the noncommutative analogue of Calder{\'o}n-Zygmund
decomposition turns out to be much harder than Gundy's
decomposition. Although we shall justify this below in further
detail, the main reason is that singular integral operators do not
localize the support of the function on which it acts, something
that happens for instance with martingale transforms or martingale
square functions.

Let us now formulate the noncommutative Calder{\'o}n-Zygmund
decomposition. If $f \in L_1(\Mn)_+$ and $\lambda \in \R_+$, we
consider the disjoint projections $(p_k)_{k \in \Z}$ given by
Cuculescu's construction. Let $p_{\infty}$ denote the projection
onto the ortho-complement of the range of $\sum_k p_k$. In
particular, using the terminology $\widehat{\Z} = \Z \cup
\{\infty\}$ we find the relation
$$\sum_{k \in \widehat{\Z}}^{\null} \hskip1pt p_k =
\mathbf{1}_{\Mn}.$$ Then, the \emph{good} and \emph{bad} parts are
given by $$g = \sum_{i,j \in \widehat{\Z}}^{\null} p_i f_{i \vee
j} p_j \quad \mbox{and} \quad b = \sum_{i,j \in
\widehat{\Z}}^{\null} p_i (f - f_{i \vee j}) p_j,$$ with $i \vee j
= \max (i,j)$. We will show how this generalizes the classical
decomposition.

\noindent \textbf{4. Main weak type inequality.} Let $\Delta$
denote the diagonal of $\R^n \times \R^n$. We will write in what
follows $T$ to denote a linear map $\mathcal{S} \to \mathcal{S}'$
from test functions to distributions which is associated to a
given kernel $k: \R^{2n} \setminus \Delta \to \C$. This means that
for any smooth test function $f$ with compact support, we have
$$Tf(x) = \int_{\R^n} k(x,y) f(y) \, dy \quad \mbox{for all} \quad
x \notin \mbox{supp} \hskip1pt f.$$ Given two points $x,y \in
\R^n$, the distance $|x-y|$ between $x$ and $y$ will be taken for
convenience with respect to the $\ell_\infty(n)$ metric. As usual,
we impose size and smoothness conditions on the kernel:
\begin{itemize}
\item[a)] If $x,y \in \R^n$, we have $$|k(x,y)| \ \lesssim \
\frac{1}{|x-y|^n}.$$

\item[b)] There exists $0 < \gamma \le 1$ such that
$$\begin{array}{rcl} \big| k(x,y) - k(x',y) \big| & \lesssim &
\displaystyle \frac{|x-x'|^\gamma}{|x-y|^{n+\gamma}} \quad
\mbox{if} \quad |x-x'| \le \frac12 \hskip1pt |x-y|, \\ [10pt]
\big| k(x,y) - k(x,y') \big| & \lesssim & \displaystyle
\frac{|y-y'|^\gamma}{|x-y|^{n+\gamma}} \quad \mbox{if} \quad
|y-y'| \le \frac12 \hskip1pt |x-y|.
\end{array}$$
\end{itemize}
We will refer to this $\gamma$ as the Lipschitz smoothness
parameter of the kernel.

\begin{TheoA} \label{MainThA}
Let $T$ be a generalized Calder{\'o}n-Zygmund operator associated to a
kernel satisfying the size and smoothness estimates above. Assume
that $T$ is bounded on $L_q$ for some $1 < q < \infty$. Then,
given any $f \in L_1(\Mn)$, the estimate below holds for some
constant $\mathrm{c}_{n,\gamma}$ depending only on the dimension
$n$ and the Lipschitz smoothness parameter $\gamma$
$$\sup_{\lambda > 0} \lambda \, \varphi \Big\{ |Tf| > \lambda
\Big\} \le \mathrm{c}_{n,\gamma} \, \|f\|_1.$$ In particular,
given $1 < p < \infty$ and $f \in L_p(\Mn)$, we find $$\|T \! f
\|_p \le \mathrm{c}_{n,\gamma} \, \frac{p^2}{p-1} \, \|f\|_p.$$
\end{TheoA}

The expression $\sup_{\lambda > 0} \lambda \, \varphi \big\{ |Tf|
> \lambda \big\}$ is just a slight abuse of notation to denote the
noncommutative weak $L_1$ norm, to be rigorously defined below. We
find it though more intuitive, since it is reminiscent of the
classical terminology. Theorem A provides a positive answer to our
problem for any singular integral associated to a generalized
Calder{\'o}n-Zygmund kernel satisfying the size/smoothness conditions
imposed above. Moreover, the asymptotic behavior of the constants
as $p \to 1$ and $p \to \infty$ is optimal. Independently, Tao Mei
has recently obtained another argument for this which does not
include the weak type inequality \cite{TaoMei}. We shall present
it at the end of the paper, since we shall use it indirectly to
obtain weak type inequalities for singular integrals associated to
operator-valued kernels. In the language of operator space theory
and following Pisier's characterization \cite{Pis,Pis2} of
complete boundedness we immediately obtain:

\begin{Intcor}
Let $T$ be a generalized $L_q$-bounded and $\gamma$-Lipschitz
Calder{\'o}n-Zygmund operator. Let us equip $L_p$ with its natural
operator space structure. Then, the cb-norm of $T: L_p \to L_p$ is
controlled by $$\mathrm{c}_{n,\gamma} \, \frac{p^2}{p-1}.$$ Thus,
the growth rate as $p \to 1$ or $p \to \infty$ coincides with the
Banach space case.
\end{Intcor}

\noindent Before going on, a few remarks are in order:

\textbf{a)} It is standard to reduce the proof of Theorem A to the
case $q=2$.

\textbf{b)} The reader might think that our hypothesis on
Lipschitz smoothness for the first variable is unnecessary to
obtain the weak type inequality and that only smoothness with
respect to the second variable is needed. Namely, this is the case
in the classical theory. It is however not the case in this paper
because the use of certain almost orthogonality arguments (see
below) forces us to apply both kinds of smoothness. We refer to
Remark \ref{Onlypoint} for the specific point where the
$x$-Lipschitz smoothness is applied and to Remark
\ref{RemWeakHypo} for more in depth discussion on the conditions
imposed on the kernel.

\textbf{c)} In the classical case $\mathrm{E}_\lambda$ is a
perfectly delimited region of $\R^n$. In particular, we may
construct the dilation $9 \hskip1pt \mathrm{E}_\lambda = \bigcup_j
9 \hskip1pt Q_j$. This set is useful to estimate the bad part $b$
since it has two crucial properties. First, it is small because
$|9 \, \mathrm{E}_\lambda| \sim |\mathrm{E}_\lambda|$ and
$\mathrm{E}_\lambda$ satisfies the Hardy-Littlewood weak maximal
inequality. Second, its complement is far away from
$\mathrm{E}_\lambda$ (the support of $b$) so that $Tb$ restricted
to $\R^n \setminus 9 \, \mathrm{E}_\lambda$ avoids the singularity
of the kernel. The problem that we find in the noncommutative case
is that $\mathrm{E}_\lambda$ is no longer a region in $\R^n$.
Indeed, given a dyadic cube $Q$ and a positive $f \in L_1$, we
have either $f_Q > \lambda$ or not and this dichotomy completely
determines the set $\mathrm{E}_\lambda$. However, for $f \in
L_1(\Mn)_+$ the average $f_Q$ is a positive operator (not a
positive number) and the dichotomy disappears since the condition
$f_Q > \lambda$ is only satisfied in part of the spectrum of
$f_Q$. This difficulty is inherent to the noncommutativity and is
motivated by the lack of a total order in the positive cone of
$\M$. It also produces difficulties to define noncommutative
maximal functions \cite{Cu,J1}, a problem that required the recent
theory of operator spaces for its solution and is in the heart of
the matter. Our construction of the right noncommutative analog
$\zeta$ of $\R^n \setminus 9 \hskip1pt \mathrm{E}_\lambda$ is a
key step in this paper, see Lemma \ref{keylem} below. Here it is
relevant to recall that, quite unexpectedly (in contrast with the
classical case) we shall need the projection $\zeta$ to deal with
\emph{both} the good and the bad parts.

\textbf{d)} Another crucial difference with the classical setting
and maybe the hardest point to overcome is the lack of estimates
i) and ii) above in the noncommutative framework. Indeed, given $f
\in L_1(\Mn)_+$ we only have such estimates for the diagonal terms
$$\summ_k p_k f_k p_k \quad \mbox{and} \quad \summ_k p_k (f-f_k)
p_k.$$ A more detailed discussion on this topic is given in
Appendix B below. Let us now explain how we face the lack of the
classical inequalities. Since $\1_\Mn - \zeta$ is the
noncommutative analog of $9 \hskip1pt \mathrm{E}_\lambda$ which is
\emph{small} as explained above, we can use the noncommutative
Hardy-Littlewood weak maximal inequality to reduce our problem to
estimate the terms $\zeta \hskip1pt T(g) \hskip1pt \zeta$ and
$\zeta \hskip1pt T(b) \hskip1pt \zeta$. A very naive and formally
incorrect way to explain what to do here is the following. Given a
fixed positive integer $s$, we find \emph{something like}
\begin{eqnarray*} \Big\| \zeta \, T \Big( \sum_{|i-j|=s} p_i
f_{i \vee j} p_j \Big) \zeta \Big\|_2 & \lesssim & s \hskip1pt
2^{- \gamma s} \, \Big\| \summ_k p_k f_k p_k \Big\|_2,
\\ \Big\| \zeta \, T \Big( \sum_{|i-j|=s} p_i (f
- f_{i \vee j}) p_j \Big) \zeta \Big\|_1 & \lesssim & s \hskip1pt
2^{-\gamma s} \, \Big\| \summ_k p_k (f-f_k) p_k \Big\|_1,
\end{eqnarray*}
where $\gamma$ is the Lipschitz smoothness parameter of the
kernel. In other words, we may estimate the action of $\zeta
\hskip1pt T (\hskip1pt \cdot \hskip1pt) \zeta$ on the terms in the
$s$-th upper and lower diagonals by $s \hskip1pt 2^{- \gamma s}$
times the \emph{corresponding size} of the main diagonal. Then,
recalling that i) and ii) hold on the diagonal, it is standard to
complete the argument. We urge however the reader to understand
this just as a motivation (not as a claim) since the argument is
quite more involved than this. For instance, we will need to
replace the off-diagonal terms of $g$ by other $g_{k,s}$'s
satisfying
$$\sum_{k,s} g_{k,s} = \sum_{i \neq j} p_i f_{i \vee j} p_j.$$

\vskip-5pt

A rough way of rephrasing this phenomenon is to say that
Calder{\'o}n-Zygmund operators are \emph{almost diagonal} when acting
on operator-valued functions. In other contexts, this almost
diagonal nature has already appeared in the literature. The
wavelet proof of the $T1$ theorem \cite{MC} exhibits this property
of singular integrals with respect to the Haar system in $\R^n$.
This also applies in the context of Clifford analysis \cite{Mi}.
Moreover, some deep results in \cite{Ch,DRdF} (which we will
comment below) use this almost diagonal nature as a key idea. It
is also worthy of mention that these difficulties do not appear in
\cite{PR}. The reason is that the operators for which Gundy's
decomposition is typically applied (martingale transforms or
martingale square functions) do not move the support of the
original function/operator. This means that the action of $T$ over
the off-diagonal terms is essentially supported by $\1_\Mn -
\zeta$. Consequently, these terms are controlled by means of the
noncommutative analog of Doob's maximal weak type inequality, see
\cite{PR} for further details. As a byproduct, we observe that the
pseudo-localization principle which we present below is not needed
in \cite{PR}.

\noindent \textbf{5. Pseudo-localization.} A key point in our
argument is the behavior of singular integrals acting on the
off-diagonal terms $p_i f_{i \vee j} p_j$ and $p_i (f-f_{i \vee
j}) p_j$, as a function of the parameter $s = |i-j|$ in a region
$\zeta \approx \R^n \setminus 9 \hskip1pt \mathrm{E}_\lambda$
which is in some sense far away from their (left and right)
support. The idea we need to exploit relies on the following
principle: more regularity of the kernel of $T$ implies a faster
decay of $T \! f$ far away from the support of $f$. That is why
the $b$-terms $p_i (f-f_{i \vee j}) p_j$ are better than the
$g$-terms $p_i f_{i \vee j} p_j$. Indeed, the cancellation of $f -
f_{i \vee j}$ allows to subtract a piecewise constant function
from the kernel (in the standard way) to apply the smoothness
properties of it and obtain suitable $L_1$ estimates. However, the
off-diagonal $g$-terms are not mean-zero (at least at first sight)
with respect to $\int_{\R^n}$ and we need more involved tools to
prove this pseudo-localization property in the $L_2$ metric. For
the sake of clarity and since our result might be of independent
interest even in the classical theory, we state it for
scalar-valued functions. The way we apply it in our noncommutative
setting will be clarified along the text, see Theorem
\ref{pseudolocal} for the noncommutative form of this principle.

Since we are assuming that $T$ is bounded on $L_2$, we may further
assume by homogeneity that it is of norm $1$. In the sequel, we
will only consider $L_2$-normalized Calder{\'o}n-Zygmund operators.
Our result is related to the following problem.

\begin{localprob}
Given $f: \R^n \to \C$ in $L_2$ and $0 < \delta < 1$, find the
sets $\Sigma_{f,\delta}$ such that the inequality below holds for
all normalized Calder{\'o}n-Zygmund operator satisfying the imposed
size/smoothness conditions
$$\Big( \int_{\R^n \setminus \Sigma_{f \! , \delta}} |T \! f(x)|^2
dx \Big)^{\frac12} \, \le \, \delta \hskip1pt \Big( \int_{\R^n}
|f(x)|^2 dx \Big)^{\frac12}.$$
\end{localprob}

Given $f: \R^n \to \C$ in $L_2$, let $f_k$ and $df_k$ denote the
$k$-th condition expectation of $f$ with respect to the standard
dyadic filtration and its corresponding $k$-th martingale
difference. That is, we have $df_k = \sum_{Q \in \Q_k} \big( f_Q -
f_{\widehat{Q}} \big) 1_Q.$ Let $\mathcal{R}_k$ be the class of
sets in $\R^n$ being the union of a family of cubes in $\Q_k$.
Given such an $\mathcal{R}_k$-set $\Omega = \bigcup_j Q_j$, we
shall work with the dilations $9 \Omega = \bigcup_j 9 Q_j$, where
$9 Q$ denotes the $9$-concentric father of $Q$. We shall prove the
following result.

\begin{localtheo}
Let us fix a positive integer $s$. Given a function $f$ in $L_2$
and any integer $k$, we define $\Omega_k$ to be the smallest
$\mathcal{R}_k$-set containing the support of $df_{k+s}$. If we
further consider the set $$\Sigma_{f,s} = \bigcup_{k \in \Z} 9
\hskip1pt \Omega_k,$$ then we have the localization estimate
$$\Big( \int_{\R^n \setminus \hskip1pt
\Sigma_{f \! ,s}} |T \! f(x)|^2 \, dx \Big)^{\frac12} \ \le \
\mathrm{c}_{n,\gamma} \hskip1pt s \hskip1pt 2^{- \gamma s/4} \Big(
\int_{\R^n} |f(x)|^2 dx \Big)^{\frac12},$$ for any
$L_2$-normalized Calder{\'o}n-Zygmund operator with Lipschitz
parameter $\gamma$.
\end{localtheo}

Given any integer $k \in \Z$, we are considering the smallest set
$\Omega_k$ containing $\mathrm{supp} \hskip2pt df_{k+s}$ and
belonging to an $s$ times \emph{coarser} topology. This procedure
gives rise to an apparently artificial \emph{shift condition}
$$\mathrm{supp} \hskip2pt df_{k+s} \subset \Omega_k$$ which gives
a measure of how much we should enlarge $\Sigma_f = \bigcup_k
\mathrm{supp} \hskip2pt df_k$ at every scale. However, this
condition (or its noncommutative analog) is quite natural in our
setting since it is satisfied by the off-diagonal terms of $g$,
precisely those for which our previous tools did not work. In the
classical/commutative setting there are some natural situations
for which our result applies and some others which limit the
applicability of it. For instance, at first sight our result is
only applicable for functions $f$ satisfying $f_m =0$ for some
integer $m$. There are also some other natural questions such as
an $L_p$ analog or our result or an equivalent formulation using a
Littlewood-Paley decomposition, instead of martingale differences.
For the sake of clarity in our exposition, we prove the result in
the body of the paper and we postpone these further comments to
Appendix A below.

The proof of this result reduces to a \emph{shifted form of the
$T1$ theorem} in a sense to be explained below. In particular,
almost orthogonality methods are essential in our approach.
Compared to the standard proofs of the $T1$ theorem, with wavelets
\cite{MC} or more generally with approximations of the identity
\cite{St}, we need to work in a dyadic/martingale setting forced
by the role of Cuculescu's construction in this paper. This
produces a lack of smoothness in the functions we work with,
requiring quite involved estimates to obtain almost orthogonality
results. An apparently new aspect of our estimates is the
asymmetry of our bounds when applying Schur lemma, see Remark
\ref{Asimetria} for more details.

Let us briefly comment the relation of our result with two papers
by Christ \cite{Ch} and Duoandikoetxea/Rubio de Francia
\cite{DRdF}. Although both papers already exploited the almost
diagonal nature of Calder{\'o}n-Zygmund operators, only
convolution-type singular integrals are considered and no
localization result is pursued there. Being more specific, a
factor $2^{-\gamma s}$ is obtained in \cite{Ch} for the bad part
of Calder{\'o}n-Zygmund decomposition. As explained above, we need to
produce this factor for the good part. This is very unusual (or
even new) in the literature. Nevertheless, the way we have stated
our pseudo-localization result shows that the key property is the
shift condition $\mathrm{supp} \hskip2pt df_{k+s} \subset
\Omega_k$, regardless we work with good or bad parts. On the other
hand, in \cite{DRdF} Littlewood-Paley theory and the commutativity
produced by the use of convolution operators is used to obtain
related estimates in $L_p$ with $p \neq 2$. In particular, almost
orthogonality does not play any role there. The lack of a suitable
noncommutative Littlewood-Paley theory and our use of generalized
Calder{\'o}n-Zygmund operators make their argument not applicable
here.

\noindent \textbf{6. Operator-valued kernels.} At the end of the
paper we extend our main results to certain Calder{\'o}n-Zygmund
operators associated to kernels $k: \R^{2n} \setminus \Delta \to
\M$ satisfying the canonical size/smoothness conditions. In other
words, we replace the absolute value by the norm in $\M$:
\begin{itemize}
\item[a)] If $x,y \in \R^n$, we have $$\|k(x,y)\|_\M \ \lesssim \
\frac{1}{|x-y|^n}.$$

\item[b)] There exists $0 < \gamma \le 1$ such that
$$\begin{array}{rcl} \big\| k(x,y) - k(x',y) \big\|_\M & \lesssim &
\displaystyle \frac{|x-x'|^\gamma}{|x-y|^{n+\gamma}} \quad
\mbox{if} \quad |x-x'| \le \frac12 \hskip1pt |x-y|, \\ [10pt]
\big\| k(x,y) - k(x,y') \big\|_\M & \lesssim & \displaystyle
\frac{|y-y'|^\gamma}{|x-y|^{n+\gamma}} \quad \mbox{if} \quad
|y-y'| \le \frac12 \hskip1pt |x-y|.
\end{array}$$
\end{itemize}
Unfortunately, not every such kernel satisfies the analog of
Theorem A. Namely, we shall construct (using classical
Littlewood-Paley methods) a simple kernel satisfying the size and
smoothness conditions above and giving rise to a Calder{\'o}n-Zygmund
operator bounded on $L_2(\Mn)$ but not on $L_p(\Mn)$ for $1 < p <
2$. However, a detailed inspection of our proof of Theorem A and a
few auxiliary results will show that the key condition (together
with the size/smoothness hypotheses on the kernel) for the
operator $T$ is to be an $\M$-bimodule map. Of course, this always
holds in the context of Theorem A. When dealing with
operator-valued kernels this is false in general, but it holds for
instance when dealing with \emph{standard} Calder{\'o}n-Zygmund
operators $$T \! f(x) \, = \, \xi f(x) \, + \, \lim_{\varepsilon
\to 0} \int_{|x-y| > \varepsilon} k(x,y) \hskip1pt f(y) \, dy$$
associated to a commuting kernel $k: \R^{2n} \setminus \Delta \to
\mathcal{Z_M}$, with $\mathcal{Z_M} = \M \cap \M'$ standing for
the center of $\M$. Note that we are only requiring the
$\M$-bimodule property to hold on the singular integral part,
since the multiplier part is always well-behaved as far as $\xi
\in \Mn$. Note also that when $\M$ is a factor, any commuting
kernel must be scalar-valued and we go back to Theorem A.

\begin{TheoB} \label{MainThB}
Let $T$ be a generalized Calder{\'o}n-Zygmund operator associated to
an operator-valued kernel $k: \R^{2n} \setminus \Delta \to \M$
satisfying the imposed size/smoothness conditions. Assume that $T$
is an $\M$-bimodule map bounded on $L_q(\Mn)$ for some $1 < q <
\infty$. Then, the following weak type inequality holds for some
constant $\mathrm{c}_{n,\gamma}$ depending only on the dimension
$n$ and the Lipschitz smoothness parameter $\gamma$
$$\sup_{\lambda > 0} \lambda \, \varphi \Big\{ |Tf| > \lambda
\Big\} \le \mathrm{c}_{n,\gamma} \, \|f\|_1.$$ In particular,
given $1 < p < \infty$ and $f \in L_p(\Mn)$, we find $$\|T \! f
\|_p \le \mathrm{c}_{n,\gamma} \, \frac{p^2}{p-1} \, \|f\|_p.$$
\end{TheoB}

The strong $L_p$ inequalities stated in Theorem B do not follow
from a UMD-type argument as it happened with Theorem A. In
particular, these $L_p$ estimates seem to be new and independently
obtained by Tao Mei as pointed above.

\noindent \textbf{7. Appendices.} We conclude the paper with two
appendices. A further analysis on pseudo-localization is given in
Appendix A. This mainly includes remarks related to our result,
some conjectures on possible generalizations and a corollary on
the rate of decreasing of the $L_2$ mass of a singular integral
far away from the support of the function on which it acts. In
Appendix B we study the noncommutative form of Calder{\'o}n-Zygmund
decomposition in further detail. In particular, we give some
weighted inequalities for the good and bad parts which generalize
the classical $L_1$ and $L_2$ estimates satisfied by these
functions. The sharpness of our estimates remains as an open
interesting question.

\noindent \textbf{Remark.} The value of the constant
$\mathrm{c}_{n,\gamma}$ will change from one instance to another.

\noindent \textbf{Acknowledgement.} I would like to thank J.M.
Martell, F. Soria and Q. Xu for discussions related to the content
of this paper and specially to Tao Mei for keeping me up to date
on his related work.

\section{Noncommutative integration} \label{S1}

We begin with a quick survey of definitions and results on
noncommutative $L_p$ spaces and related topics that will be used
along the paper. All or most of it will be well-known to experts
in the field. The right framework for a noncommutative analog of
measure theory and integration is von Neumann algebra theory. We
refer to \cite{KR,Ta} for a systematic study of von Neumann
algebras and to the recent survey by Pisier/Xu \cite{PX2} for a
detailed exposition of noncommutative $L_p$ spaces.

\subsection{Noncommutative $L_p$}

A \emph{von Neumann algebra} is a weak-operator closed
$\mathrm{C}^*$-algebra. By the Gelfand-Naimark-Segal theorem, any
von Neumann algebra $\M$ can be embedded in the algebra
$\mathcal{B}(\mathcal{H})$ of bounded linear operators on some
Hilbert space $\mathcal{H}$. In what follows we will identify $\M$
with a subalgebra of $\mathcal{B(H)}$. The positive cone $\M_+$ is
the set of positive operators in $\M$. A \emph{trace} $\tau: \M_+
\to [0,\infty]$ on $\M$ is a linear map satisfying the tracial
property $\tau(a^*a) = \tau(aa^*)$. It is said to be \emph{normal}
if $\sup_\alpha \tau(a_\alpha) = \tau(\sup_\alpha a_\alpha)$ for
any bounded increasing net $(a_\alpha)$ in $\M_+$; it is
\emph{semifinite} if for any non-zero $a \in \M_+$, there exists
$0 < a' \le a$ such that $\tau(a') < \infty$ and it is
\emph{faithful} if $\tau(a) = 0$ implies $a = 0$. Taking into
account that $\tau$ plays the role of the integral in measure
theory, all these properties are quite familiar. A von Neumann
algebra $\M$ is called \emph{semifinite} whenever it admits a
normal semifinite faithful (\emph{n.s.f.} in short) trace $\tau$.
Except for a brief comment in Remark \ref{TypeIII} below we shall
always work with semifinite von Neumann algebras. Recalling that
any operator $a$ can be written as a linear combination $a_1 - a_2
+ ia_3 - ia_4$ of four positive operators, we can extend $\tau$ to
the whole algebra $\M$. Then, the tracial property can be restated
in the familiar way $\tau(ab) = \tau(ba)$ for all $a,b \in \M$.

According to the GNS construction, it is easily seen that the
noncommutative analogs of measurable sets (or equivalently
characteristic functions of those sets) are orthogonal
projections. Given $a \in \M_+$, the support projection of $a$ is
defined as the least projection $q$ in $\M$ such that $qa = a =
aq$ and will be denoted by $\mbox{supp} \hskip1pt a$. Let
$\mathcal{S}_+$ be the set of all $a \in \M_+$ such that
$\tau(\mbox{supp} \hskip1pt a) < \infty$ and set $\mathcal{S}$ to
be the linear span of $\mathcal{S}_+$. If we write $|x|$ for the
operator $(x^*x)^{\frac12}$, we can use the spectral measure
$\gamma_{|x|}: \R_+ \to \mathcal{B}(\mathcal{H})$ of the operator
$|x|$ to define
$$|x|^p = \int_{\R_+} s^p \, d \gamma_{|x|}(s) \quad
\mbox{for} \quad 0 < p < \infty.$$ We have $x \in \mathcal{S}
\Rightarrow |x|^p \in \mathcal{S}_+ \Rightarrow \tau(|x|^p) <
\infty$. If we set $\|x\|_p = \tau( |x|^p )^{\frac1p}$, it turns
out that $\| \ \|_p$ is a norm in $\mathcal{S}$ for $1 \le p <
\infty$ and a $p$-norm for $0 < p < 1$. Using that $\mathcal{S}$
is a $w^*$-dense $*$-subalgebra of $\M$, we define the
\emph{noncommutative $L_p$ space} $L_p(\M)$ associated to the pair
$(\M, \tau)$ as the completion of $(\mathcal{S}, \| \ \|_p)$. On
the other hand, we set $L_\infty(\M) = \M$ equipped with the
operator norm. Many of the fundamental properties of classical
$L_p$ spaces like duality, real and complex interpolation... can
be transferred to this setting. The most important properties for
our purposes are the following:
\begin{itemize}
\item H\"older inequality. If $1/r = 1/p+1/q$, we have $\|ab\|_r
\le \|a\|_p \|b\|_q$.

\item The trace $\tau$ extends to a continuous functional on
$L_1(\M)$: $|\tau(x)| \le \|x\|_1$.
\end{itemize}
We refer to \cite{PX2} for a definition of $L_p$ over
non-semifinite von Neumann algebras.

\subsection{Noncommutative symmetric spaces}

Let $$\M' = \Big\{ b \in \mathcal{B}(\mathcal{H}) \ \big| \, ab =
ba \ \mbox{for all} \ a \in \M \Big\}$$ be the commutant of $\M$.
A closed densely-defined operator on $\mathcal{H}$ is
\emph{affiliated} with $\M$ when it commutes with every unitary
$u$ in the commutant $\M'$. Recall that $\M = \M''$ and this
implies that every $a \in \M$ is affiliated with $\M$. The
converse fails in general since we may find unbounded operators.
If $a$ is a densely defined self-adjoint operator on $\mathcal{H}$
and $a = \int_{\R} s \hskip1pt d \gamma_a(s)$ is its spectral
decomposition, the spectral projection $\int_{\mathcal{R}} d
\gamma_a(s)$ will be denoted by $\chi_\mathcal{R}(a)$. An operator
$a$ affiliated with $\mathcal{M}$ is \emph{$\tau$-measurable} if
there exists $s > 0$ such that $$\tau \big( \chi_{(s,\infty)}
(|a|) \big) = \tau \big\{ |a| > s \big\} < \infty.$$ The
\emph{generalized singular-value} $\mu(a): \R_+ \to \R_+$ is
defined by
$$\mu_t (a) = \inf \Big\{ s
> 0 \, \big| \ \tau \big\{ |x| > s \big\} \le t \Big\}.$$
This provides us with a noncommutative analogue of the so-called
non-increasing rearrangement of a given function. We refer to
\cite{FK} for a detailed exposition of the function $\mu(a)$ and
the corresponding notion of convergence in measure.

If $L_0(\M)$ denotes the $*$-algebra of $\tau$-measurable
operators, we have the following equivalent definition of $L_p$
$$L_p(\M) = \Big\{a \in L_0(\M) \, \big| \ \Big( \int_{\R_+}
\mu_t(a)^p \, dt \Big)^{\frac1p} < \infty \Big\}.$$ The same
procedure applies to symmetric spaces. Given the pair $(\M,\tau)$,
let $\mathrm{X}$ be a rearrangement invariant quasi-Banach
function space on the interval $(0, \tau(\mathbf{1}_\M))$. The
\emph{noncommutative symmetric space} $\mathrm{X}(\mathcal{M})$ is
defined by
$$\mathrm{X}(\mathcal{M}) = \Big\{a \in L_0(\mathcal{M}) \,
\big| \ \mu(a) \in \mathrm{X} \Big\} \quad \text{with} \quad
\left\| a \right\|_{\mathrm{X}(\M)} = \|\mu(a)\|_{\mathrm{X}}.$$
It is known that $\mathrm{X}(\mathcal{M})$ is a Banach (resp.
quasi-Banach) space whenever $\mathrm{X}$ is a Banach (resp.
quasi-Banach) function space. We refer the reader to
\cite{DDdP,X1} for more in depth discussion of this construction.
Our interest in this paper is restricted to noncommutative
$L_p$-spaces and \emph{noncommutative weak $L_1$-spaces}.
Following the construction of symmetric spaces of measurable
operators, the noncommutative weak $L_1$-space
$L_{1,\infty}(\mathcal{M})$,  is defined as the set of all $a$ in
$L_0(\mathcal{M})$ for which the quasi-norm
\[ \left\|a\right\|_{1,\infty} = \sup_{t > 0} \, t \hskip1pt
\mu_t(x) = \sup_{\lambda > 0} \, \lambda \hskip1pt \tau \Big\{ |x|
> \lambda \Big\}\] is finite. As in the commutative case, the
noncommutative weak $L_1$ space satisfies a quasi-triangle
inequality that will be used below with no further reference.
Indeed, the following inequality holds for $a_1, a_2 \in
L_{1,\infty}(\M)$ $$\lambda \, \tau \Big\{ |a_1+a_2| > \lambda
\Big\} \le \lambda \, \tau \Big\{ |a_1| > \lambda/2 \Big\} +
\lambda \, \tau \Big\{ |a_2| > \lambda/2 \Big\}.$$

\subsection{Noncommutative martingales}

Consider a von Neumann subalgebra (a weak$^*$ closed
$*$-subalgebra) $\mathcal{N}$ of $\mathcal{M}$. A
\emph{conditional expectation} $\mathcal{E}: \mathcal{M} \to
\mathcal{N}$ from $\mathcal{M}$ onto $\mathcal{N}$ is a positive
contractive projection. The conditional expectation $\mathcal{E}$
is called \emph{normal} if the adjoint map $\mathcal{E}^*$
satisfies $\mathcal{E}^*(\mathcal{M}_*) \subset \mathcal{N}_*$. In
this case, there is a map $\mathcal{E}_*: \mathcal{M}_*
\rightarrow \mathcal{N}_*$ whose adjoint is $\mathcal{E}$. Note
that such normal conditional expectation exists if and only if the
restriction of $\tau$ to the von Neumann subalgebra $\mathcal{N}$
remains semifinite, see e.g. Theorem 3.4 in \cite{Ta}. Any such
conditional expectation is trace preserving (i.e. $\tau \circ
\mathcal{E} = \tau$) and satisfies the bimodule property
\[ \mathcal{E}(a_1 b \hskip1pt a_2) = a_1 \mathcal{E}(b) \hskip1pt
a_2 \quad \mbox{for all} \quad a_1, a_2 \in \mathcal{N} \
\mbox{and} \ b \in \mathcal{M}.\]

Let $(\mathcal{M}_k)_{k \ge 1}$  be an increasing sequence of von
Neumann subalgebras of $\mathcal{M}$ such that the union of the
$\mathcal{M}_k$'s is weak$^*$ dense in $\mathcal{M}$. Assume that
for every $k \ge 1$, there is a normal conditional expectation
$\mathcal{E}_k: \mathcal{M} \to \mathcal{M}_k$. Note that for
every $1 \le p < \infty$ and $k \ge 1$, $\mathcal{E}_k$ extends to
a positive contraction $\mathcal{E}_k: L_p(\mathcal{M}) \to
L_p(\mathcal{M}_k)$. A \emph{noncommutative martingale} with
respect to the filtration $(\mathcal{M}_k)_{k \ge 1}$ is a
sequence $a = (a_k)_{k \ge 1}$ in $L_1(\mathcal{M})$ such that
\[
\mathcal{E}_j(a_k) = a_j \quad \mbox{for all} \quad 1 \le j \le k
< \infty.
\]
If additionally $a \subset L_p(\mathcal{M})$ for some $1 \le p \le
\infty$ and $\|a\|_p = \sup_{k \ge 1} \|a_k\|_p < \infty$, then
$a$ is called an \emph{$L_p$-bounded martingale}. Given a
martingale $a = (a_k)_{k \ge 1}$, we assume the convention that
$a_0 = 0$. Then, the martingale difference sequence $da =(da_k)_{k
\ge 1}$ associated to $x$ is defined by $da_k = a_k - a_{k-1}$.

The next result due to Cuculescu \cite{Cu} was the first known
result in the theory and will be crucial in this paper. It can be
viewed as a noncommutative analogue of the classical weak type
$(1,1)$ boundedness of Doob's maximal function.

\begin{Cuculescutheo}
Suppose $a = (a_1, a_2, \ldots)$ is a positive $L_1$ martingale
relative to the filtration $(\mathcal{M}_k)_{k \ge 1}$ and let
$\lambda$ be a positive number. Then there exists a decreasing
sequence of projections
$$q(\lambda)_1, q(\lambda)_2, q(\lambda)_3, \ldots$$ in
$\mathcal{M}$ satisfying the following properties
\begin{itemize}
\item[i)] $q(\lambda)_k$ commutes with $q(\lambda)_{k-1} a_k
q(\lambda)_{k-1}$ for each $k \ge 1$.

\item[ii)] $q(\lambda)_k$ belongs to $\mathcal{M}_k$ for each $k
\ge 1$ and $q(\lambda)_k a_k q(\lambda)_k \le \lambda \hskip1pt
q(\lambda)_k$.

\item[iii)] The following estimate holds $$\tau \Big(
\mathbf{1}_\M - \bigwedge_{k \ge 1} q(\lambda)_k \Big) \le
\frac{1}{\lambda} \hskip1pt \sup_{k \ge 1} \|a_k\|_1.$$
\end{itemize}
Explicitly, we set $q(\lambda)_0 = \1_\M$ and define $q(\lambda)_k
= \chi_{(0,\lambda]}(q(\lambda)_{k-1} a_k q(\lambda)_{k-1})$.
\end{Cuculescutheo}

The theory of noncommutative martingales has achieved considerable
progress in recent years. The renewed interest on this topic
started from the fundamental paper of Pisier and Xu \cite{PX1},
where they introduced a new functional analytic approach to study
Hardy spaces and the Burkholder-Gundy inequalities for
noncommutative martingales. Shortly after, many classical
inequalities have been transferred to the noncommutative setting.
A noncommutative analogue of Doob's maximal function \cite{J1},
the noncommutative John-Nirenberg theorem \cite{JM}, extensions of
Burkholder inequalities for conditioned square functions
\cite{JX1} and related weak type inequalities \cite{R2,R3,R4}; see
\cite{PR} for a simpler approach to some of them.

\section{A pseudo-localization principle}

Let us now proceed with the proof of the pseudo-localization
principle stated in the Introduction. In the course of it we will
see the link with a shifted form of the $T1$ theorem, which is
formulated in a dyadic martingale setting. Since we are concerned
with its applications to our noncommutative problem, we leave a
more in depth analysis of our result to Appendix A below.

\subsection{Three auxiliary results}
\label{3AR}

We need some well-known results that live around David-Journ{\'e}'s
$T1$ theorem. Cotlar lemma is very well-known and its proof can be
found in \cite{Duo,St}. We include the proof of Schur lemma, since
our statement and proof is non-standard, see Remark
\ref{Asimetria} below for details. The localization estimate at
the end follows from \cite{MC}. We give the proof for
completeness.

\begin{cotlar}
Let $\mathcal{H}$ be a Hilbert space and let us consider a family
$(T_k)_{k \in \Z}$ of bounded operators on $\mathcal{H}$ with
finitely many non-zero $T_k$'s. Assume that there exists a
summable sequence $(\alpha_k)_{k \in \Z}$ such that
$$\max \Big\{ \big\| T_i^* T_j^{\null} \big\|_{\mathcal{B(H)}},
\big\| T_i^{\null} T_j^* \big\|_{\mathcal{B(H)}} \Big\} \, \le \,
\alpha_{i-j}^2$$ for all $i,j \in \Z$. Then we automatically have
$$\Big\| \summ_k T_k \Big\|_{\mathcal{B(H)}} \, \le \, \summ_k
\alpha_k.$$
\end{cotlar}

\begin{schur}
Let $T$ be given by $$T \! f(x) = \int_{\R^n} k(x,y) \hskip1pt
f(y) \, dy.$$ Let us define the Schur integrals associated to $k$
\begin{eqnarray*}
\mathcal{S}_1(x) & = & \int_{\R^n} \big| k(x,y) \big| \, dy,
\\ \mathcal{S}_2(y) & = & \int_{\R^n} \big| k(x,y)
\big| \, dx.
\end{eqnarray*}
Assume that both $\mathcal{S}_1$ and $\mathcal{S}_2$ belong to
$L_\infty$. Then, $T$ is bounded on $L_2$ and
$$\|T\|_{\mathcal{B}(L_2)} \le \sqrt{ \big\| \mathcal{S}_1 \big\|_\infty
\big\| \mathcal{S}_2 \big\|_\infty^{\null} }.$$
\end{schur}

\dem By the Cauchy-Schwarz inequality we obtain
\begin{eqnarray*}
\lefteqn{\hskip-25pt \Big( \int_{\R^n} \Big| \int_{\R^n} k(x,y)
\hskip1pt f(y) \, dy \Big|^2 dx \Big)^{\frac12}} \\ & \le & \Big(
\int_{\R^n} \Big[ \int_{\R^n} |k(x,y)| \hskip1pt |f(y)| \, dy
\Big]^2 dx \Big)^{\frac12} \\ & \le & \Big( \int_{\R^n} \Big[
\int_{\R^n} |k(x,y)| \, dy \Big] \Big[ \int_{\R^n} |k(x,y)|
\hskip1pt |f(y)|^2 \, dy \Big] \, dx \Big)^{\frac12} \\ & \le &
\sqrt{\big\| \mathcal{S}_1 \big\|_\infty} \, \Big( \int_{\R^n}
\int_{\R^n} |k(x,y)| \hskip1pt |f(y)|^2 \, dy \, dx
\Big)^{\frac12} \\ & \le & \sqrt{\big\| \mathcal{S}_1
\big\|_\infty \big\| \mathcal{S}_2 \big\|_\infty} \, \Big(
\int_{\R^n} |f(y)|^2 \, dy \Big)^{\frac12}.
\end{eqnarray*}
\fin

\begin{remark} \label{Asimetria}
\emph{Typically, Schur lemma is formulated as
$$\|T\|_{\mathcal{B}(L_2)} \le \frac{1}{2} \Big( \big\| \mathcal{S}_1
\big\|_\infty + \big\| \mathcal{S}_2 \big\|_\infty \Big),$$ see
e.g. \cite{MC,St}. This might happen because we usually have
$\|\mathcal{S}_1\|_\infty \sim \|\mathcal{S}_2\|_\infty$, by
certain symmetry in the estimates. In particular, the cases for
which the arithmetic mean does not help but the geometric mean
does are very rare in the literature, or even (as far as we know)
not existent! However, motivated by a lack of symmetry in our
estimates, this is exactly the case in this paper.}
\end{remark}

\begin{localest} Assume that $$|k(x,y)| \lesssim
\frac{1}{|x-y|^n} \quad \mbox{for all} \quad x,y \in \R^n.$$ Let
$T$ be a Calder{\'o}n-Zygmund operator associated to the kernel $k$
and assume that $T$ is $L_2$-normalized. Then, given $x_0 \in
\R^n$ and $r_1, r_2 \in \R_+$ with $r_2 > 2 \hskip1pt r_1$, the
estimate below holds for any pair $f,g$ of bounded scalar-valued
functions respectively supported by $\mathsf{B}_{r_1}(x_0)$ and
$\mathsf{B}_{r_2}(x_0)$
$$\big| \big\langle T \! f, g \big\rangle \big| \le \mathrm{c}_n
\hskip1pt r_1^n \hskip1pt \log(r_2/r_1) \hskip1pt \|f\|_\infty
\|g\|_\infty.$$
\end{localest}

\dem Let us write $\mathsf{B}$ for the ball
$\mathsf{B}_{3r_1/2}(x_0)$ and let us consider a smooth function
$\rho$ which is identically $1$ on $\mathsf{B}$ and identically
$0$ outside $\mathsf{B}_{2r_1}(x_0)$. Set $\eta = 1-\rho$ so that
we may decompose
$$\big\langle T \! f, g \big\rangle = \big\langle T \! f, \rho g
\big\rangle + \big\langle T \! f, \eta g \big\rangle.$$ For the
first term we have
\begin{eqnarray*}
\big| \big\langle T \! f, \rho g \big\rangle \big| & \le & \|T \!
f\|_2 \|\rho g\|_2 \ \le \ \|f\|_2 \|\rho g\|_2 \\
& \le & \|f\|_\infty \|g\|_\infty \sqrt{\big| \mathrm{supp} f
\big| \big| \mathrm{supp} (\rho g) \big|} \ \le \ \mathrm{c}_n
\hskip1pt r_1^n \hskip1pt \|f\|_\infty \|g\|_\infty.
\end{eqnarray*}
On the other hand, for the second term we have
$$|\langle T \! f, \eta g \rangle| = \Big| \int_{\mathsf{B}_{r_2}(x_0)
\setminus \mathsf{B}} \Big( \int_{\mathsf{B}_{r_1}(x_0)} k(x,y)
\hskip1pt f(y) \, dy \Big) \hskip2pt \overline{\eta g}(x) \, dx
\Big|.$$ The latter integral is clearly bounded by
$$\|f\|_\infty \|g\|_\infty \hskip1pt
\int_{\Omega} \frac{dx \, dy}{|x-y|^n}$$ with $\Omega =
(\mathsf{B}_{r_2}(x_0) \setminus \mathsf{B}) \times
\mathsf{B}_{r_1}(x_0)$. However, it is easily checked that an
upper bound for the double integral given above is provided by
$\mathrm{c}_n \hskip1pt r_1^n \hskip1pt \log(r_2/r_1)$, where
$\mathrm{c}_n$ is a constant depending only on $n$. This completes
the proof. \fin

\subsection{Shifted $T1$ theorem}

By the conditions imposed on $T$ in the Introduction, it is clear
that its adjoint $T^*$ is an $L_2$-normalized Calder{\'o}n-Zygmund
operator with kernel $k^*(x,y) = \overline{k(y,x)}$ satisfying the
same size and smoothness estimates. This implies that $T^*1$
(understood in a weak sense, see e.g. \cite{St} for details)
belongs to $\mathrm{BMO}$, the space of functions with bounded
mean oscillation. In addition, if $\Delta_j = \mathsf{E}_j -
\mathsf{E}_{j-1}$ denotes the dyadic martingale difference
operator, it is also well known that for any $\rho \in
\mathrm{BMO}$ the dyadic paraproduct against $\rho$
$$\Pi_\rho(f) = \sum_{j=-\infty}^\infty \Delta_j(\rho)
\mathsf{E}_{j-1}(f)$$ is bounded on $L_2$. Here it is necessary to
know how $\mathrm{BMO}$ is related to its dyadic version
$\mathrm{BMO}_d$, see \cite{GJ} and \cite{Mei} for details. It is
clear that $\Pi_\rho(1) = \rho$ and the adjoint of $\Pi_\rho$ is
given by the operator
$$\Pi_\rho^*(f) = \sum_{j=-\infty}^\infty \mathsf{E}_{j-1} \big(
\overline{\Delta_j(\rho)} f \big).$$ Thus, since $T^*1 \in
\mathrm{BMO}$ we may write $$T = T_0 + \Pi_{T^*1}^*.$$ According
to our previous considerations, both $T_0$ and $\Pi_{T^*1}^*$ are
Calder{\'o}n-Zygmund operators bounded on $L_2$ and their kernels
satisfy the standard size and smoothness conditions imposed on $T$
with the same Lipschitz smoothness parameter $\gamma$, see
\cite{St} for the latter assertion. Moreover, the operator $T_0$
now satisfies $T_0^*1=0$. Now we use that $T_0^*1$ is the weak$^*$
limit of a sequence $(T_0^* \rho_k)_{k \ge 1}$ in $\mathrm{BMO}$,
where the $\rho_k$'s are increasing bump functions which converge
to $1$. In particular, the relation below holds for any $f \in
H_1$
\begin{equation} \label{T*1=0}
\int_{\R^n}^{\null} T_0 f (x) \, dx = 0.
\end{equation}
Indeed, we have $\langle T_0 f,1 \rangle = \langle f, T_0^*1
\rangle = 0$. The use of paraproducts is exploited in the $T1$
theorem to produce the cancellation condition \eqref{T*1=0}, which
is a key assumption to make Cotlar lemma effective in this
setting. The paraproduct term is typically estimated using
Carleson's lemma, although we will not need it here. What we shall
do is to prove that our theorem for $T_0$ and $\Pi_{T^*1}^*$
reduces to prove a shifted form of the $T1$ theorem. In this
paragraph we only deal with $T_0$.

Let $T$ be a generalized Calder{\'o}n-Zygmund operator as in the
statement of our result and assume that $T$ satisfies the
cancellation condition \eqref{T*1=0}, so that there is no need to
use the notation $T_0$ in what follows. Let us write $$\R^n
\setminus \Sigma_{f,s} = \bigcap_{k \in \Z} \Theta_k \qquad
\mbox{with} \qquad \R^n \setminus \Theta_k = 9 \hskip1pt
\Omega_k.$$ Denote by $\mathsf{E}_k$ the $k$-th dyadic conditional
expectation and by $\Delta_k$ the martingale difference operator
$\mathsf{E}_k - \mathsf{E}_{k-1}$, so that $\mathsf{E}_k(f) = f_k$
and $\Delta_k(f) = df_k$. Recall that $\Omega_k$ and $\Theta_k$
are $\mathcal{R}_k$-sets. In particular, the action of multiplying
by the characteristic functions $1_{\Omega_k}$ or $1_{\Theta_k}$
commutes with $\mathsf{E}_j$ for all $j \ge k$. Then we consider
the following decomposition \vskip-5pt $$1_{\R^n \setminus
\Sigma_{f,s}} T \! f = 1_{\R^n \setminus \Sigma_{f,s}} \Big(
\summ_k \mathsf{E}_k T \Delta_{k+s} 1_{\Omega_k} + \summ_k
(id-\mathsf{E}_k) 1_{\Theta_k} T 1_{\Omega_k} \Delta_{k+s} \Big)
(f).$$ \vskip3pt \noindent Note that we have used here the shift
condition $\mathrm{supp} \hskip1pt df_{k+s} \subset \Omega_k$ as
well as the commutation relations mentioned above in conjunction
with $\R^n \! \setminus \! \Sigma_{f,s} \subset \Theta_k$. Next we
observe that $1_{\Theta_k} T 1_{\Omega_k} = 1_{\Theta_k} T_{4
\cdot 2^{-k}} 1_{\Omega_k}$, where $T_{\varepsilon}$ denotes the
truncated singular integral formally given by $$T_\varepsilon f
(x) = \int_{|x-y| > \varepsilon} k(x,y) \hskip1pt f(y) \, dy.$$
Indeed, we have
$$1_{\Theta_k} T 1_{\Omega_k} f (x) \, = \, 1_{\Theta_k}(x)
\sum_{\begin{subarray}{c} Q \in \Q_k \\ Q \hskip1pt \cap \hskip1pt
\Omega_k \neq \emptyset \end{subarray}} 1_{\R^n \setminus 9Q}(x)
\int_Q k(x,y) \hskip1pt f(y) \, dy,$$ from where the claimed
identity follows, since we have $$\mathrm{dist}(Q, \R^n \!
\setminus \! 9Q) = 4 \cdot 2^{-k}$$ for all $Q \in \Q_k$. Taking
all these considerations into account, we deduce $$1_{\R^n
\setminus \Sigma_{f,s}} T \! f = 1_{\R^n \setminus \Sigma_{f,s}}
\Big( \summ_k \mathsf{E}_k T \Delta_{k+s} + \summ_k
(id-\mathsf{E}_k) T_{4 \cdot 2^{-k}} \Delta_{k+s} \Big) (f).$$ In
particular, our problem reduces to estimate the norm in
$\mathcal{B}(L_2)$ of \vskip-5pt
$$\Phi_s = \summ_k \mathsf{E}_k T \Delta_{k+s} \quad \mbox{and}
\quad \Psi_s = \summ_k (id-\mathsf{E}_k) T_{4 \cdot 2^{-k}}
\Delta_{k+s}.$$ \vskip3pt \noindent Both $\Phi_s$ and $\Psi_s$ are
reminiscent of well-known operators (in a sense $\mathsf{E}_k T$
and $T_{4 \cdot 2^{-k}}$ behave here in the same way) appearing in
the proof of the $T1$ theorem by David and Journ{\'e} \cite{DJ}.
Indeed, what we find (in the context of dyadic martingales) is
exactly the $s$-shifted analogs meaning that we replace $\Delta_k$
by $\Delta_{k+s}$. In summary, we have proved that under the
assumption that cancellation condition \eqref{T*1=0} holds our
main result reduces to the proof of the theorem below.

\begin{shifttheo}
Let $T$ be an $L_2$-normalized Calder{\'o}n-Zygmund operator with
Lipschitz parameter $\gamma$. Assume that $T^*1 = 0$ or, in other
words, that we have $\int_{\R^n}^{\null} T \! f(x) \, dx = 0$ for
any $f \in H_1$. Then, we have
$$\|\Phi_s\|_{\mathcal{B}(L_2)} = \Big\| \summ_k \mathsf{E}_k T
\Delta_{k+s} \Big\|_{\mathcal{B}(L_2)} \le \mathrm{c}_{n,\gamma}
\hskip1pt s \hskip1pt 2^{- \gamma s/4}.$$ Moreover, regardless the
value of $T^*1$ we also have $$\|\Psi_s\|_{\mathcal{B}(L_2)} =
\Big\| \summ_k (id-\mathsf{E}_k) T_{4 \cdot 2^{-k}} \Delta_{k+s}
\Big\|_{\mathcal{B}(L_2)} \le \mathrm{c}_{n,\gamma} \hskip1pt
2^{-\gamma s/2}.$$
\end{shifttheo}

\begin{remark}
\emph{For some time, our hope was to estimate $$\Big\| \summ_k
T_{4 \cdot 2^{-k}} \Delta_{k+s} \Big\|_{\mathcal{B}(L_2)}$$ since
we believed that the truncation of order $2^{-k}$ in conjunction
with the action of $\Delta_{k+s}$ was enough to produce the right
decay. Note that our pseudo-localization result could also be
deduced from this estimate. However, the cancellation produced by
the paraproduct decomposition in $\Phi_s$ and by the presence of
the term $id - \mathsf{E}_k$ in $\Psi_s$ play an essential role in
the argument.}
\end{remark}

\subsection{Paraproduct argument}
\label{PA}

Now we show how the estimate of the paraproduct term also reduces
to the shifted $T1$ theorem stated above. Indeed, let us write
$\Pi$ instead of $\Pi_{T^*1}^*$ to simplify the notation. Then, as
we did above, it is straightforward to see that \vskip-5pt
$$1_{\R^n \setminus \Sigma_{f,s}} \Pi f =
1_{\R^n \setminus \Sigma_{f,s}} \Big( \summ_k \mathsf{E}_k \Pi
\Delta_{k+s} 1_{\Omega_k} + \summ_k (id-\mathsf{E}_k) \Pi_{4 \cdot
2^{-k}} \Delta_{k+s} \Big) (f).$$ \vskip3pt \noindent Recalling
one more time that $\Pi$ is an $L_2$-bounded generalized
Calder{\'o}n-Zygmund operator satisfying the same size and smoothness
conditions as $T$, the estimate for the second operator $$\Big\|
\summ_k (id-\mathsf{E}_k) \Pi_{4 \cdot 2^{-k}} \Delta_{k+s}
\Big\|_{\mathcal{B}(L_2)} \, \le \, \mathrm{c}_{n,\gamma}
\hskip1pt 2^{- \gamma s/2}$$ follows from the second assertion of
the shifted $T1$ theorem. Here it is essential to note that the
hypothesis $T^*1 = 0$ is not needed for $\Psi_s$. Therefore, it
only remains to estimate the first operator. However, we claim
that $1_{\R^n \setminus \Sigma_{f,s}}  \summ_k \mathsf{E}_k \Pi
\Delta_{k+s} 1_{\Omega_k} f$ is identically zero. Let us prove
this assertion. We have $$1_{\R^n \setminus \Sigma_{f,s}} \summ_k
\mathsf{E}_k \Pi \Delta_{k+s} 1_{\Omega_k} f \, = \, 1_{\R^n
\setminus \Sigma_{f,s}} \summ_k \mathsf{E}_k \summ_j
\mathsf{E}_{j-1} \Big( \overline{\Delta_j(T^*1)} \hskip1pt
1_{\Omega_k} \hskip1pt df_{k+s} \Big).$$ If we fix the integer
$k$, all the $j$-terms on the second sum above vanish except for
the term associated to $j=k+s$. Indeed, if $j < k+s$ we use
$\mathsf{E}_{j-1} = \mathsf{E}_{j-1} \mathsf{E}_{k+s-1}$ and
obtain
\begin{eqnarray*}
\mathsf{E}_{j-1} \Big( \overline{\Delta_j(T^*1)} \hskip1pt
1_{\Omega_k} \hskip1pt df_{k+s} \Big) & = & \mathsf{E}_{j-1} \Big(
\mathsf{E}_{k+s-1} \big( \overline{\Delta_j(T^*1)} \hskip1pt
1_{\Omega_k} \hskip1pt df_{k+s} \big) \Big) \\ & = &
\mathsf{E}_{j-1} \Big( \overline{\Delta_j(T^*1)} \hskip1.5pt
\mathsf{E}_{k+s-1} \hskip1pt (1_{\Omega_k} \hskip1pt df_{k+s})
\Big) \, = \, 0.
\end{eqnarray*}
If $j > k+s$ we have $$\mathsf{E}_{j-1} \Big(
\overline{\Delta_j(T^*1)} \hskip1pt 1_{\Omega_k} \hskip1pt
df_{k+s} \Big) \, = \, \mathsf{E}_{j-1} \big(
\overline{\Delta_j(T^*1)} \big) \hskip1pt 1_{\Omega_k} \hskip1pt
df_{k+s} \, = \, 0.$$ In particular, we obtain the following
identity
\begin{eqnarray*}
1_{\R^n \setminus \Sigma_{f,s}} \summ_k \mathsf{E}_k \Pi
\Delta_{k+s} 1_{\Omega_k} f & = & 1_{\R^n \setminus \Sigma_{f,s}}
\summ_k \mathsf{E}_k \Big( \overline{\Delta_{k+s}(T^*1)} \hskip1pt
1_{\Omega_k} \hskip1pt df_{k+s} \Big) \\ & = & 1_{\R^n \setminus
\Sigma_{f,s}} \summ_k 1_{\Omega_k} \hskip1pt \mathsf{E}_k \Big(
\overline{\Delta_{k+s}(T^*1)} \hskip1pt df_{k+s} \Big) = 0.
\end{eqnarray*}
The last identity follows from the fact that $\Omega_k \subset
\Sigma_{f,s}$ and $\R^n \setminus \Sigma_{f,s}$ are disjoint.

\subsection{Estimating the norm of $\Phi_s$} \label{SSPhi}

Now we estimate the operator norm of the sum $\Phi_s$ under the
assumption that the cancellation condition \eqref{T*1=0} holds for
$T$. We begin by identifying the kernel of the operators appearing
in $\Phi_s$. Let us denote by $k_{\mathsf{e},k}$ and
$k_{\delta,k+s}$ the kernels of $\mathsf{E}_k$ and $\Delta_{k+s}$
respectively. The kernel of the operator $\mathsf{E}_k T
\Delta_{k+s}$ is then given by $$k_{s,k}(x,y) \, = \, \int_{\R^n
\times \R^n}^{\null} k_{\mathsf{e},k}(x,w) \hskip1pt k(w,z)
\hskip1pt k_{\delta, k+s} (z,y) \, dw \, dz.$$ It is
straightforward to verify that
\begin{eqnarray*}
& \displaystyle k_{\mathsf{e},k} (x,w) \, = \, 2^{nk} \sum_{R \in
\Q_k} 1_{R \times R} (x,w), & \\ & \displaystyle k_{\delta,k+s}
(z,y) \, = \, 2^{n(k+s)} \sum_{Q \in \Q_{k+s}} \Big( 1_{Q \times
Q} (z,y) - \frac{1}{2^n} \hskip1pt 1_{Q \times \widehat{Q}} (z,y)
\Big). &
\end{eqnarray*}
Given $x,y \in \R^n$, define $R_x$ to be the only cube in $\Q_k$
containing $x$, while $Q_y$ will stand for the only cube in
$\Q_{k+s}$ containing $y$. Moreover, let $Q_2, Q_3, \ldots,
Q_{2^n}$ be the remaining cubes in $\Q_{k+s}$ sharing dyadic
father with $Q_y$. Let us introduce the following functions
\begin{eqnarray*}
\phi_{R_x}(w) & = & \frac{1}{|R_x|} \, 1_{R_x}(w), \\
\psi_{\widehat{Q}_y}(z) & = & \frac{1}{|\widehat{Q}_y|} \,
\sum_{j=2}^{2^n} 1_{Q_y}(z) - 1_{Q_j}(z).
\end{eqnarray*}
Then the kernel $k_{s,k}(x,y)$ can be written as follows
\begin{equation} \label{kernel}
k_{s,k}(x,y) = \left\langle T(\psi_{\widehat{Q}_y}), \phi_{R_x}
\right\rangle.
\end{equation}
Notice that $\psi_{\widehat{Q}_y} \in H_1$ since it is a linear
combination of atoms.

\subsubsection{Schur type estimates}

In this paragraph we give pointwise estimates for the kernels
$k_{s,k}$ and use them to obtain upper bounds of the Schur
integrals associated to them. Both will be used below to produce
Cotlar type estimates.

\begin{lemma} \label{prelimest}
The following estimates hold\hskip1pt$:$
\begin{itemize}
\item[a)] If $y \in \R^n \setminus 3 R_x$, we have
$$\big| k_{s,k}(x,y) \big| \le \mathrm{c}_n \hskip1pt 2^{-
\gamma (k+s)} \frac{1}{|x-y|^{n+\gamma}}.$$

\item[b)] If $y \in 3 R_x \setminus R_x$, we have $$\big|
k_{s,k}(x,y) \big| \le \mathrm{c}_{n,\gamma} \hskip1pt 2^{- \gamma
(k+s)} 2^{nk} \min \Bigg\{ \int_{R_x} \frac{dw}{|w
-\mathrm{c}_y|^{n+\gamma}}, s 2^{\gamma (k+s)} \Bigg\}.$$

\item[c)] Similarly, if $y \in R_x$ we have $$\hskip14pt \big|
k_{s,k}(x,y) \big| \le \mathrm{c}_{n,\gamma} \hskip1pt 2^{- \gamma
(k+s)} 2^{nk} \min \Bigg\{ \int_{\R^n \setminus R_x} \frac{dw}{|w
-\mathrm{c}_y|^{n+\gamma}}, s 2^{\gamma (k+s)} \Bigg\}.$$
\end{itemize}
The constant $\mathrm{c}_{n,\gamma}$ only depends on $n$ and
$\gamma;$ $\mathrm{c}_y$ denotes the center of the cube
$\widehat{Q}_y$.
\end{lemma}

\dem We proceed in several steps.

\noindent \textbf{The first estimate.} Using $$\int_{\R^n}
\psi_{\widehat{Q}_y}(z) \, dz = 0,$$ we obtain the following
identity where $\mathrm{c}_y$ denotes the center of
$\widehat{Q}_y$
$$\big| k_{s,k}(x,y) \big| \, = \,
\Big| \int_{R_x \times \widehat{Q}_y} \phi_{R_x}(w) \big[ k(w,z) -
k(w,\mathrm{c}_y) \big] \psi_{\widehat{Q}_y}(z) \, dw \, dz
\Big|.$$ Since $|z-\mathrm{c}_y| \le \frac12 |w-\mathrm{c}_y|$ for
$(w,z) \in R_x \times \widehat{Q}_y$, Lipschitz smoothness gives
$$\big| k_{s,k}(x,y) \big| \, \le \, \int_{R_x \times
\widehat{Q}_y} \phi_{R_x}(w)
\frac{|z-\mathrm{c}_y|^\gamma}{|w-\mathrm{c}_y|^{n+\gamma}}
|\psi_{\widehat{Q}_y}(z)| \, dw \, dz.$$ Then, we use
$|w-\mathrm{c}_y| \ge \frac13 |x-y|$ and $|z-\mathrm{c}_y| \le
2^{-(k+s)}$ for $(w,z) \in R_x \times \widehat{Q}_y$
$$\big| k_{s,k}(x,y) \big| \, \le \, \mathrm{c}_n \frac{2^{-\gamma
(k+s)}}{|x-y|^{n+\gamma}} \int_{R_x \times \widehat{Q}_y}
\phi_{R_x}(w) |\psi_{\widehat{Q}_y}(z)| \, dw \, dz \le
\mathrm{c}_n \frac{2^{- \gamma (k+s)}}{|x-y|^{n+\gamma}}.$$

\noindent \textbf{The second estimate.} By \eqref{T*1=0}, we have
$$\int_{\R^n} T(\psi_{\widehat{Q}_y})(w) \, dw = 0.$$
Using this cancellation, we shall use the following relations:
\begin{itemize}
\item If $y \notin R_x \Rightarrow \widehat{Q}_y \not\subset R_x$
and $\displaystyle \big| k_{s,k}(x,y) \big| \, = \,
\frac{1}{|R_x|} \, \Big| \int_{R_x} T(\psi_{\widehat{Q}_y})(w) \,
dw \Big|.$

\vskip5pt

\item If $y \in R_x \Rightarrow \widehat{Q}_y \subset R_x$ and
$\displaystyle \big| k_{s,k}(x,y) \big| \, = \, \frac{1}{|R_x|} \,
\Big| \int_{\R^n \setminus R_x} T(\psi_{\widehat{Q}_y})(w) \, dw
\Big|.$
\end{itemize}
In the first case, we may have
\begin{itemize}
\item[b1)] $3 \widehat{Q}_y \cap R_x = \emptyset$,

\item[b2)] $3 \widehat{Q}_y \cap R_x \neq \emptyset$.
\end{itemize}
If $3 \widehat{Q}_y \cap R_x = \emptyset$, we may use Lipschitz
smoothness as above to obtain
\begin{eqnarray*}
\big| k_{s,k}(x,y) \big| & \le & \frac{1}{|R_x|} \, \int_{R_x
\times \widehat{Q}_y}
\frac{|z-\mathrm{c}_y|^\gamma}{|w-\mathrm{c}_y|^{n+\gamma}}
|\psi_{\widehat{Q}_y}(z)| \, dw \, dz, \\ & \le & \mathrm{c}_n
2^{- \gamma (k+s)} \hskip1pt 2^{nk} \int_{R_x}
\frac{dw}{|w-\mathrm{c}_y|^{n+\gamma}}.
\end{eqnarray*}
On the other hand, if $3 \widehat{Q}_y \cap R_x \neq \emptyset$ we
use the latter estimate on $R_x \setminus 3 \widehat{Q}_y$
\begin{eqnarray*}
\big| k_{s,k}(x,y) \big| & \le & \mathrm{c}_n \hskip1pt 2^{-
\gamma (k+s)} \hskip1pt 2^{nk} \int_{R_x \setminus 3
\widehat{Q}_y} \frac{dw}{|w-\mathrm{c}_y|^{n+\gamma}} \\ & + &
\frac{\mathrm{c}_n}{|R_x| \hskip1pt |\widehat{Q}_y|} \int_{(R_x
\cap \hskip1pt 3 \widehat{Q}_y) \times \widehat{Q}_y} |k(w,z)| \,
dw \, dz.
\end{eqnarray*}
We claim that the second term on the right is dominated by the
first one, up to a constant $\mathrm{c}_{n,\gamma}$ depending only
on $n$ and $\gamma$. Indeed, let us write $\delta_z =
\mathrm{dist} (z, \partial \widehat{Q}_y)$ with $\partial \Omega$
denoting the boundary of $\Omega$. The size estimate for the
kernel gives
$$\frac{1}{|R_x| \hskip1pt |\widehat{Q}_y|} \int_{(R_x \cap
\hskip1pt 3 \widehat{Q}_y) \times \widehat{Q}_y} \big|k(w,z)\big|
\, dw \, dz \le 2^{nk} \hskip1pt 2^{n(k+s)} \int_{\widehat{Q}_y}
\int_{R_x \cap \hskip1pt 3 \widehat{Q}_y} \frac{dw}{|w-z|^n} \, dz
.$$

\noindent
\begin{picture}(360,200)(-180,-100)
\linethickness{1.8pt}
    \put(0,-75){\line(0,1){160}}
    \put(0,-75){\line(1,0){100}}

\linethickness{.8pt}
    \put(0,-15){\line(0,1){30}}
    \put(0,15){\line(-1,0){30}}
    \put(-30,15){\line(0,-1){30}}
    \put(-30,-15){\line(1,0){30}}

    \put(60,-55){\huge $R_x$}

    \put(-28,-10){$\widehat{Q}_y$}
    \put(-52,-35){$3 \widehat{Q}_y$}
    \put(-21.5,3){\mbox{\huge $\cdot$}}

    \qbezier[15](30,-45)(30,0)(30,45)
    \qbezier[15](30,45)(-15,45)(-60,45)
    \qbezier[15](-60,45)(-60,0)(-60,-45)
    \qbezier[15](-60,-45)(-15,-45)(30,-45)

    \put(-19.2,7.4){\circle{14.3}}
    \put(-19.2,7.4){\circle{144.5}}
    \put(-19.2,7.4){\line(2,1){6.5}}
    \put(-19.2,7.4){\line(2,-1){64.5}}
    \put(-21.1,9.7){\tiny $\alpha$}
    \put(-24,3.7){\tiny $z$}
    \put(19.2,-20.4){\tiny $\beta$}
\end{picture}

\null

\vskip-35pt

\null
\begin{center}
\textsc{Figure I} \\ We have $\alpha = \delta_z$ and $\beta \le 2
\sqrt{n} 2^{-(k+s-1)}$
\end{center}

\noindent According to Figure I, we easily see
that
\begin{eqnarray*} \lefteqn{\hskip-10pt 2^{nk} \hskip1pt 2^{n(k+s)}
\int_{\widehat{Q}_y} \int_{R_x \cap \hskip1pt 3 \widehat{Q}_y}
\frac{dw}{|w-z|^n} \, dz} \\ & \le & \mathrm{c}_n 2^{nk} \hskip1pt
2^{n(k+s)} \int_{\widehat{Q}_y} \Big[ \int_{\mathrm{S}_{n-1}}
\Big( \int_{\delta_z}^{2 \sqrt{n} \hskip1pt 2^{-(k+s-1)}}
\frac{dr}{r} \Big) \, d \sigma \Big] \, dz \\ [3pt] & \le &
\mathrm{c}_n 2^{nk} \hskip1pt 2^{n(k+s)} \int_{\widehat{Q}_y} \log
\Big( \frac{2 \sqrt{n} \hskip1pt 2^{-(k+s-1)}}{\delta_z} \Big) \,
dz \cdot \sigma(\mathrm{S}_{n-1}) \\ [1pt] & \sim & \mathrm{c}_n
2^{nk} \hskip1pt 2^{n(k+s)} \int_0^{\sqrt{n} / 2^{k+s}} \log \Big(
\frac{4 \sqrt{n} \hskip1pt 2^{-(k+s)}}{\sqrt{n} \hskip1pt
2^{-(k+s)} -r} \Big) \hskip1pt r^{n-1} \, dr \ \le \ \mathrm{c}_n
2^{nk}.
\end{eqnarray*}

\noindent This gives rise to $$\big|k_{s,k}(x,y)\big| \le
\mathrm{c}_n \hskip1pt 2^{- \gamma (k+s)} \hskip1pt 2^{nk}
\int_{R_x} \frac{dw}{|w-\mathrm{c}_y|^{n+\gamma}} \, + \,
\mathrm{c}_n 2^{nk} \hskip1pt 1_{\mathcal{U}_{s,k}^x}(y),$$ where
the set $\mathcal{U}_{s,k}^x$ is defined by $$\mathcal{U}_{s,k}^x
= \Big\{ y \in \R^n \setminus R_x \, \big| \ \mathrm{dist}(y,
\partial R_x) < 2^{-(k+s-1)} \Big\}.$$ However, it is easily seen
that for $y \in \mathcal{U}_{s,k}^x$ we have
$$2^{- \gamma (k+s)} \int_{R_x}
\frac{dw}{|w-\mathrm{c}_y|^{n+\gamma}} \, \ge \, \mathrm{c}_n 2^{-
\gamma (k+s)} \int_{\mathrm{S}_{n-1}} \Big(
\int_{2^{-(k+s)}}^{2^{-(k+s)} + 2^{-k}} \frac{dr}{r^{1+\gamma}}
\Big) \, d \sigma \ge \mathrm{c}_{n,\gamma}.$$ In particular, we
deduce our claim and so $$\big|k_{s,k}(x,y)\big| \le
\mathrm{c}_{n,\gamma} \hskip1pt 2^{- \gamma (k+s)} \hskip1pt
2^{nk} \int_{R_x} \frac{dw}{|w-\mathrm{c}_y|^{n+\gamma}}.$$ In the
second case $\widehat{Q}_y \subset R_x$, we may have
\begin{itemize}
\item[c1)] $3 \widehat{Q}_y \cap (\R^n \setminus R_x) =
\emptyset$,

\item[c2)] $3 \widehat{Q}_y \cap (\R^n \setminus R_x) \neq
\emptyset$.
\end{itemize}
The argument in this case is entirely similar. Indeed, if the
intersection is empty we use Lipschitz smoothness one more time
and the same argument as above gives
\begin{eqnarray*}
\big|k_{s,k}(x,y)\big| & \le & 2^{- \gamma (k+s)} \hskip1pt 2^{nk}
\int_{\R^n \setminus R_x} \frac{dw}{|w-\mathrm{c}_y|^{n+\gamma}}.
\end{eqnarray*}
If the intersection is not empty, the inequality
$$\frac{1}{|R_x| \hskip1pt |\widehat{Q}_y|} \int_{((\R^n \setminus
R_x) \cap \hskip1pt 3 \widehat{Q}_y) \times \widehat{Q}_y}
|k(w,z)| \, dw \, dz \, \le \, \mathrm{c}_n 2^{nk}$$ can be proved
as above. This gives rise to the estimate
$$\big|k_{s,k}(x,y)\big| \le \mathrm{c}_n \hskip1pt
2^{- \gamma (k+s)} \hskip1pt 2^{nk} \int_{\R^n \setminus R_x}
\frac{dw}{|w-\mathrm{c}_y|^{n+\gamma}} \, + \, \mathrm{c}_n 2^{nk}
\hskip1pt 1_{\mathcal{V}_{s,k}^x}(y),$$ where the set
$\mathcal{V}_{s,k}^x$ is defined by $$\mathcal{V}_{s,k}^x = \Big\{
y \in R_x \, \big| \ \mathrm{dist}(y, \partial R_x) < 2^{-(k+s-1)}
\Big\}.$$ Now we use that for $y \in \mathcal{V}_{s,k}^x$ we have
$$2^{- \gamma (k+s)} \int_{\R^n \setminus R_x}
\frac{dw}{|w-\mathrm{c}_y|^{n+\gamma}} \, \ge \, \mathrm{c}_n 2^{-
\gamma (k+s)} \int_{\mathrm{S}_{n-1}} \Big(
\int_{2^{-(k+s)}}^{\infty} \frac{dr}{r^{1+\gamma}} \Big) \, d
\sigma \ge \mathrm{c}_{n,\gamma}.$$ Our estimates prove the first
halves of inequalities b) and c) above.

\noindent \textbf{The third estimate.} It remains to prove that
$$\big| k_{s,k}(x,y) \big| = \Big| \left\langle
T(\psi_{\widehat{Q}_y}), \phi_{R_x} \right\rangle \Big| \le
\mathrm{c}_n s \hskip1pt 2^{nk}.$$ Since $y \in 3 R_x$, the
localization estimate in Paragraph \ref{3AR} gives
\begin{eqnarray*}
\big| k_{s,k}(x,y) \big| & \le & \mathrm{c}_n \hskip1pt
\ell(\widehat{Q}_y)^n \hskip1pt \log \Big(
\frac{\ell(3R_x)}{\ell(\widehat{Q}_y)} \Big) \hskip1pt
\|\phi_{R_x}\|_\infty \|\psi_{\widehat{Q}_y}\|_\infty
\\ & = & \mathrm{c}_n \hskip1pt |\widehat{Q}_y| \hskip2pt
\log \big( 3 \hskip1pt 2^{s-1} \big) \hskip1pt \frac{1}{|R_x|}
\hskip1pt \frac{2^n-1}{|\widehat{Q}_y|} \ \le \ \mathrm{c}_n
\hskip1pt s \hskip1pt 2^{nk}.
\end{eqnarray*}
We have used that $T$ is assumed to be $L_2$-normalized. The proof
is complete. \fin

\begin{lemma} \label{prelimest2} Let us define
$$\begin{array}{rclcl}
\mathcal{S}^1_{s,k}(x) & = & \displaystyle \int_{\R^n} \big|
k_{s,k}(x,y) \big| \, dy, \\ [10pt] \mathcal{S}^2_{s,k}(y) & = &
\displaystyle \int_{\R^n} \big| k_{s,k}(x,y) \big| \, dx.
\end{array}$$
Then, there exists a constant $\mathrm{c}_{n,\gamma}$ depending
only on $n,\gamma$ such that
$$\begin{array}{rclcl}
\mathcal{S}^1_{s,k}(x) & \le & \displaystyle
\frac{\mathrm{c}_{n,\gamma} \hskip1pt s}{2^{\gamma s}} &
\mbox{for all} \quad (x,k) \in \R^n \!\! \times \Z, \\
[10pt] \mathcal{S}^2_{s,k}(y) & \le & \hskip1pt
\mathrm{c}_{n,\gamma} \hskip1pt s & \mbox{for all} \quad \hskip1pt
(y,k) \in \R^n \!\! \times \Z.
\end{array}$$
\end{lemma}

\dem We estimate $\mathcal{S}^1_{s,k}$ and $\mathcal{S}^2_{s,k}$
in turn.

\noindent \textbf{Estimate of $\mathcal{S}^1_{s,k}(x)$.} Given $x
\in \R^n$, define the cube $R_x$ as above. Then we decompose the
integral defining $\mathcal{S}^1_{s,k}(x)$ into three regions
according to Lemma \ref{prelimest} and estimate each one
independently. Using Lemma \ref{prelimest} a) we find
\begin{equation} \label{part1a}
\int_{\R^n \setminus 3 R_x} \big| k_{s,k}(x,y) \big| \, dy \ \le \
\mathrm{c}_n \hskip1pt 2^{- \gamma (k+s)} \int_{\R^n \setminus 3
R_x} \frac{dy}{|x-y|^{n+\gamma}} \ \le \ \mathrm{c}_n \hskip1pt
2^{- \gamma s}.
\end{equation}
On the other hand, the first estimate in Lemma \ref{prelimest} b)
gives
$$\int_{3R_x \setminus R_x} \big| k_{s,k}(x,y) \big| \, dy \ \le
\ \mathrm{c}_{n,\gamma} \hskip1pt 2^{- \gamma (k+s)} \hskip1pt
2^{nk} \int_{3R_x \setminus R_x} \int_{R_x}
\frac{1}{|w-\mathrm{c}_y|^{n+\gamma}} \, dw \, dy.$$ Now we set
$\delta_w = \mathrm{dist}(w, \partial R_x)$ for $w \in R_x$. Then
we clearly have $$\widehat{\delta}_w \equiv \delta_w + 2^{-(k+s)}
\le \delta_w + \mathrm{dist}(\mathrm{c}_y, \partial R_x) \le
|w-\mathrm{c}_y|.$$

\vskip-15pt

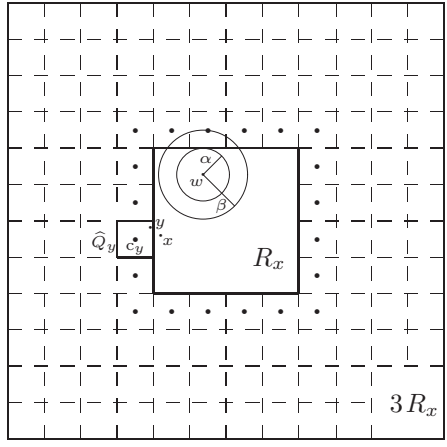
\noindent
\begin{picture}(360,200)(-180,-100)
\linethickness{.6pt}

    \put(-27.5,-27.5){\line(0,1){55}}
    \put(-27.5,-27.5){\line(1,0){55}}
    \put(27.5,27.5){\line(0,-1){55}}
    \put(27.5,27.5){\line(-1,0){55}}

    \put(-82.5,-82.5){\line(0,1){165}}
    \put(-82.5,-82.5){\line(1,0){165}}
    \put(82.5,82.5){\line(0,-1){165}}
    \put(82.5,82.5){\line(-1,0){165}}

    \put(-27.5,0){\line(-1,0){13.75}}
    \put(-41.25,0){\line(0,-1){13.75}}
    \put(-41.25,-13.75){\line(1,0){13.75}}

\linethickness{.1pt}

    \put(62,-72.5){$3 \hskip1pt R_x$}
    \put(10,-17){$R_x$}
    \put(-38,-10.875){\tiny{$\mathrm{c}_y$}}
    \put(-10,15){$\cdot$}
    \put(-27,-2){\tiny{$y$}}
    \put(-30,-5){$\cdot$}
    \put(-8.7,17.5){\circle{20}}
    \put(-8.7,17.5){\circle{33.75}}
    \put(-14,13){\tiny{$w$}}
    \put(-10,22){\tiny{$\alpha$}}
    \put(-4,4){\tiny{$\beta$}}
    \put(-51,-10){\tiny{$\widehat{Q}_y$}}
    \put(-24,-9){\tiny{$x$}}
    \put(-26,-8){$\cdot$}
    \put(-8.7,17.5){\path(0,0)(7.071,7.071)}
    \put(-8.7,17.5){\path(0,0)(11.932,-11.932)}
    \put(-35,-98){$(\alpha,\beta) = (\delta_w, \widehat{\delta}_{w})$}

    \multiput(-82.5,-68.75)(10,0){14}{\line(1,0){5}}
    \multiput(-82.5,-55)(10,0){17}{\line(1,0){5}}
    \multiput(-82.5,-41.25)(10,0){17}{\line(1,0){5}}
    \multiput(-82.5,41.25)(10,0){17}{\line(1,0){5}}
    \multiput(-82.5,55)(10,0){17}{\line(1,0){5}}
    \multiput(-82.5,68.75)(10,0){17}{\line(1,0){5}}

    \multiput(-68.75,-82.5)(0,10){17}{\line(0,1){5}}
    \multiput(-55,-82.5)(0,10){17}{\line(0,1){5}}
    \multiput(-41.25,-82.5)(0,10){7}{\line(0,1){5}}
    \multiput(-41.25,82.5)(0,-10){8}{\line(0,-1){5}}
    \multiput(41.75,-82.5)(0,10){17}{\line(0,1){5}}
    \multiput(55,-82.5)(0,10){17}{\line(0,1){5}}
    \multiput(68.75,82.5)(0,-10){14}{\line(0,-1){5}}

    \multiput(-27.5,-82.5)(0,10){6}{\line(0,1){5}}
    \multiput(-13.75,-82.5)(0,10){6}{\line(0,1){5}}
    \multiput(0,-82.5)(0,10){6}{\line(0,1){5}}
    \multiput(13.75,-82.5)(0,10){6}{\line(0,1){5}}
    \multiput(27.5,-82.5)(0,10){6}{\line(0,1){5}}

    \multiput(-27.5,82.5)(0,-10){6}{\line(0,-1){5}}
    \multiput(-13.75,82.5)(0,-10){6}{\line(0,-1){5}}
    \multiput(0,82.5)(0,-10){6}{\line(0,-1){5}}
    \multiput(13.75,82.5)(0,-10){6}{\line(0,-1){5}}
    \multiput(27.5,82.5)(0,-10){6}{\line(0,-1){5}}

    \multiput(-82.5,-27.5)(10,0){6}{\line(1,0){5}}
    \multiput(-82.5,-13.75)(10,0){4}{\line(1,0){5}}
    \multiput(-82.5,0)(10,0){4}{\line(1,0){5}}
    \multiput(-82.5,13.75)(10,0){6}{\line(1,0){5}}
    \multiput(-82.5,27.5)(10,0){6}{\line(1,0){5}}

    \multiput(82.5,-27.5)(-10,0){6}{\line(-1,0){5}}
    \multiput(82.5,-13.75)(-10,0){6}{\line(-1,0){5}}
    \multiput(82.5,0)(-10,0){6}{\line(-1,0){5}}
    \multiput(82.5,13.75)(-10,0){6}{\line(-1,0){5}}
    \multiput(82.5,27.5)(-10,0){6}{\line(-1,0){5}}

    \multiput(-36.675,-38.775)(13.75,0){6}{\mbox{\huge $\cdot$}}
    \multiput(-36.675,29.975)(13.75,0){6}{\mbox{\huge $\cdot$}}
    \multiput(-36.675,16.225)(0,-13.75){4}{\mbox{\huge $\cdot$}}
    \multiput(32.075,16.225)(0,-13.75){4}{\mbox{\huge $\cdot$}}
\end{picture}

\begin{center} \textsc{Figure II} \\ Even if $x$ and $y$ are
close, we have a $(\widehat{\delta}_w - \delta_w)$ \hskip-1pt --
\hskip-1pt margin
\end{center}

\vskip5pt

\noindent In particular, we find (see Figure II above)
$$\int_{3R_x \setminus R_x} \frac{dy}{|w-\mathrm{c}_y|^{n+\gamma}}
\, \lesssim \,
\frac{|\mathsf{B}_{\widehat{\delta}_w}(w)|}{\widehat{\delta}_w^{n+\gamma}}
+ \int_{\R^n \setminus \mathsf{B}_{\widehat{\delta}_w}(w)}
\frac{dy}{|w-y|^{n+\gamma}} \, \sim \,
1/\widehat{\delta}_w^\gamma.$$ This provides us with the estimate
\begin{eqnarray*}
\int_{3R_x \setminus R_x} \int_{R_x}
\frac{1}{|w-\mathrm{c}_y|^{n+\gamma}} \, dw \, dy \, & \le & \,
\mathrm{c}_n \int_{\mathrm{S}_{n-1}} \int_0^{2^{-k}}
\frac{r^{n-1}}{(2^{-(k+s)} + 2^{-k} - r)^\gamma} \, dr \, d
\sigma.
\end{eqnarray*}
Using $t = 2^{-k} + 2^{-(k+s)} - r$ and the bound $r \le 2^{-k}$
\begin{eqnarray*}
\int_{3R_x \setminus R_x} \int_{R_x}
\frac{1}{|w-\mathrm{c}_y|^{n+\gamma}} \, dw \, dy & \le &
\mathrm{c}_n 2^{-(n-1)k} \int_{2^{-(k+s)}}^{2^{-k}}
\frac{dt}{t^\gamma} \
\\ & \le & \mathrm{c}_n \begin{cases} s \hskip1pt 2^{-(n-1)k} & \mathrm{if} \
\gamma = 1, \\ \mathrm{c}_\gamma \hskip1pt 2^{-nk} \hskip1pt
2^{\gamma k} & \mathrm{if} \ 0 < \gamma < 1.
\end{cases}
\end{eqnarray*}
In summary, combining our estimates we have obtained
\begin{equation} \label{part1b}
\int_{3R_x \setminus R_x} \big| k_{s,k}(x,y) \big| \, dy \ \le \
\mathrm{c}_{n,\gamma} \hskip1pt s \hskip1pt 2^{- \gamma s}.
\end{equation}
It remains to control the integral over $R_x$. By Lemma
\ref{prelimest} c) $$\int_{R_x} \big| k_{s,k}(x,y) \big| \, dy \
\le \ \ \mathrm{c}_{n,\gamma} \hskip1pt 2^{- \gamma (k+s)}
\hskip1pt 2^{nk} \int_{R_x} \int_{\R^n \setminus R_x}
\frac{1}{|w-\mathrm{c}_y|^{n+\gamma}} \, dw \, dy.$$ For any given
$y \in R_x$, we set again $$\delta_{\mathrm{c}_y} =
\mathrm{dist}(\mathrm{c}_y, \partial R_x) \ge 2^{-(k+s)}.$$
Arguing as above, we may use polar coordinates to obtain
\begin{eqnarray*}
\int_{R_x} \int_{\R^n \setminus R_x}
\frac{1}{|w-\mathrm{c}_y|^{n+\gamma}} \, dw \, dy & \le &
\int_{R_x} \Big( \int_{\mathrm{S}_{n-1}}
\int_{\delta_{\mathrm{c}_y}}^{\infty} \frac{r^{n-1}}{r^{n+\gamma}}
\, dr \, d \sigma \Big) \, dy \ \sim \ \int_{R_x}
\frac{dy}{\delta_{\mathrm{c}_y}^\gamma} \\ & \sim &
\int_{\mathrm{S}_{n-1}} \int_0^{2^{-k}-2^{-(k+s)}}
\frac{r^{n-1}}{(2^{-k} - r)^\gamma} \, dr \, d \sigma \\ & + &
\int_{\mathrm{S}_{n-1}} \int_{2^{-k} - 2^{-(k+s)}}^{2^{-k}}
\frac{r^{n-1}}{2^{- \gamma (k+s)}} \, dr \, d \sigma.
\end{eqnarray*}
The first integral is estimated as above
$$\int_{\mathrm{S}_{n-1}} \int_0^{2^{-k}-2^{-(k+s)}}
\frac{r^{n-1}}{(2^{-k} - r)^\gamma} \, dr \, d \sigma \ \le \
\mathrm{c}_n
\begin{cases}
s \hskip1pt 2^{-(n-1)k} & \mathrm{if} \ \gamma = 1, \\
\mathrm{c}_\gamma \hskip1pt 2^{-nk} \hskip1pt 2^{\gamma k} &
\mathrm{if} \ 0 < \gamma < 1, \end{cases}$$ as for the second we
obtain an even better bound. Indeed, we have
\begin{eqnarray*}
\int_{\mathrm{S}_{n-1}} \int_{2^{-k} - 2^{-(k+s)}}^{2^{-k}}
\frac{r^{n-1}}{2^{- \gamma (k+s)}} \, dr \, d \sigma & \sim &
2^{\gamma (k+s)} \hskip1pt \Big( 2^{-nk} - \big[ 2^{-k} -
2^{-(k+s)} \big]^n \Big) \\ & = & 2^{\gamma (k+s)} \hskip1pt
2^{-nk} \Big( 1 - \big[ 1 - 2^{-s} \big]^n \Big) \\ & \le &
2^{\gamma (k+s)} \hskip1pt 2^{-nk} \sum_{j=1}^n {{n}\choose{j}}
2^{-sj} \\ & \le & \mathrm{c}_n \hskip1pt 2^{-nk} \hskip1pt
2^{\gamma k}.
\end{eqnarray*}
Writing all together we finally get
\begin{equation} \label{part1c}
\int_{R_x} \big| k_{s,k}(x,y) \big| \, dy \ \le \
\mathrm{c}_{n,\gamma} \hskip1pt s \hskip1pt 2^{- \gamma s}.
\end{equation}
According to \eqref{part1a}, \eqref{part1b} and \eqref{part1c} we
obtain the upper bound $\mathcal{S}^1_{s,k}(x) \le
\mathrm{c}_{n,\gamma} \hskip1pt s \hskip1pt 2^{- \gamma s}$.

\noindent \textbf{Estimate of $\mathcal{S}^2_{s,k}(y)$.} Given a
fixed point $y$, we consider a partition $\R^n = \Omega_1 \cup
\Omega_2$ where $\Omega_1$ is the set of points $x$ such that $y
\notin 3 R_x$ and $\Omega_2 = \R^n \setminus \Omega_1$. In the
region $\Omega_1$ we may proceed as in \eqref{part1a}. On the
other hand, inside $\Omega_2$ and according to Lemma
\ref{prelimest} we know that $|k_{s,k}(x,y)| \le
\mathrm{c}_{n,\gamma} \hskip1pt s \hskip1pt 2^{nk}$. This means
that we have
$$\mathcal{S}^2_{s,k}(y) \le \mathrm{c}_{n,\gamma} \Big( 2^{-
\gamma s} + |\Omega_2| \hskip1pt s \hskip1pt 2^{nk} \Big) =
\mathrm{c}_{n,\gamma} \Big( 2^{- \gamma s} + |3R_y| \hskip1pt s
\hskip1pt 2^{nk} \Big) \le \mathrm{c}_{n,\gamma} \hskip1pt s.$$
This upper bound holds for all $(y,k) \in \R^n \!\! \times \Z$.
Hence, the proof is complete. \fin

\subsubsection{Cotlar type estimates}

Let us write $\Lambda_{s,k}$ for $\mathsf{E}_k T \Delta_{k+s}$.
According to the pairwise orthogonality of martingale differences,
we have $\Lambda_{s,i}^{\null} \Lambda_{s,j}^* = 0$ whenever $i
\neq j$. In particular, it follows from Cotlar lemma that it
suffices to control the norm of the operators $\Lambda_{s,i}^*
\Lambda_{s,j}^{\null}$. Explicitly, our estimate for $\Phi_s$
stated in the shifted $T1$ theorem will be deduced from $$\big\|
\Lambda_{s,i}^* \Lambda_{s,j}^{\null} \big\|_{\mathcal{B}(L_2)} \,
\le \, \mathrm{c}_{n,\gamma} \hskip1pt s^2 \hskip1pt 2^{- \gamma
s/2} \alpha_{i-j}^2$$ for some summable sequence $(\alpha_k)_{k
\in \Z}$. The kernel of $\Lambda_{s,i}^* \Lambda_{s,j}^{\null}$ is
given by
$$k_{i,j}^s(x,y) \, = \, \int_{\R^n} \overline{k_{s,i}(z,x)}
\hskip1pt k_{s,j}(z,y) \, dz.$$ Before proceeding with our
estimates we need to point out another cancellation property which
easily follows from \eqref{T*1=0}. Given $r > 0$ and a point $y
\in \R^n$, let $f(z) = 1_{\mathsf{B}_r(y)}(z) /
|\mathsf{B}_r(y)|$. Then it is clear that $\Delta_{k+s} f =
df_{k+s}$ is in $H_1$ since it can be written as a linear
combination of atoms. According to our cancellation condition
\eqref{T*1=0} we find
$$\int_{\R^n} \mathsf{E}_k T \Delta_{k+s} f(x) \, dx = \int_{\R^n}
T df_{k+s} (x) \, dx = 0.$$ In terms of the kernels, this identity
is written as
$$\int_{\R^n} \Big( \int_{\R^n} k_{s,k}(x,z) f(z) \, dz \Big) \,
dx = 0.$$ Using Fubini theorem (our estimates in Lemma
\ref{prelimest} ensure the integrability) and taking the limit as
$r \to 0$, the Lebesgue differentiation theorem implies the
following identity, which holds for almost every point $y$
\begin{equation} \label{cancellation2}
\int_{\R^n} k_{s,k}(x,y) \, dx = 0.
\end{equation}
This holds for all $k \in \Z$ and we deduce
\begin{eqnarray*}
k_{i,j}^s(x,y) & = & \int_{\R^n} \overline{k_{s,i}(z,x)} \hskip1pt
\Big( k_{s,j}(z,y) - k_{s,j}(x,y) \Big) \, dz \\ & = & \int_{\R^n}
\! \Big( \overline{k_{s,i}(z,x)} - \overline{k_{s,i}(y,x)} \Big)
\hskip1pt k_{s,j}(z,y) \ dz.
\end{eqnarray*}
In order to estimate the kernels $k_{i,j}^s$, we use the first or
the second expression above according to whether $i \ge j$ or not.
Since the estimates are entirely similar we shall assume in what
follows that $i \ge j$ and work in the sequel with the first
expression above. Moreover, given $w \in \R^n$ we shall write all
through out this paragraph $R_w$ for the only cube in $\Q_j$
containing $w$. Then, since $R_z = R_x$ whenever $z \in R_x$, it
follows from \eqref{kernel} that
\begin{eqnarray*}
k_{i,j}^s(x,y) & = & \int_{\R^n \setminus 3 R_x}
\overline{k_{s,i}(z,x)} \hskip1pt \Big( k_{s,j}(z,y) -
k_{s,j}(x,y) \Big) \, dz \\ & + & \int_{3R_x \setminus R_x}
\overline{k_{s,i}(z,x)} \hskip1pt \Big( k_{s,j}(z,y) -
k_{s,j}(x,y) \Big) \, dz.
\end{eqnarray*}
If $\alpha_{i,j}^s(x,y)$ and $\beta_{i,j}^s(x,y)$ are the first
and second terms above, let
\begin{eqnarray*}
\mathcal{S}_{i,j,s}^{1,\alpha}(x) & = & \int_{\R^n} \big|
\alpha_{i,j}^s(x,y) \big| \, dy, \\
\mathcal{S}_{i,j,s}^{2,\alpha}(y) & = & \int_{\R^n} \big|
\alpha_{i,j}^s(x,y) \big| \, dx, \\
\mathcal{S}_{i,j,s}^{1,\beta}(x) & = & \int_{\R^n} \big|
\beta_{i,j}^s(x,y) \big| \, dy, \\
\mathcal{S}_{i,j,s}^{2,\beta}(y) & = & \int_{\R^n} \big|
\beta_{i,j}^s(x,y) \big| \, dx.
\end{eqnarray*}
According to Schur lemma from Paragraph \ref{3AR}, we obtain the
upper bound
\begin{equation} \label{schurbound}
\big\| \Lambda_{s,i}^* \Lambda_{s,j}^{\null}
\big\|_{\mathcal{B}(L_2)} \, \le \, \sqrt{ \Big( \big\|
\mathcal{S}_{i,j,s}^{1,\alpha} \big\|_\infty + \big\|
\mathcal{S}_{i,j,s}^{1,\beta} \big\|_\infty \Big) \Big( \big\|
\mathcal{S}_{i,j,s}^{2,\alpha} \big\|_\infty + \big\|
\mathcal{S}_{i,j,s}^{2,\beta} \big\|_\infty \Big)}.
\end{equation}

\begin{lemma} \label{cotlarest1}
We have $$\max \Big\{ \big\| \mathcal{S}_{i,j,s}^{1,\alpha}
\big\|_\infty, \big\| \mathcal{S}_{i,j,s}^{1,\beta} \big\|_\infty
\Big\} \, \le \, \mathrm{c}_{n,\gamma} \hskip1pt s \hskip1pt (s +
|i-j|) \hskip1pt 2^{- \gamma s}.$$
\end{lemma}

\dem According to Lemma \ref{prelimest}, we know that $$\big|
k_{s,i}(z,x) \big| \le \mathrm{c}_n \hskip1pt 2^{- \gamma (i+s)}
\frac{1}{|x-z|^{n+\gamma}}$$ whenever $z \notin 3 \hskip1pt R_x$.
\hskip1pt Moreover, Lemma \ref{prelimest2} gives $$\int_{\R^n}
\big| k_{s,j}(z,y) - k_{s,j}(x,y)\big| \, dy \, \le \,
\mathrm{c}_{n,\gamma} \hskip1pt s \hskip1pt 2^{- \gamma s}.$$ If
we combine the two estimates above, we obtain
\begin{eqnarray*}
\mathcal{S}_{i,j,s}^{1,\alpha}(x) & \le & \int_{\R^n} \Big(
\int_{\R^n \setminus 3 R_x} \big| k_{s,i}(z,x) \big| \,
\big| k_{s,j}(z,y) - k_{s,j}(x,y)\big| \, dz \Big) \, dy \\
& = & \int_{\R^n \setminus 3 R_x} \big| k_{s,i}(z,x) \big| \,
\Big( \int_{\R^n} \big| k_{s,j}(z,y) - k_{s,j}(x,y)\big| \, dy
\Big) \, dz \\ & \le & \mathrm{c}_{n,\gamma} \hskip1pt s \hskip1pt
2^{-\gamma s} 2^{-\gamma (i+s)} \int_{\R^n \setminus 3 R_x}
\frac{dz}{|x-z|^{n+\gamma}} \ \le \ \mathrm{c}_{n,\gamma}
\hskip1pt s \hskip1pt 2^{- 2 \gamma s} 2^{- \gamma |i-j|}.
\end{eqnarray*}
The last inequality uses the assumption $i \ge j$, so that $i-j =
|i-j|$. The estimate above holds for all $x \in \R^n$. Hence,
since $\mathrm{c}_{n,\gamma} s \hskip1pt 2^{- 2 \gamma s} 2^{-
\gamma |i-j|}$ is much smaller than $\mathrm{c}_{n,\gamma} s
\hskip1pt (s + |i-j|) \hskip1pt 2^{- \gamma s}$, it is clear that
the first function satisfies the thesis. Let us now proceed with
the second function. To that aim we observe that $k_{s,j}(z,y)$ is
$j$-measurable as a function in $z$, meaning that $\mathsf{E}_j
(k_{s,j}(\, \cdot,y))(z) = k_{s,j}(z,y)$. This follows from
\eqref{kernel}. In particular, the same holds for the function
$$1_{3 R_x \setminus R_x} (z) \Big(
k_{s,j}(z,y) - k_{s,j}(x,y) \Big).$$ Therefore, using the integral
invariance of conditional expectations
\begin{eqnarray*}
\beta_{i,j}^s(x,y) & = & \int_{3R_x \setminus R_x}
\overline{\mathsf{E}_j(k_{s,i}(\, \cdot ,x))(z)} \hskip1pt \Big(
k_{s,j}(z,y) - k_{s,j}(x,y) \Big) \, dz \\ & = & \sum_{R \sim R_x}
\int_R \overline{\frac{1}{|R|} \int_R^{\null} k_{s,i}(w,x) \, dw}
\hskip1pt \Big( k_{s,j}(z,y) - k_{s,j}(x,y) \Big) \, dz,
\end{eqnarray*}
where $R \sim R_x$ is used to denote that $R$ is a neighbor of
$R_x$ in $\Q_j$. That is, the neighbors of $R_x$ form a partition
of $3 \hskip1pt R_x \setminus R_x$ formed by $3^n - 1$ cubes in
$\Q_j$. If $\mathrm{c}_R$ denotes the center of $R$, we use that
$k_{s,j}(z,y) = k_{s,j}(\mathrm{c}_R,y)$ for $z \in R$ and obtain
the estimate
\begin{equation} \label{estbeta}
\big| \beta_{i,j}^s(x,y) \big| \, \le \, \sum_{R \sim R_x} \Big|
\int_R k_{s,i} (w,x) \, dw \Big| \, \big| k_{s,j}(\mathrm{c}_R,y)
- k_{s,j}(x,y) \big|.
\end{equation}
This, combined with Lemma \ref{prelimest2}, produces
\begin{equation} \label{s1beta}
\mathcal{S}_{i,j,s}^{1,\beta}(x) \, \le \,  \mathrm{c}_{n,\gamma}
\hskip1pt s \hskip1pt 2^{-\gamma s} \sum_{R \sim R_x} \Big| \int_R
k_{s,i} (w,x) \, dw \Big|.
\end{equation}
Let us now estimate the integral. If $w \in S_w \in \Q_i$ and $x
\in O_x \in \Q_{i+s}$
\begin{eqnarray*}
\int_R k_{s,i} (w,x) \, dw & = & \int_R \Big\langle T
\psi_{\widehat{O}_x}, \phi_{S_w} \Big\rangle \, dw \\ & = &
\sum_{S \subset R, S \in \Q_i} \int_S T \psi_{\widehat{O}_x} (z)
\, dz \ = \ \Big\langle T \psi_{\widehat{O}_x}, 1_R \Big\rangle.
\end{eqnarray*}
Now we use the localization estimate from Paragraph \ref{3AR} to
obtain
$$\Big| \Big\langle T \psi_{\widehat{O}_x}, 1_R \Big\rangle \Big|
\, \le \, \mathrm{c}_n \hskip1pt \ell(\widehat{O}_x)^n \hskip1pt
\log \Big( \frac{\ell(3R)}{\ell(\widehat{O}_x)} \Big) \hskip1pt
\|\psi_{\widehat{O}_x}\|_\infty \|1_R\|_\infty \, \le \,
\mathrm{c}_n \big( s + |i-j| \big).$$ Since there are $3^n - 1$
neighbors, this estimate completes the proof with \eqref{s1beta}.
\fin

\begin{lemma} \label{cotlarest2}
We have $$\max \Big\{ \big\| \mathcal{S}_{i,j,s}^{2,\alpha}
\big\|_\infty, \big\| \mathcal{S}_{i,j,s}^{2,\beta} \big\|_\infty
\Big\} \, \le \, \mathrm{c}_n s^2 \hskip1pt \big( 1 + |i-j| \big)
\hskip1pt 2^{- \gamma |i-j|}.$$
\end{lemma}

\dem Once again, Lemma \ref{prelimest} gives
$$\big| k_{s,i}(z,x) \big| \le \mathrm{c}_n \hskip1pt 2^{- \gamma
(i+s)} \frac{1}{|x-z|^{n+\gamma}} \qquad \mbox{for} \qquad z
\notin 3 \hskip1pt R_x.$$ This, together with Fubini theorem
produces
\begin{eqnarray*}
\mathcal{S}_{i,j,s}^{2,\alpha}(y) & \le & \mathrm{c}_n \int_{\R^n}
\Big( \int_{\R^n \setminus 3 R_x} \frac{2^{- \gamma
(i+s)}}{|x-z|^{n+\gamma}} \,
\big| k_{s,j}(z,y) - k_{s,j}(x,y)\big| \, dz \Big) \, dx \\
& \le & \mathrm{c}_n \int_{\R^n} \Big( \int_{\R^n \setminus
\mathsf{B}_{2^{-j}}(z)} \frac{2^{- \gamma
(i+s)}}{|x-z|^{n+\gamma}} \, dx \Big) \, \big| k_{s,j}(z,y) \big|
\, dz \\ & + & \mathrm{c}_n \int_{\R^n} \Big( \int_{\R^n \setminus
\mathsf{B}_{2^{-j}}(x)} \frac{2^{- \gamma
(i+s)}}{|x-z|^{n+\gamma}} \, dz \Big) \, \big| k_{s,j}(x,y) \big|
\, dx \\ & = & \mathrm{c}_n \hskip1pt 2^{- \gamma s} 2^{- \gamma
|i-j|} \int_{\R^n} \big| k_{s,j}(w,y) \big| \, dw.
\end{eqnarray*}
Now, according to Lemma \ref{prelimest2} we now that the integral
on the right is bounded by $\mathrm{c}_{n,\gamma} s$ for all $y$
in $\R^n$. Therefore, the $L_\infty$ norm of the first function is
much smaller than our upper bound. Let us now estimate the second
function. If we proceed as in Lemma \ref{cotlarest1} and use
\eqref{estbeta}, we find
$$\mathcal{S}_{i,j,s}^{2,\beta}(y) \, \le \, \int_{\R^n}
\sum_{R \sim R_x} \Big| \int_R k_{s,i} (w,x) \, dw \Big| \, \big|
k_{s,j}(\mathrm{c}_R,y) - k_{s,j}(x,y) \big| \, dx.$$ Now we need
a different estimate for the integral of $k_{s,i}( \, \cdot, x)$
over the neighbor cubes $R$ of $R_x$. Indeed, combining the
pointwise estimates obtained in Lemma \ref{prelimest} it easily
follows that
\begin{equation} \label{globalest}
\big| k_{s,i}(w,x) \big| \, \le \, \mathrm{c}_{n,\gamma} \,
\frac{s \hskip1pt 2^{ni}}{\big( 1 + 2^i |x-w| \big)^{n+\gamma}}
\qquad \mbox{for all} \qquad (w,x) \in \R^n \times \R^n.
\end{equation}
If we set $\delta_x = \mathrm{dist}(x, \partial R_x) \le
\mathrm{dist}(x,
\partial R)$, we get
\begin{eqnarray*}
\Big| \int_R k_{s,i} (w,x) \, dw \Big| & \le &
\mathrm{c}_{n,\gamma} \hskip1pt s \int_R \frac{2^{ni}}{\big( 1 +
2^i |x-w| \big)^{n+\gamma}} \, dw \\ & \le & \mathrm{c}_{n,\gamma}
\hskip1pt s \int_{\mathrm{S}_{n-1}} \Big( \int_{\delta_x}^{\infty}
\frac{2^{ni} r^{n-1}}{\big( 1 + 2^i r \big)^{n+\gamma}} \, dr
\Big) \, d \sigma \\ & = & \mathrm{c}_{n,\gamma} \hskip1pt s
\int_{2^i \delta_x}^{\infty} \frac{z^{n-1}}{\big( 1 + z
\big)^{n+\gamma}} \, dz \ \le \ \mathrm{c}_{n,\gamma} \hskip1pt s
\, \frac{1}{\big( 1 + 2^i \delta_x \big)^\gamma}.
\end{eqnarray*}
Using \eqref{globalest} for $k_{s,j}$, we have
$$\mathcal{S}_{i,j,s}^{2,\beta}(y) \le \mathrm{c}_{n,\gamma}
\hskip1pt s^2 \hskip1pt 2^{nj} \hskip1pt \Upsilon(i,j,\gamma)$$
where the term $\Upsilon(i,j,\gamma)$ is given by
$$\int_{\R^n} \sum_{R \sim R_x} \frac{1}{\big( 1 + 2^i
\delta_x \big)^\gamma} \, \Bigg( \frac{1}{\big( 1 + 2^j
|\mathrm{c}_R - y| \big)^{n+\gamma}} + \frac{1}{\big( 1 + 2^j |x -
y| \big)^{n+\gamma}} \Bigg) \, dx.$$ It is straightforward to see
that it suffices to estimate the integral
\begin{equation} \label{simplificationest}
\int_{\R^n} \frac{1}{\big( 1 + 2^i \delta_x \big)^\gamma} \,
\frac{1}{\big( 1 + 2^j |x - y| \big)^{n+\gamma}} \, dx.
\end{equation}
Indeed, both functions inside the big bracket above are comparable
and the sum $\sum_{R \sim R_x}$ can be deleted since it only
provides an extra factor of $3^n -1$. Now, the main idea to
estimate \eqref{simplificationest} is to observe that the two
functions in the integrand are nearly independent inside any
dyadic cube of $\Q_j$. Let us be more explicit, we have
\begin{eqnarray*}
\lefteqn{\hskip-15pt \int_{\R^n} \frac{1}{\big( 1 + 2^i \delta_x
\big)^\gamma} \, \frac{1}{\big( 1 + 2^j |x - y| \big)^{n+\gamma}}
\, dx} \\ & \le & \sum_{R \in \Q_j} \int_R \frac{1}{\big( 1 + 2^i
\delta_x \big)^{\gamma}} \, dx \ \ \frac{1}{\big( 1 + 2^j
\hskip1pt \mathrm{dist}(R, R_y) \big)^{n+\gamma}} \\ & \le &
\sup_{R \in \Q_j} \Big( \int_R \frac{1}{\big( 1 + 2^i \delta_x
\big)^{\gamma}} \, dx \Big) \sum_{R \in \Q_j} \ \ \frac{1}{\big( 1
+ 2^j \hskip1pt \mathrm{dist}(R, R_y) \big)^{n+\gamma}} \\ & \sim
& \sup_{R \in \Q_j} \Big( \int_R \frac{1}{\big( 1 + 2^i \delta_x
\big)^{\gamma}} \, dx \Big) \ \int_{\R^n} \frac{2^{nj}}{\big( 1 +
2^j \hskip1pt |x - y| \big)^{n+\gamma}} \, dx.
\end{eqnarray*}
The integral on the right is majorized by an absolute constant.
Moreover, recalling that $\delta_x$ stands for $\mathrm{dist}(x,
\partial R_x)$ and that $R_x = R$ for any $x \in R \in \Q_j$, it
is clear that the integral on the left does not depend on the
chosen cube $R$, so that the supremum is unnecessary. To estimate
this integral we set $\lambda = 1 - 2^{j-i}$
\begin{eqnarray*}
\int_{R} \frac{dx}{\big( 1 + 2^i \delta_x \big)^\gamma} & \sim &
\int_{\mathrm{S}_{n-1}} \Big( \int_0^{2^{-j}} \frac{r^{n-1}}{\big(
1 + 2^i (2^{-j} - r) \big)^\gamma} \, dr \Big) \, d \sigma \\ &
\sim & \int_0^{\lambda 2^{-j}} \frac{r^{n-1}}{\big( 1 + 2^i
(2^{-j} - r) \big)^\gamma} \, dr + \int_{\lambda 2^{-j}}^{2^{-j}}
\frac{r^{n-1}}{\big( 1 + 2^i (2^{-j} - r) \big)^\gamma} \, dr.
\end{eqnarray*}
The first integral is majorized by
\begin{eqnarray*}
2^{- \gamma i} \int_{0}^{\lambda 2^{-j}} \frac{r^{n-1}}{(2^{-j} -
r)^\gamma} \, dr & \le & 2^{- \gamma i} \lambda^{n-1} 2^{-nj} 2^j
\int_{0}^{\lambda 2^{-j}} \frac{dr}{(2^{-j} - r)^\gamma} \\ & \le
& 2^{-nj}
\begin{cases}
|i-j| 2^{-|i-j|} & \mbox{if} \ \gamma = 1,\\
\mathrm{c}_\gamma 2^{-\gamma|i-j|} & \mbox{if} \ 0 < \gamma < 1.
\end{cases}
\end{eqnarray*}
The second integral is majorized by
$$\int_{\lambda 2^{-j}}^{2^{-j}} r^{n-1} \,
dr \le 2^{- nj} 2^j (2^{-j} - \lambda 2^{-j}) = 2^{-nj}
2^{-|i-j|}.$$ Combining our estimates we finally get
$$\mathcal{S}_{i,j,s}^{2,\beta}(y) \le
\mathrm{c}_{n,\gamma} \hskip1pt s^2 \hskip1pt \big( 1 + |i-j|
\big) \hskip1pt 2^{- \gamma |i-j|}.$$ Since the last estimate
holds for all $y \in \R^n$, the proof is complete. \fin

\noindent \textbf{Conclusion.} According to \eqref{schurbound},
Lemmas \ref{cotlarest1} and \ref{cotlarest2} give
$$\big\| \Lambda_{s,i}^* \Lambda_{s,j}^{\null}
\big\|_{\mathcal{B}(L_2)} \le \mathrm{c}_{n,\gamma} \hskip1pt
\sqrt{ s^3 \hskip1pt (s + |i-j|)^2 \hskip1pt 2^{- \gamma s}
\hskip1pt 2^{- \gamma |i-j|}} \le \mathrm{c}_{n,\gamma} \hskip1pt
s^2 \hskip1pt 2^{- \gamma s/2} \, \alpha_{i-j}^2,$$ where
$\alpha_k = (1+|k|)^{\frac12} \hskip1pt 2^{- \gamma |k|/4}$. In
particular, Cotlar lemma provides the estimate
$$\|\Phi_s\|_{\mathcal{B}(L_2)} \, = \, \Big\| \summ_k
\mathsf{E}_k T \Delta_{k+s} \Big\|_{\mathcal{B}(L_2)} \, \le \,
\mathrm{c}_{n,\gamma} \hskip1pt s \hskip1pt 2^{- \gamma s/4}
\summ_k \alpha_k \, = \, \mathrm{c}_{n,\gamma} \hskip1pt s
\hskip1pt 2^{- \gamma s/4}.$$

\subsection{Estimating the norm of $\Psi_s$} \label{SSPsi}

We finally estimate the operator norm of $\Psi_s$. This will
complete the proof of our pseudo-localization principle. We shall
adapt some of the notation introduced in the previous paragraph.
Namely, we shall now write $\Lambda_{s,k}$ when referring to the
operator $(id - \mathsf{E}_k) T_{4 \cdot 2^{-k}} \Delta_{k+s}$ and
$k_{s,k}(x,y)$ will be reserved for its kernel. Arguing as above
it is simple to check that we have
$$k_{s,k}(x,y) = T_{4 \cdot 2^{-k}} \psi_{\widehat{Q}_y} (x) -
\Big\langle T_{4 \cdot 2^{-k}} \psi_{\widehat{Q}_y}, \phi_{R_x}
\Big\rangle.$$ We shall use the terminology
\begin{eqnarray*}
k_{s,k}^1(x,y) & = & T_{4 \cdot 2^{-k}} \psi_{\widehat{Q}_y} (x),
\\ k_{s,k}^2(x,y) & = & \Big\langle T_{4 \cdot 2^{-k}}
\psi_{\widehat{Q}_y}, \phi_{R_x} \Big\rangle.
\end{eqnarray*}

\subsubsection{Schur type estimates}

\begin{lemma} \label{pointwisetruncated}
Let us consider the sets $$\mathcal{W}_{s,k}^x \, = \, \Big\{ w
\in \R^n \, \big| \ 4 \cdot 2^{-k} - 2^{-(k+s-1)} \le |x-w| < 4
\cdot 2^{-k} + 2^{-(k+s-1)} \Big\}.$$ Then, the following
pointwise estimate holds
$$\big| k_{s,k}(x,y) \big| \, \le \, \mathrm{c}_n 1_{\R^n
\setminus \mathsf{B}_{2 \cdot 2^{-k}}(x)}(y) \hskip1pt \Big(
\frac{2^{- \gamma (k+s)}}{|x-y|^{n+\gamma}} + 2^{nk} \hskip1pt
1_{\mathcal{W}_{s,k}^x}(y) \Big).$$
\end{lemma}

\dem We have $$k_{s,k}^1(x,y) \, = \, \int_{\widehat{Q}_y} 1_{\R^n
\setminus \mathsf{B}_{4 \cdot 2^{-k}}(x)}(z) \hskip1pt k(x,z)
\hskip1pt \psi_{\widehat{Q}_y}(z) \, dz.$$ If $|x-y| \le 3 \cdot
2^{-k}$ we have $$|x-z| \le |x-y| + |y-z| \le |x-y| + 2^{-(k+s-1)}
\le 4 \cdot 2^{-k}$$ since $z \in \widehat{Q}_y$. In particular,
we obtain $$k_{s,k}^1(x,y) = 0 \qquad \mbox{whenever} \qquad |x-y|
\le 3 \cdot 2^{-k}.$$ If $|x-y| > 5 \cdot 2^{-k}$, then we have
for $z \in \widehat{Q}_y$
$$|x-z| \ge |x-y| - |z-y| \ge |x-y| - 2^{-(k+s-1)}
> 4 \cdot 2^{-k}.$$ Thus, we can argue in the usual way and obtain
\begin{eqnarray*}
\big| k_{s,k}^1(x,y) \big| & = & \Big| \int_{\widehat{Q}_y} \Big(
k(x,z) - k(x,\mathrm{c}_y) \Big) \hskip1pt \psi_{\widehat{Q}_y}(z)
\, dz \Big| \\ & \le & \mathrm{c}_n \, \frac{2^{- \gamma
(k+s)}}{|x-y|^{n+\gamma}} \, \int_{\widehat{Q}_y} \big|
\psi_{\widehat{Q}_y} (z) \big| \, dz \ \le \ \mathrm{c}_n \,
\frac{2^{- \gamma (k+s)}}{|x-y|^{n+\gamma}}.
\end{eqnarray*}
If $3 \cdot 2^{-k} < |x-y| \le 5 \cdot 2^{-k}$, we write
$k_{s,k}^1(x,y)$ as a sum of two integrals
\begin{eqnarray*}
k_{s,k}^1(x,y) & = & \int_{\R^n} 1_{\R^n \setminus \mathsf{B}_{4
\cdot 2^{-k}}(x)}(z) \hskip1pt \Big( k(x,z) - k(x,\mathrm{c}_y)
\Big) \hskip1pt \psi_{\widehat{Q}_y}(z) \, dz \\ & + & \int_{\R^n}
k(x,\mathrm{c}_y) \hskip1pt \Big( 1_{\R^n \setminus \mathsf{B}_{4
\cdot 2^{-k}}(x)}(z) - 1_{\R^n \setminus \mathsf{B}_{4 \cdot
2^{-k}}(x)}(\mathrm{c}_y) \Big) \hskip1pt \psi_{\widehat{Q}_y}(z)
\, dz \\ [8pt] & = &  \mathsf{A}_1 + \mathsf{B}_1.
\end{eqnarray*}
Here we have used that $\psi_{\widehat{Q}_y}$ is mean-zero.
Lipschitz smoothness gives once more $$|\mathsf{A}_1| \ \le \
\mathrm{c}_n \, \frac{2^{- \gamma (k+s)}}{|x-y|^{n+\gamma}}.$$ To
estimate $\mathsf{B}_1$ we use the size condition on the kernel
$$|\mathsf{B}_1| \, \le \, \frac{\mathrm{c}_n}{|x-\mathrm{c}_y|^n}
\, \frac{1}{|\widehat{Q}_y|} \int_{\widehat{Q}_y} \big| 1_{\R^n
\setminus \mathsf{B}_{4 \cdot 2^{-k}}(x)}(z) - 1_{\R^n \setminus
\mathsf{B}_{4 \cdot 2^{-k}}(x)}(\mathrm{c}_y) \big| \, dz.$$

\noindent
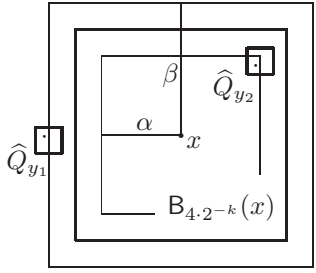
\begin{picture}(360,120)(-180,-60)
\linethickness{0.8pt}
    \put(-2.3,-4.4){$\mbox{\huge$\cdot$}$}
    \put(2,-5){$x$}
    \put(-40,-40){\line(1,0){80}}
    \put(-40,-40){\line(0,1){80}}
    \put(40,40){\line(-1,0){80}}
    \put(40,40){\line(0,-1){80}}

    \put(-55,-7){\line(1,0){10}}
    \put(-55,-7){\line(0,1){10}}
    \put(-45,3){\line(-1,0){10}}
    \put(-45,3){\line(0,-1){10}}

    \put(-53,-3){$\cdot$}

    \put(35,33){\line(-1,0){10}}
    \put(35,33){\line(0,-1){10}}
    \put(25,23){\line(1,0){10}}
    \put(25,23){\line(0,1){10}}

    \put(27,24){$\cdot$}

    \put(-66,-13){$\widehat{Q}_{y_1}$}
    \put(12,15){$\widehat{Q}_{y_2}$}

\linethickness{.2pt}

    \put(40,40){\line(0,-1){80}}

    \put(-30,-30){\line(1,0){20}}
    \put(-30,-30){\line(0,1){60}}
    \put(30,30){\line(-1,0){60}}
    \put(30,30){\line(0,-1){45}}

    \put(0,0){\line(-1,0){30}}
    \put(0,0){\line(0,1){50}}

    \put(-17,2){$\alpha$}
    \put(-7,20){$\beta$}

    \put(-5,-30){$\mathsf{B}_{4\cdot 2^{-k}}(x)$}

    \put(-50,-50){\line(1,0){100}}
    \put(-50,-50){\line(0,1){100}}
    \put(50,50){\line(-1,0){100}}
    \put(50,50){\line(0,-1){100}}

    \put(-108,-70){$(\alpha,\beta) = \big( 4 \cdot 2^{-k} -
    2^{-(k+s-1)}, 4 \cdot 2^{-k} + 2^{-(k+s-1)}\big)$}
\end{picture}

\null

\vskip-15pt

\null
\begin{center}
\textsc{Figure III} \\ If $y \notin \mathcal{W}_{s,k}^x =
\mathsf{B}_{\beta}(x) \setminus \mathsf{B}_\alpha(x)$, we have
$\widehat{Q}_y \cap \partial \hskip1pt\mathsf{B}_{4 \cdot
2^{-k}}(x) = \emptyset$
\end{center}

\vskip5pt

\noindent Since $3 \cdot 2^{-k} < |x-y| \le 5 \cdot 2^{-k}$, we
have $\mathrm{c}_n |x-\mathrm{c}_y|^{-n} \sim \mathrm{c}_n
2^{nk}$. Moreover, the only $z$'s for which the integrand above is
not zero are those with $(z,\mathrm{c}_y)$ lying at different
sides of $\partial \hskip1pt \mathsf{B}_{4 \cdot 2^{-k}}(x)$. This
can only happen when $y \in \mathcal{W}_{s,k}^x$ and we get
$$|\mathsf{B}_1| \le \mathrm{c}_n 2^{nk}
\frac{1_{\mathcal{W}_{s,k}^x(y)}}{|\widehat{Q}_y|}
\int_{\widehat{Q}_y} \big| 1_{\R^n \setminus \mathsf{B}_{4 \cdot
2^{-k}}(x)}(z) - 1_{\R^n \setminus \mathsf{B}_{4 \cdot
2^{-k}}(x)}(\mathrm{c}_y) \big| \, dz \le \mathrm{c}_n 2^{nk}
\hskip1pt 1_{\mathcal{W}_{s,k}^x(y)}.$$ Combining our estimates
obtained so far we get
\begin{equation} \label{k1est}
\big| k_{s,k}^1(x,y) \big| \, \le \, \mathrm{c}_n 1_{\R^n
\setminus \mathsf{B}_{3 \cdot 2^{-k}}(x)}(y) \hskip1pt \Big(
\frac{2^{- \gamma (k+s)}}{|x-y|^{n+\gamma}} + 2^{nk} \hskip1pt
1_{\mathcal{W}_{s,k}^x}(y) \Big).
\end{equation}
Let us now study pointwise estimates for the kernel
$$k_{s,k}^2(x,y) \, = \, \frac{1}{|R_x|} \int_{R_x} \Big(
\int_{\widehat{Q}_y} 1_{\R^n \setminus \mathsf{B}_{4 \cdot
2^{-k}}(w)}(z) \hskip1pt k(w,z) \hskip1pt \psi_{\widehat{Q}_y}(z)
\, dz \Big) \, dw.$$ If $|x-y| \le 2 \cdot 2^{-k}$ we have $$|w-z|
\le |w-x| + |x-y| + |y-z| \le 4 \cdot 2^{-k}$$ for all $(w,z) \in
R_x \times \widehat{Q}_y$. This gives $$k_{s,k}^2(x,y) = 0 \ \,
\quad \mbox{whenever} \ \, \quad |x-y| \le 2 \cdot 2^{-k}.$$ If
$|x-y| > 6 \cdot 2^{-k}$, then we have for $(w,z) \in R_x \times
\widehat{Q}_y$
$$|w-z| \ge |x-y| - |x-w| - |z-y| > 4 \cdot 2^{-k}.$$ Therefore,
we obtain as usual the estimate
\begin{eqnarray*}
\big| k_{s,k}^2(x,y) \big| & \le & \frac{1}{|R_x|} \int_{R_x}
\Big( \int_{\widehat{Q}_y} \big| k(w,z) - k(w,\mathrm{c}_y) \big|
\, \big| \psi_{\widehat{Q}_y}(z) \big| \, dz \Big) \, dw \\ & \le
& \frac{2^{- \gamma (k+s)}}{|R_x|} \int_{R_x} \frac{1}{|w -
\mathrm{c}_y|^{n+\gamma}} \Big( \int_{\widehat{Q}_y} \big|
\psi_{\widehat{Q}_y}(z) \big| \, dz \Big) \, dw \\ & \le &
\mathrm{c}_n \frac{2^{- \gamma (k+s)}}{|R_x|} \int_{R_x}
\frac{dw}{|w - \mathrm{c}_y|^{n+\gamma}} \ = \ \mathrm{c}_n
\frac{2^{- \gamma (k+s)}}{|x - y|^{n+\gamma}}.
\end{eqnarray*}
When $2 \cdot 2^{-k} < |x-y| \le 6 \cdot 2^{-k}$ we have the two
integrals
\begin{eqnarray*}
\lefteqn{k_{s,k}^2(x,y)} \\ & = & \frac{1}{|R_x|} \int_{R_x \times
\widehat{Q}_y} 1_{\R^n \setminus \mathsf{B}_{4 \cdot
2^{-k}}(w)}(z) \hskip1pt \Big( k(w,z) - k(w,\mathrm{c}_y) \Big)
\hskip1pt \psi_{\widehat{Q}_y}(z) \, dw \, dz \\ & + &
\frac{1}{|R_x|} \int_{R_x \times \widehat{Q}_y} k(w,\mathrm{c}_y)
\Big( 1_{\R^n \setminus \mathsf{B}_{4 \cdot 2^{-k}}(w)}(z) -
1_{\R^n \setminus \mathsf{B}_{4 \cdot 2^{-k}}(w)}(\mathrm{c}_y)
\Big) \psi_{\widehat{Q}_y}(z) \, dw \, dz \\ [8pt] & = &
\mathsf{A}_2 + \mathsf{B}_2.
\end{eqnarray*}
By Lipschitz smoothness, we may estimate $\mathsf{A}_2$ by
$$|\mathsf{A}_2| \, \le \, \mathrm{c}_n \hskip1pt
\frac{2^{- \gamma (k+s)}}{|R_x| \, |\widehat{Q}_y|} \int_{R_x
\times \widehat{Q}_y} \frac{1_{\R^n \setminus \mathsf{B}_{4 \cdot
2^{-k}}(w)}(z)}{|w-z|^{n+\gamma}} \, dw \, dz \, \le \,
\mathrm{c}_n 2^{nk} \hskip1pt 2^{- \gamma s}$$ since $4 \cdot
2^{-k} \le |w-z| \le |w-x| + |x-y| + |y-z| \le 8 \cdot 2^{-k}$.
That is
$$|\mathsf{A}_2| \, \le \, \mathrm{c}_n \frac{2^{- \gamma
(k+s)}}{|x-y|^{n+\gamma}} \qquad \big( 2 \cdot 2^{-k} < |x-y| \le
6 \cdot 2^{-k} \big).$$ To estimate $\mathsf{B}_2$ we first
observe that
$$|w-\mathrm{c}_y| \ge |x-y| - |x-w| - |\mathrm{c}_y-y| \ge
(2 - 1 - \mbox{$\frac12$}) \hskip1pt 2^{-k} = \mbox{$\frac12$}
\hskip1pt 2^{-k}.$$ Then we apply the size estimate for the kernel
and Fubini theorem
\begin{eqnarray*}
|\mathsf{B}_2| & \le & \frac{1}{|R_x| \hskip1pt |\widehat{Q}_y|}
\int_{R_x \times \widehat{Q}_y} \frac{\big| 1_{\R^n \setminus
\mathsf{B}_{4 \cdot 2^{-k}}(w)}(z) - 1_{\R^n \setminus
\mathsf{B}_{4 \cdot 2^{-k}}(w)}(\mathrm{c}_y)
\big|}{|w-\mathrm{c}_y|^n} \, dw \, dz \\ & \le &
\frac{\mathrm{c}_n 2^{nk}}{|R_x| \hskip1pt |\widehat{Q}_y|}
\int_{\widehat{Q}_y} \Big( \int_{R_x} \big| 1_{\R^n \setminus
\mathsf{B}_{4 \cdot 2^{-k}}(w)}(z) - 1_{\R^n \setminus
\mathsf{B}_{4 \cdot 2^{-k}}(w)}(\mathrm{c}_y) \big| \, dw \Big) \,
dz.
\end{eqnarray*}
In the integral inside the brackets, the points $z$ and
$\mathrm{c}_y$ are fixed. Moreover, since $z \in \widehat{Q}_y$ we
know that $|z-\mathrm{c}_y| \le 2^{-(k+s)}$. Therefore, we find
that the only $w$'s for which the integrand of the inner integral
is not zero live in
$$\mathcal{W}_{s+1,k}^{\mathrm{c}_y} = \Big\{ w \in \R^n \, \big| \ 4
\cdot 2^{-k} - 2^{-(k+s)} \le |w - \mathrm{c}_y| < 4 \cdot 2^{-k}
+ 2^{-(k+s)} \Big\}.$$ This automatically gives the estimate
$$|\mathsf{B}_2| \, \le \,
\frac{\mathrm{c}_n 2^{nk}}{|R_x|}
|\mathcal{W}_{s+1,k}^{\mathrm{c}_y}| \, \le \, \mathrm{c}_n 2^{nk}
\hskip1pt 2^{-s} \, \le \, \mathrm{c}_n \frac{2^{- \gamma
(k+s)}}{|x-y|^{n+\gamma}}.$$ Our partial estimates so far produce
the global estimate
\begin{equation} \label{k2est}
\big| k_{s,k}^2(x,y) \big| \, \le \, \mathrm{c}_n 1_{\R^n
\setminus \mathsf{B}_{2 \cdot 2^{-k}}(x)}(y) \hskip1pt \frac{2^{-
\gamma (k+s)}}{|x-y|^{n+\gamma}}.
\end{equation}
The assertion then follows from a combination of inequalities
\eqref{k1est} and \eqref{k2est}. \fin

\begin{lemma} \label{Schurtruncated}
Let us define
\begin{eqnarray*}
\mathcal{S}_{s,k}^1(x) & = & \int_{\R^n} \big| k_{s,k}(x,y) \big|
\, dy, \\ \mathcal{S}_{s,k}^2(y) & = & \int_{\R^n} \big|
k_{s,k}(x,y) \big| \, dx.
\end{eqnarray*}
Then there exists a constant $\mathrm{c}_n$ such that $$\max
\Big\{ \mathcal{S}_{s,k}^1(x), \mathcal{S}_{s,k}^2(y) \Big\} \,
\le \, \mathrm{c}_n 2^{- \gamma s}.$$
\end{lemma}

\dem According to Lemma \ref{pointwisetruncated} and
$|\mathcal{W}_{s,k}^x| \le \mathrm{c}_n \hskip1pt 2^{-nk} 2^{-s}$
$$\mathcal{S}_{s,k}^1(x) \, \le \, \mathrm{c}_n \int_{\R^n \setminus
\mathsf{B}_{2 \cdot 2^{-k}}(x)} \Big( \frac{2^{- \gamma
(k+s)}}{|x-y|^{n+\gamma}} + 2^{nk} 1_{\mathcal{W}_{s,k}^x} (y)
\Big) \, dy \, \le \, \mathrm{c}_n 2^{- \gamma s}.$$ The same
argument applies for $\mathcal{S}_{s,k}^2(y)$, since we have
$1_{\mathcal{W}_{s,k}^x}(y) = 1_{\mathcal{W}_{s,k}^y}(x)$. \fin

\subsubsection{Cotlar type estimates}

We have again $\Lambda_{s,i}^{\null} \Lambda_{s,j}^* = 0$ for $i
\neq j$, so that we are reduced (by Cotlar lemma) to estimate the
norms of $\Lambda_{s,i}^* \Lambda_{s,j}^{\null}$ in
$\mathcal{B}(L_2)$. The kernel of $\Lambda_{s,i}^*
\Lambda_{s,j}^{\null}$ is given by $$k_{i,j}^s(x,y) = \int_{\R^n}
\overline{k_{s,i}(z,x)} \hskip1pt k_{s,j}(z,y) \, dz.$$ Taking
$f(z) = 1_{\mathsf{B}_r(y)}(z) / |\mathsf{B}_r(y)|$, we note
$$\frac{1}{|\mathsf{B}_r(y)|} \int_{\mathsf{B}_r(y)} \Big(
\int_{\R^n} k_{s,k}(x,z) \, dx \Big) \, dz = \int_{\R^n} (id -
\mathsf{E}_k) T_{4 \cdot 2^{-k}} \Delta_{k+s} \hskip1pt f(x) \, dx
= 0,$$ due to the integral invariance of conditional expectations.
Taking the limit as $r \to 0$, we deduce from Lebesgue
differentiation theorem that the cancellation condition
\eqref{cancellation2} also holds for our new kernels
$k_{s,k}(x,y)$ and for a.e. $y \in \R^n$. In particular, the same
discussion as above leads us to use \eqref{cancellation2} in one
way or another according to $i \ge j$ or viceversa. Both cases can
be estimated in the same way. Thus we assume in what follows that
$i \ge j$ and use the expression
\begin{eqnarray*}
k_{i,j}^s(x,y) & = & \int_{\R^n} \overline{k_{s,i}(z,x)} \Big(
k_{s,j}(z,y) - k_{s,j}(x,y) \Big) \, dz \\ & = & \int_{\R^n
\setminus \mathsf{B}_{2 \cdot 2^{-j}}(x)} \overline{k_{s,i}(z,x)}
\Big( k_{s,j}(z,y) - k_{s,j}(x,y) \Big) \, dz \\ & + &
\int_{\mathsf{B}_{2 \cdot 2^{-j}}(x) \setminus \mathsf{B}_{2 \cdot
2^{-i}}(x)} \overline{k_{s,i}(z,x)} \Big( k_{s,j}(z,y) -
k_{s,j}(x,y) \Big) \, dz.
\end{eqnarray*}
Observe that the integrand vanishes for $z$ in $\mathsf{B}_{2
\cdot 2^{-i}}(x)$ since $k_{s,i}(z,x)$ does, according to Lemma
\ref{pointwisetruncated}. Let us write $\alpha_{i,j}^s(x,y)$ and
$\beta_{i,j}^s(x,y)$ for the first and second terms on the right.
Then (as before) we need to estimate the quantity $$\sqrt{ \Big(
\big\| \mathcal{S}_{i,j,s}^{1,\alpha} \big\|_\infty + \big\|
\mathcal{S}_{i,j,s}^{1,\beta} \big\|_\infty \Big) \Big( \big\|
\mathcal{S}_{i,j,s}^{2,\alpha} \big\|_\infty + \big\|
\mathcal{S}_{i,j,s}^{2,\beta} \big\|_\infty \Big)},$$ where the
$\mathcal{S}$ functions are given by
\begin{eqnarray*}
\mathcal{S}_{i,j,s}^{1,\alpha}(x) & = & \int_{\R^n} \big|
\alpha_{i,j}^s(x,y) \big| \, dy, \\
\mathcal{S}_{i,j,s}^{2,\alpha}(y) & = & \int_{\R^n} \big|
\alpha_{i,j}^s(x,y) \big| \, dx, \\
\mathcal{S}_{i,j,s}^{1,\beta}(x) & = & \int_{\R^n} \big|
\beta_{i,j}^s(x,y) \big| \, dy, \\
\mathcal{S}_{i,j,s}^{2,\beta}(y) & = & \int_{\R^n} \big|
\beta_{i,j}^s(x,y) \big| \, dx.
\end{eqnarray*}

\begin{lemma} \label{cotlarest3}
We have $$\max \Big\{ \big\| \mathcal{S}_{i,j,s}^{1,\alpha}
\big\|_\infty, \big\| \mathcal{S}_{i,j,s}^{1,\beta} \big\|_\infty
\Big\} \, \le \, \mathrm{c}_n 2^{- 2 \gamma s}.$$
\end{lemma}

\dem According to Lemma \ref{pointwisetruncated}, we know that
$$\big| k_{s,i}(z,x) \big| \le \mathrm{c}_n \Big(
\frac{2^{- \gamma(i+s)}}{|x-z|^{n+\gamma}} + 2^{ni} \hskip1pt
1_{\mathcal{W}_{s,i}^x}(z) \Big)$$ for all $z \in \R^n \setminus
\mathsf{B}_{2 \cdot 2^{-j}}(x)$. Moreover, Lemma
\ref{Schurtruncated} gives $$\int_{\R^n} \big| k_{s,j}(z,y) -
k_{s,j}(x,y) \big| \, dy \le \mathrm{c}_n 2^{- \gamma s}.$$
Combining these estimates we find an $L_\infty$ bound for
$\mathcal{S}_{i,j,s}^{1,\alpha}$
\begin{eqnarray*}
\mathcal{S}_{i,j,s}^{1,\alpha}(x) & \le & \mathrm{c}_n 2^{- \gamma
s} \int_{\R^n \setminus \mathsf{B}_{2 \cdot 2^{-j}}(x)} \Big(
\frac{2^{- \gamma(i+s)}}{|x-z|^{n+\gamma}} + 2^{ni} \hskip1pt
1_{\mathcal{W}_{s,i}^x}(z) \Big) \, dz \\ [8pt] & \le &
\mathrm{c}_n 2^{- 2 \gamma s} \hskip1pt 2^{- \gamma |i-j|} +
\mathrm{c}_n 2^{- \gamma s} 2^{ni} \big| (\R^n \setminus
\mathsf{B}_{2 \cdot 2^{-j}}(x)) \hskip1pt \cap \hskip1pt
\mathcal{W}_{s,i}^x \big|.
\end{eqnarray*}
We claim that $\mathcal{S}_{i,j,s}^{1,\alpha}(x) \le \mathrm{c}_n
2^{- 2 \gamma s} \hskip1pt 2^{- \gamma |i-j|}$. Indeed, if the
intersection above is empty there is nothing to prove. If it is
not empty, the following inequality must hold $$4 \cdot 2^{-i} +
2^{-(i+s-1)} > 2 \cdot 2^{-j}.$$ This implies that we can only
have $i=j$ or $i=j+1$ and hence $$2^{- \gamma s} 2^{ni} \big|
(\R^n \setminus \mathsf{B}_{2 \cdot 2^{-j}}(x)) \hskip1pt \cap
\hskip1pt \mathcal{W}_{s,i}^x \big| \le 2^{- \gamma s} 2^{ni}
\big| \mathcal{W}_{s,i}^x \big| \le 2^{- 2 \gamma s} \sim 2^{- 2
\gamma s} \hskip1pt 2^{- \gamma |i-j|}.$$ Therefore, the first
function clearly satisfies the thesis. Let us now analyze the
second function. To that aim we proceed exactly as above in
$\mathsf{B}_{2 \cdot 2^{-j}}(x) \setminus \mathsf{B}_{2 \cdot
2^{-i}}(x)$ and obtain
\begin{eqnarray*}
\mathcal{S}_{i,j,s}^{1,\beta}(x) & \le & \mathrm{c}_n 2^{- \gamma
s} \int_{\mathsf{B}_{2 \cdot 2^{-j}}(x) \setminus \mathsf{B}_{2
\cdot 2^{-i}}(x)} \Big( \frac{2^{- \gamma(i+s)}}{|x-z|^{n+\gamma}}
+ 2^{ni} \hskip1pt 1_{\mathcal{W}_{s,i}^x}(z) \Big) \, dz \\ [8pt]
& \le & \mathrm{c}_n 2^{- 2 \gamma s} + \mathrm{c}_n 2^{- \gamma
s} 2^{ni} \big| (\R^n \setminus \mathsf{B}_{2 \cdot 2^{-j}}(x))
\hskip1pt \cap \hskip1pt \mathcal{W}_{s,i}^x \big| \ \le \
\mathrm{c}_n 2^{- 2 \gamma s}.
\end{eqnarray*} \fin

\begin{lemma} \label{cotlarest4}
We have $$\max \Big\{ \big\| \mathcal{S}_{i,j,s}^{2,\alpha}
\big\|_\infty, \big\| \mathcal{S}_{i,j,s}^{2,\beta} \big\|_\infty
\Big\} \, \le \, \mathrm{c}_{n,\gamma} \hskip1pt \big( 1 + |i-j|
\big) \hskip1pt 2^{- \gamma |i-j|}.$$
\end{lemma}

\dem For the first function we have
\begin{eqnarray*}
\mathcal{S}_{i,j,s}^{2,\alpha} (y) & \le & \mathrm{c}_n
\int_{\R^n} \Big[ \int_{\R^n \setminus \mathsf{B}_{2 \cdot
2^{-j}}(x)} \Big( \frac{2^{- \gamma(i+s)}}{|x-z|^{n+\gamma}} +
2^{ni} 1_{\mathcal{W}_{s,i}^x}(z) \Big) \big| k_{s,j}(z,y) \big|
\, dz \Big] \, dx \\ & + & \mathrm{c}_n \int_{\R^n} \Big[
\int_{\R^n \setminus \mathsf{B}_{2 \cdot 2^{-j}}(x)} \Big(
\frac{2^{- \gamma(i+s)}}{|x-z|^{n+\gamma}} + 2^{ni}
1_{\mathcal{W}_{s,i}^x}(z) \Big) \big| k_{s,j}(x,y) \big| \, dz
\Big] \, dx \\ & = & \mathrm{c}_n \int_{\R^n} \Big[ \int_{\R^n
\setminus \mathsf{B}_{2 \cdot 2^{-j}}(z)} \Big( \frac{2^{-
\gamma(i+s)}}{|x-z|^{n+\gamma}} + 2^{ni}
1_{\mathcal{W}_{s,i}^z}(x) \Big) \, dx \Big] \big| k_{s,j}(z,y)
\big| \, dz \\ & + & \mathrm{c}_n \int_{\R^n} \Big[ \int_{\R^n
\setminus \mathsf{B}_{2 \cdot 2^{-j}}(x)} \Big( \frac{2^{-
\gamma(i+s)}}{|x-z|^{n+\gamma}} + 2^{ni}
1_{\mathcal{W}_{s,i}^x}(z) \Big) \, dz \Big] \big| k_{s,j}(x,y)
\big| \, dx, \\ [8pt] & \le & \mathrm{c}_n 2^{- 3 \gamma s}
\hskip1pt 2^{- \gamma |i-j|}.
\end{eqnarray*}
Here we have used Lemma \ref{pointwisetruncated}, Fubini theorem
and $1_{\mathcal{W}_{s,i}^x}(z) = 1_{\mathcal{W}_{s,i}^z}(x)$. The
last inequality follows arguing as in Lemma \ref{cotlarest3}. Let
us now estimate the second $\mathcal{S}$ function. We may assume
$i \neq j$ because otherwise $\mathcal{S}_{i,j,s}^{2,\beta} = 0$.
Let us decompose
$$\big| k_{s,j}(z,y) - k_{s,j}(x,y) \big| \, \le \, \mathsf{A} +
\mathsf{B},$$ where these terms are given by
\begin{eqnarray*}
\mathsf{A} & = & \big| k_{s,j}^1(z,y) - k_{s,j}^1(x,y) \big|, \\
\mathsf{B} & = & \big| k_{s,j}^2(z,y) - k_{s,j}^2(x,y) \big|.
\end{eqnarray*}
Moreover, we further decompose the $\mathsf{A}$-term into
\begin{eqnarray*}
\mathsf{A} & \le & \int_{\widehat{Q}_y} 1_{\R^n \setminus
\mathsf{B}_{4 \cdot 2^{-j}}(z)}(w) \hskip1pt \big| k(z,w) - k(x,w)
\big| \hskip1pt \big| \psi_{\widehat{Q}_y}(w) \big| \, dw \\ & + &
\int_{\widehat{Q}_y} \big| k(x,w) \big| \hskip1pt \big| 1_{\R^n
\setminus \mathsf{B}_{4 \cdot 2^{-j}}(z)}(w) - 1_{\R^n \setminus
\mathsf{B}_{4 \cdot 2^{-j}}(x)}(w) \big| \hskip1pt \big|
\psi_{\widehat{Q}_y}(w) \big| \, dw \ = \ \mathsf{A}_1 +
\mathsf{A}_2.
\end{eqnarray*}
This gives rise to
\begin{eqnarray*}
\lefteqn{\mathcal{S}_{i,j,s}^{2,\beta}(y)} \\ & \le & \int_{\R^n}
\Big( \int_{\mathsf{B}_{2 \cdot 2^{-j}}(x) \setminus \mathsf{B}_{2
\cdot 2^{-i}}(x)} \big| k_{s,i}(z,x) \big| \hskip1pt \big(
\mathsf{A}_1 + \mathsf{A}_2 + \mathsf{B} \big) \, dz \Big) \, dx \\
& \le & \mathrm{c}_n \int_{\R^n} \Big( \int_{\mathsf{B}_{2 \cdot
2^{-j}}(x) \setminus \mathsf{B}_{2 \cdot 2^{-i}}(x)} \Big(
\frac{2^{- \gamma (i+s)}}{|x-z|^{n+\gamma}} + 2^{ni}
1_{\mathcal{W}_{s,i}^x}(z) \Big) \hskip1pt \big( \mathsf{A}_1 +
\mathsf{A}_2 + \mathsf{B} \big) \, dz \Big) \, dx \\ [8pt] & = &
\mathcal{A}_1 + \mathcal{A}_2 + \mathcal{B}.
\end{eqnarray*}

\noindent \textbf{The $\mathcal{A}_1$-term.} We have
\begin{eqnarray*}
\mathsf{A}_1 & \le & \int_{\widehat{Q}_y} 1_{\R^n \setminus
\mathsf{B}_{4 \cdot 2^{-j}}(z)}(w) \hskip1pt
\frac{|x-z|^\gamma}{|z-w|^{n+\gamma}} \hskip1pt \big|
\psi_{\widehat{Q}_y}(w) \big| \, dw \\ & \le & \frac{\mathrm{c}_n
2^{nj} \hskip1pt 2^{\gamma j}}{|\widehat{Q}_y|}
\int_{\widehat{Q}_y} \frac{|x-z|^\gamma}{\big( 1 + 2^j |z-w|
\big)^{n+\gamma}} \, dw \ \sim \ \mathrm{c}_n \frac{2^{nj}
\hskip1pt 2^{\gamma j} \hskip1pt |x-z|^\gamma}{\big( 1 + 2^j |z-y|
\big)^{n+\gamma}}.
\end{eqnarray*}
Lipschitz smoothness is applicable since $z \in \mathsf{B}_{2
\cdot 2^{-j}}(x) \setminus \mathsf{B}_{2 \cdot 2^{-i}}(x)$. We
then have
\begin{eqnarray*}
\mathcal{A}_1 & \le & \mathrm{c}_n \int_{\R^n} \Big(
\int_{\mathsf{B}_{2 \cdot 2^{-j}}(x) \setminus \mathsf{B}_{2 \cdot
2^{-i}}(x)} \frac{2^{- \gamma (i+s)}}{|x-z|^{n+\gamma}} \,
\frac{2^{nj} \hskip1pt 2^{\gamma j} \hskip1pt |x-z|^\gamma}{\big(
1 + 2^j |z-y| \big)^{n+\gamma}} \, dz \Big) \, dx \\ & + &
\mathrm{c}_n \int_{\R^n} \Big( \int_{\mathsf{B}_{2 \cdot
2^{-j}}(x) \setminus \mathsf{B}_{2 \cdot 2^{-i}}(x)} 2^{ni}
1_{\mathcal{W}_{s,i}^x}(z) \ \frac{2^{nj} \hskip1pt 2^{\gamma j}
\hskip1pt |x-z|^\gamma}{\big( 1 + 2^j |z-y| \big)^{n+\gamma}} \,
dz \Big) \, dx \\ [8pt] & = & \mathcal{A}_{11} + \mathcal{A}_{12}.
\end{eqnarray*}
The estimate of $\mathcal{A}_{11}$ is standard
\begin{eqnarray*}
\mathcal{A}_{11} & = & \mathrm{c}_n \int_{\R^n} \frac{2^{- \gamma
(i+s)} \hskip1pt 2^{nj} \hskip1pt 2^{\gamma j}}{\big( 1 + 2^j
|z-y| \big)^{n+\gamma}} \Big( \int_{\mathsf{B}_{2 \cdot 2^{-j}}(z)
\setminus \mathsf{B}_{2 \cdot 2^{-i}}(z)} \frac{dx}{|x-z|^n} \Big)
\, dz \\ & = & \mathrm{c}_n |i-j| \int_{\R^n} \frac{2^{- \gamma
(i+s)} \hskip1pt 2^{nj} \hskip1pt 2^{\gamma j}}{\big( 1 + 2^j
|z-y| \big)^{n+\gamma}} \, dz \ \sim \ \mathrm{c}_n |i-j|
\hskip1pt 2^{- \gamma s} \hskip1pt 2^{- \gamma |i-j|}.
\end{eqnarray*}
The term $\mathcal{A}_{12}$ can be written as follows
$$\mathcal{A}_{12} \ = \ \mathrm{c}_n \int_{\R^n} \frac{2^{ni}
\hskip1pt 2^{nj} \hskip1pt 2^{\gamma j}}{\big( 1 + 2^j |z-y|
\big)^{n+\gamma}} \Big( \int_{\mathsf{B}_{2 \cdot 2^{-j}}(z)
\setminus \mathsf{B}_{2 \cdot 2^{-i}}(z)} |x-z|^\gamma
1_{\mathcal{W}_{s,i}^z}(x) \, dx \Big) \, dz.$$ Now, the presence
of $1_{\mathcal{W}_{s,i}^z}(x)$ implies that $$|x-z| \le 4 \cdot
2^{-i} + 2^{-(i+s-1)} \le 5 \cdot 2^{-i}.$$ Therefore we find
$$\mathcal{A}_{12} \ \le \ \mathrm{c}_n 2^{- \gamma |i-j|}
\int_{\R^n} \frac{2^{ni} \hskip1pt 2^{nj} \hskip1pt \big|
\mathcal{W}_{s,i}^z \big|}{\big( 1 + 2^j |z-y| \big)^{n+\gamma}}
\, dz \ \le \ \mathrm{c}_n \hskip1pt 2^{-s} \hskip1pt 2^{- \gamma
|i-j|}.$$ This means that $\mathcal{A}_{11}$ dominates
$\mathcal{A}_{12}$ and we conclude
\begin{equation} \label{estcalA1}
\mathcal{A}_1 \le \mathrm{c}_n |i-j| \hskip1pt 2^{- \gamma s}
\hskip1pt 2^{- \gamma |i-j|}.
\end{equation}

\noindent \textbf{The $\mathcal{A}_2$-term.} Consider the
symmetric difference $$\mathcal{Z}_{x,z}^j \, = \, \mathsf{B}_{4
\cdot 2^{-j}}(x) \bigtriangleup \mathsf{B}_{4 \cdot 2^{-j}}(z) =
\Big( \mathsf{B}_{4 \cdot 2^{-j}}(x) \setminus \mathsf{B}_{4 \cdot
2^{-j}}(z) \Big) \cup \Big( \mathsf{B}_{4 \cdot 2^{-j}}(z)
\setminus \mathsf{B}_{4 \cdot 2^{-j}}(x) \Big).$$ Then we clearly
have $$\mathsf{A}_2 \ = \ \int_{\widehat{Q}_y \cap \hskip1pt
\mathcal{Z}_{x,z}^j} \big| k(x,w) \big| \hskip1pt \big|
\psi_{\widehat{Q}_y}(w) \big| \, dw \ \le \ \mathrm{c}_n 2^{nj}
\frac{\big| \widehat{Q}_y \cap \mathcal{Z}_{x,z}^j \big|}{\big|
\widehat{Q}_y \big|},$$ where the $2^{nj}$ comes from the size
condition on the kernel and the inequality $$|x-w| \ge
\mathrm{dist}(x, \partial \mathcal{Z}_{x,z}^j) \ge 4 \cdot 2^{-j}
- |x-z| \ge 2 \cdot 2^{-j},$$ which holds for any $w \in
\mathcal{Z}_{x,z}^j$ and $z \in \mathsf{B}_{2 \cdot 2^{-j}}(x)
\setminus \mathsf{B}_{2 \cdot 2^{-i}}(x)$. This allows us to write
\begin{eqnarray*}
\mathcal{A}_2 & \le & \mathrm{c}_n 2^{nj} \int_{\R^n} \Big(
\int_{\mathsf{B}_{2 \cdot 2^{-j}}(x) \setminus \mathsf{B}_{2 \cdot
2^{-i}}(x)} \frac{2^{- \gamma (i+s)}}{|x-z|^{n+\gamma}} \hskip1pt
\frac{\big| \widehat{Q}_y \cap \mathcal{Z}_{x,z}^j \big|}{\big|
\widehat{Q}_y \big|} \, dz \Big) \, dx \\ & + & \mathrm{c}_n
2^{nj} \int_{\R^n} \Big( \int_{\mathsf{B}_{2 \cdot 2^{-j}}(x)
\setminus \mathsf{B}_{2 \cdot 2^{-i}}(x)} 2^{ni}
1_{\mathcal{W}_{s,i}^x}(z) \frac{\big| \widehat{Q}_y \cap
\mathcal{Z}_{x,z}^j \big|}{\big| \widehat{Q}_y \big|} \, dz \Big)
\, dx
\\ [8pt] & = & \mathcal{A}_{21} + \mathcal{A}_{22}.
\end{eqnarray*}
Before proceeding with the argument, we note
\begin{itemize}
\item If $|x-y| > 7 \cdot 2^{-j}$ $$|z-w| \ge |x-y| - |x-z| -
|w-y| > 4 \cdot 2^{-j}$$ for all $(w,z) \in \widehat{Q}_y \times
\big( \mathsf{B}_{2 \cdot 2^{-j}}(x) \setminus \mathsf{B}_{2 \cdot
2^{-i}}(x) \big)$. Similarly, we have $$|x-w| \ge |x-y| - |w-y|
> 6 \cdot 2^{-j}.$$ This implies that $w \notin \mathcal{Z}_{x,z}^j$
for any $w \in \widehat{Q}_y$, so that $\widehat{Q}_y \cap
\mathcal{Z}_{x,z}^j = \emptyset$.

\item If $|x-y| < 2^{-j}$ $$|z-w| \le |z-x| + |x-y| + |y-w| < 4
\cdot 2^{-j}$$ for all $(w,z) \in \widehat{Q}_y \times \big(
\mathsf{B}_{2 \cdot 2^{-j}}(x) \setminus \mathsf{B}_{2 \cdot
2^{-i}}(x) \big)$. Similarly, we have $$|x-w| \le |x-y| + |y-w| <
2 \cdot 2^{-j}.$$ This implies that $w \notin \mathcal{Z}_{x,z}^j$
for any $w \in \widehat{Q}_y$, so that $\widehat{Q}_y \cap
\mathcal{Z}_{x,z}^j = \emptyset$.
\end{itemize}
In particular, we conclude that
\begin{eqnarray*}
\lefteqn{\mathcal{A}_{21} + \mathcal{A}_{22}} \\ & = &
\mathrm{c}_n 2^{nj} \int_{\mathsf{B}_{7 \cdot 2^{-j}}(y) \setminus
\mathsf{B}_{2^{-j}}(y)} \Big( \int_{\mathsf{B}_{2 \cdot 2^{-j}}(x)
\setminus \mathsf{B}_{2 \cdot 2^{-i}}(x)} \hskip1pt \frac{2^{-
\gamma (i+s)}}{|x-z|^{n+\gamma}} \hskip1pt \frac{\big|
\widehat{Q}_y \cap \mathcal{Z}_{x,z}^j \big|}{\big| \widehat{Q}_y
\big|} \, dz \Big) \, dx
\\ & + & \mathrm{c}_n 2^{nj} \int_{\mathsf{B}_{7 \cdot 2^{-j}}(y)
\setminus \mathsf{B}_{2^{-j}}(y)} \Big( \int_{\mathsf{B}_{2 \cdot
2^{-j}}(x) \setminus \mathsf{B}_{2 \cdot 2^{-i}}(x)} 2^{ni}
1_{\mathcal{W}_{s,i}^x}(z) \frac{\big| \widehat{Q}_y \cap
\mathcal{Z}_{x,z}^j \big|}{\big| \widehat{Q}_y \big|} \, dz \Big)
\, dx.
\end{eqnarray*}
Observe now that $\widehat{Q}_y$ behaves as a ball of radius
$2^{-(j+s)}$ while $\mathcal{Z}_{x,z}^j$ behaves like an annulus
of radius $4 \cdot 2^{-j}$ and width $|x-z|$. Therefore, the
measure of the intersection can be estimated by $$\big|
\widehat{Q}_y \cap \mathcal{Z}_{x,z}^j \big| \, \le \,
\mathrm{c}_n \min\Big\{ 2^{-(n-1)(j+s)} \hskip1pt |x-z|,
2^{-n(j+s)} \Big\}.$$ This provides us with the estimate
$$\frac{\big| \widehat{Q}_y \cap \mathcal{Z}_{x,z}^j \big|}{\big|
\widehat{Q}_y \big|} \, \le \, \mathrm{c}_n \min \Big\{ 2^{j+s}
\hskip1pt |x-z|, 1 \Big\}.$$ If $2^{-(j+s)} \le 2 \cdot 2^{-i}$
\begin{eqnarray*}
\mathcal{A}_{21} \le \mathrm{c}_n 2^{nj} \hskip1pt 2^{- \gamma
(i+s)} \int_{|x-y| < 7 \cdot 2^{-j}} \Big( \int_{|z-x| >
2^{-(j+s)}} \frac{dz}{|x-z|^{n+\gamma}} \Big) \, dx \le
\mathrm{c}_n 2^{- \gamma |i-j|}.
\end{eqnarray*}
If $2^{-(j+s)} > 2 \cdot 2^{-i}$
\begin{eqnarray*}
\mathcal{A}_{21} & \le & \mathrm{c}_n 2^{nj} \hskip1pt 2^{- \gamma
(i+s)} \int_{|x-y| < 7 \cdot 2^{-j}} \Big( \int_{|z-x| >
2^{-(j+s)}} \frac{dz}{|x-z|^{n+\gamma}} \Big) \, dx \\ [10pt] & +
& \mathrm{c}_n 2^{nj} \hskip1pt 2^{- \gamma (i+s)} \int_{|x-y| < 7
\cdot 2^{-j}} \Big( \int_{\mathsf{B}_{2^{-(j+s)}}(x) \setminus
\mathsf{B}_{2 \cdot 2^{-i}}(x)} \frac{2^{j+s}
|x-z|}{|x-z|^{n+\gamma}} \, dz \Big) \, dx \\ [8pt] & \le &
\mathrm{c}_n 2^{- \gamma |i-j|} \ + \ \mathrm{c}_n \begin{cases}
(i-j - s) \hskip1pt 2^{-|i-j|} & \mbox{if} \ \gamma = 1, \\
\mathrm{c}_\gamma 2^{- \gamma |i-j|} & \mbox{if} \ 0 < \gamma < 1.
\end{cases}
\end{eqnarray*}
This gives $\mathcal{A}_{21} \le \mathrm{c}_{n,\gamma} |i-j|
\hskip1pt 2^{- \gamma |i-j|}$. On the other hand, we also have
$$\mathcal{A}_{22} \, \le \, \mathrm{c}_n 2^{ni} \hskip1pt 2^{nj}
\hskip1pt 2^{j+s} \int_{\mathsf{B}_{7 \cdot 2^{-j}}(y) \setminus
\mathsf{B}_{2^{-j}}(y)} \Big( \int_{\mathsf{B}_{2 \cdot 2^{-j}}(x)
\setminus \mathsf{B}_{2 \cdot 2^{-i}}(x)} |x-z| \hskip1pt
1_{\mathcal{W}_{s,i}^x}(z) \, dz \Big) \, dx.$$ Since we have
$|x-z| < 5 \cdot 2^{-i}$ for $z \in \mathcal{W}_{s,i}^x$ and
$|\mathcal{W}_{s,i}^x| \le \mathrm{c}_n 2^{-ni} \hskip1pt 2^{-s}$,
we get
$$\mathcal{A}_{22} \, \le \, \mathrm{c}_n 2^{-|i-j|}.$$ Therefore,
$\mathcal{A}_{21}$ dominates $\mathcal{A}_{22}$ and we conclude
\begin{equation} \label{estcalA2}
\mathcal{A}_2 \, \le \, \mathrm{c}_{n,\gamma} \hskip1pt |i-j|
\hskip1pt 2^{- \gamma |i-j|}.
\end{equation}

\noindent \textbf{The $\mathcal{B}$-term.} As usual, we decompose
\begin{eqnarray*}
\mathcal{B} & \le & \mathrm{c}_n \int_{\R^n} \Big(
\int_{\mathsf{B}_{2 \cdot 2^{-j}}(x) \setminus \mathsf{B}_{2 \cdot
2^{-i}}(x)} \frac{2^{- \gamma (i+s)}}{|x-z|^{n+\gamma}} \hskip2pt
\mathsf{B} \, dz \Big) \, dx
\\ & + & \mathrm{c}_n \int_{\R^n} \Big( \int_{\mathsf{B}_{2 \cdot
2^{-j}}(x) \setminus \mathsf{B}_{2 \cdot 2^{-i}}(x)} 2^{ni}
1_{\mathcal{W}_{s,i}^x}(z) \hskip1pt \mathsf{B} \, dz \Big) \, dx
\ = \ \mathcal{B}_1 + \mathcal{B}_2,
\end{eqnarray*}
with $\mathsf{B} = \big| k_{s,j}^2(z,y) - k_{s,j}^2(x,y) \big|$.
We have $$\mathcal{B}_1 \, \le \, \mathrm{c}_n \int_{\R^n} \Big(
\int_{\mathsf{B}_{2 \cdot 2^{-j}}(x) \setminus \mathsf{B}_{2 \cdot
2^{-i}}(x)} \frac{2^{- \gamma s} \hskip1pt 2^{ni}}{\big(1 + 2^i
|x-z| \big)^{n+\gamma}} \hskip1pt \mathsf{B} \, dz \Big) \, dx,$$
since for $z \in \mathsf{B}_{2 \cdot 2^{-j}}(x) \setminus
\mathsf{B}_{2 \cdot 2^{-i}}(x)$ both integrands are comparable.
Recalling that $$k_{s,j}^2(x,y) \, = \, \Big\langle T_{4 \cdot
2^{-j}} \psi_{\widehat{Q}_y}, \phi_{R_x} \Big\rangle,$$ we observe
(as in our analysis of $\Phi_s$) that $\mathsf{B} =
\mathsf{E}_j(\mathsf{B})$ when regarded as a function of $z$. This
means that $\mathsf{B} = 0$ for any $z \in R_x$. This, together
with the fact that $\mathsf{B}_{2 \cdot 2^{-j}}(x) \subset 5 R_x$,
implies $$\mathcal{B}_1 \, \le \, \mathrm{c}_n \int_{\R^n} \Big(
\int_{5 R_x \setminus R_x} \frac{2^{- \gamma s} \hskip1pt
2^{ni}}{\big(1 + 2^i |x-z| \big)^{n+\gamma}} \hskip1pt \big|
k_{s,j}^2(z,y) - k_{s,j}^2(x,y) \big| \, dz \Big) \, dx.$$
Moreover, arguing as in Lemma \ref{cotlarest1}
\begin{eqnarray*}
\mathcal{B}_1 & \le & \mathrm{c}_n \int_{\R^n} \Big( \int_{5 R_x
\setminus R_x} \mathsf{E}_j \Big[ \frac{2^{- \gamma s} \hskip1pt
2^{ni}}{\big(1 + 2^i |x- \cdot \, | \big)^{n+\gamma}} \Big](z)
\hskip1pt \big| k_{s,j}^2(z,y) - k_{s,j}^2(x,y) \big| \, dz \Big)
\, dx \\ & \le & \mathrm{c}_n \int_{\R^n} \sum_{R \approx R_x}
\Big( \int_R \frac{2^{- \gamma s} \hskip1pt 2^{ni}}{\big(1 + 2^i
|x-w| \big)^{n+\gamma}} \, dw \Big) \hskip1pt \Big( \big|
k_{s,j}^2(\mathrm{c}_R,y) \big| + \big| k_{s,j}^2(x,y) \big| \Big)
\, dx.
\end{eqnarray*}
Here we write $R \approx R_x$ to denote that $R$ is a dyadic cube
in $\Q_j$ contained in $5R_x \setminus R_x$. Then we apply the
argument in Lemma \ref{cotlarest2}
$$\mathcal{B}_1 \, \le \, \mathrm{c}_n \int_{\R^n} \sum_{R \approx
R_x} \frac{2^{- \gamma s}}{\big(1 + 2^i \delta_x \big)^{\gamma}}
\hskip1pt \Big( \big| k_{s,j}^2(\mathrm{c}_R,y) \big| + \big|
k_{s,j}^2(x,y) \big| \Big) \, dx.$$ Now, it is clear from
\eqref{k2est} that we have $$\big| k_{s,j}^2(x,y) \big| \, \le \,
\mathrm{c}_n \hskip1pt \frac{2^{- \gamma s} \hskip1pt
2^{nj}}{\big( 1 + 2^j |x-y| \big)^{n+\gamma}}.$$ These estimates
in conjunction with the argument in Lemma \ref{cotlarest2} give
$$\mathcal{B}_1 \, \le \, \mathrm{c}_{n,\gamma} \hskip1pt 2^{-
2 \gamma s} \hskip1pt |i-j| \hskip1pt 2^{- \gamma |i-j|}.$$ To
estimate $\mathcal{B}_2$ we use that $\mathsf{B} = 0$ for any $z
\in R_x$ and $\mathsf{B}_{2 \cdot 2^{-j}}(x) \subset 5 R_x$
$$\mathcal{B}_2 \, \le \, \mathrm{c}_n 2^{ni} \int_{\R^n} \Big(
\int_{(5 R_x \setminus R_x) \cap \hskip1pt \mathcal{W}_{s,i}^x}
\big| k_{s,j}^2(z,y) \big|  + \big| k_{s,j}^2(x,y) \big| \, dz
\Big) \, dx.$$ Then we apply our estimate of $|k_{s,j}^2(\, \cdot
\, , \cdot \, )|$ given above for $(z,y)$ and $(x,y)$
\begin{eqnarray*}
\mathcal{B}_2 & \le & \mathrm{c}_n 2^{ni} \hskip1pt 2^{- \gamma s}
\hskip1pt 2^{nj} \int_{\R^n} \Big( \int_{(5 R_x \setminus R_x)
\cap \hskip1pt \mathcal{W}_{s,i}^x} \frac{dz}{\big( 1 + 2^j |z-y|
\big)^{n+\gamma}} \Big) \, dx
\\ & + & \mathrm{c}_n 2^{ni} \hskip1pt  2^{- \gamma s}
\hskip1pt 2^{nj} \int_{\R^n} \Big( \int_{(5 R_x \setminus R_x)
\cap \hskip1pt \mathcal{W}_{s,i}^x} \frac{dz}{\big( 1 + 2^j |x-y|
\big)^{n+\gamma}} \Big) \, dx \\ & \sim & \mathrm{c}_n 2^{ni}
\hskip1pt 2^{- \gamma s} \hskip1pt 2^{nj} \int_{\R^n} \Big(
\int_{(5 R_x \setminus R_x) \cap \hskip1pt \mathcal{W}_{s,i}^x}
\frac{dz}{\big( 1 + 2^j |x-y| \big)^{n+\gamma}} \Big) \, dx.
\end{eqnarray*}
Last equivalence follows from the presence of
$\mathcal{W}_{s,i}^x$. Next, the set $$(5 R_x \setminus R_x) \cap
\hskip1pt \mathcal{W}_{s,i}^x$$ forces $z$ to be outside $R_x$ but
at a distance of $x$ controlled by $4 \cdot 2^{-i} +
2^{-(i+s-1)}$. Thus, the only $x \in \R^n$ for which the inner
integral does not vanish are those $x$ for which $\mathrm{dist}(x,
\partial R_x) \le 4 \cdot 2^{-i} + 2^{-(i+s-1)}$. Notice that for
$|i-j| \le 3$ this suppose no restriction but for $|i-j|$ large
does. Given $R \in \Q_j$ we set
$$R_{s,i} = \Big\{ w \in R \, \big| \ \mathrm{dist}(w, \partial R)
\le 4 \cdot 2^{-i} + 2^{-(i+s-1)} \Big\}.$$

\noindent
\begin{picture}(360,120)(-180,-60)
\linethickness{0.8pt}
    \put(-20,-20){\line(1,0){40}}
    \put(-20,-20){\line(0,1){40}}
    \put(20,20){\line(-1,0){40}}
    \put(20,20){\line(0,-1){40}}

\put(-13,-10){\mbox{\footnotesize $R \! \setminus \!\! R_{s,i}$}}

    \put(-15,-15){\line(1,0){30}}
    \put(-15,-15){\line(0,1){30}}
    \put(15,15){\line(-1,0){30}}
    \put(15,15){\line(0,-1){30}}

\linethickness{.2pt}

    \put(-20,-60){\line(0,1){120}}
    \put(-60,-20){\line(1,0){120}}
    \put(-60,20){\line(1,0){120}}
    \put(20,-60){\line(0,1){120}}

    \put(-15,25){\line(1,0){30}}
    \put(-15,25){\line(0,1){30}}
    \put(15,55){\line(-1,0){30}}
    \put(15,55){\line(0,-1){30}}

    \put(25,25){\line(1,0){30}}
    \put(25,25){\line(0,1){30}}
    \put(55,55){\line(-1,0){30}}
    \put(55,55){\line(0,-1){30}}

    \put(25,-15){\line(1,0){30}}
    \put(25,-15){\line(0,1){30}}
    \put(55,15){\line(-1,0){30}}
    \put(55,15){\line(0,-1){30}}

    \put(25,-55){\line(1,0){30}}
    \put(25,-55){\line(0,1){30}}
    \put(55,-25){\line(-1,0){30}}
    \put(55,-25){\line(0,-1){30}}

    \put(-15,-55){\line(1,0){30}}
    \put(-15,-55){\line(0,1){30}}
    \put(15,-25){\line(-1,0){30}}
    \put(15,-25){\line(0,-1){30}}

    \put(-55,-55){\line(1,0){30}}
    \put(-55,-55){\line(0,1){30}}
    \put(-25,-25){\line(-1,0){30}}
    \put(-25,-25){\line(0,-1){30}}

    \put(-55,-15){\line(1,0){30}}
    \put(-55,-15){\line(0,1){30}}
    \put(-25,15){\line(-1,0){30}}
    \put(-25,15){\line(0,-1){30}}

    \put(-55,25){\line(1,0){30}}
    \put(-55,25){\line(0,1){30}}
    \put(-25,55){\line(-1,0){30}}
    \put(-25,55){\line(0,-1){30}}
\end{picture}

\null

\vskip-15pt

\null
\begin{center}
\textsc{Figure IV} \\ The factor $2^{-|i-j|}$ comes from
$|R_{s,i}| \le \mathrm{c}_n \hskip1pt 2^{-|i-j|} \hskip1pt |R|$
\end{center}

\vskip5pt

\noindent Our considerations allows us to complete our estimate as
follows
\begin{eqnarray*}
\mathcal{B}_2 & \le & \mathrm{c}_n 2^{ni} \hskip1pt 2^{- \gamma s}
\hskip1pt 2^{nj} \sum_{R \in \Q_j} \int_{R_{s,i}} \frac{\big|
\mathcal{W}_{s,i}^x \big|}{\big( 1 + 2^j |x-y| \big)^{n+\gamma}}
\, dx \\ & = & \mathrm{c}_n 2^{ni} \hskip1pt 2^{- \gamma s}
\hskip1pt 2^{nj} \sum_{R \in \Q_j} \frac{|R_{s,i}|}{|R|} \, |R| \,
\frac{1}{|R_{s,i}|} \int_{R_{s,i}} \frac{\big| \mathcal{W}_{s,i}^x
\big|}{\big( 1 + 2^j |x-y| \big)^{n+\gamma}} \, dx \\
& \le & \mathrm{c}_n 2^{-(1+\gamma)s} \hskip1pt 2^{nj} \hskip1pt
2^{-|i-j|} \sum_{R \in \Q_j} |R| \, \frac{1}{|R_{s,i}|}
\int_{R_{s,i}} \frac{1}{\big( 1 + 2^j |x-y| \big)^{n+\gamma}} \,
dx \\ & \sim & \mathrm{c}_n 2^{-(1+\gamma)s} \hskip1pt 2^{nj}
\hskip1pt 2^{-|i-j|} \int_{\R^n} \frac{dx}{\big( 1 + 2^j |x-y|
\big)^{n+\gamma}} \ \sim \ \mathrm{c}_n 2^{-(1+\gamma)s} \hskip1pt
2^{-|i-j|}.
\end{eqnarray*}
Combining our estimates for $\mathcal{B}_1$ and $\mathcal{B}_2$ we
get
\begin{equation} \label{estcalB}
\mathcal{B} \, \le \, \mathrm{c}_{n,\gamma} \hskip1pt 2^{- 2
\gamma s} \hskip1pt |i-j| \hskip1pt 2^{- \gamma |i-j|}
\end{equation}
Finally, the sum of \eqref{estcalA1}, \eqref{estcalA2} and
\eqref{estcalB} produces $$\mathcal{S}_{i,j,s}^{2,\beta}(y) \, \le
\, \mathrm{c}_{n,\gamma} \hskip1pt |i-j| \hskip1pt 2^{- \gamma
|i-j|}.$$ As we have proved that $\mathcal{S}_{i,j,s}^{2,\alpha}$
satisfies a better estimate, the proof is complete. \fin

\begin{remark} \label{Onlypoint}
\emph{The estimate given for $\mathcal{A}_1$ in the proof of Lemma
\ref{cotlarest4} above is the only point in the whole argument for
our pseudo-localization principle where the Lipschitz smoothness
with respect to the $x$ variable is used.}
\end{remark}

\noindent \textbf{Conclusion.} Now we have all the necessary
estimates to complete the argument. Namely, a direct application
of Lemmas \ref{cotlarest3} and \ref{cotlarest4} in conjunction
with Schur lemma give us the following estimate
$$\big\| \Lambda_{s,i}^* \Lambda_{s,j}^{\null}
\big\|_{\mathcal{B}(L_2)} \, \le \, \mathrm{c}_{n,\gamma}
\sqrt{2^{-2 \gamma s} \hskip1pt \big( 1 + |i-j| \big) \hskip1pt
2^{- \gamma |i-j|}} \, \le \, \mathrm{c}_{n,\gamma} \hskip1pt 2^{-
\gamma s} \hskip1pt \alpha_{i-j}^2$$ where $\alpha_k =
(1+|k|)^{\frac14} 2^{- \gamma |k|/4}$. Therefore, Cotlar lemma
provides
$$\|\Psi_s\|_{\mathcal{B}(L_2)} \, = \, \Big\| \summ_k (id -
\mathsf{E}_k) T \Delta_{k+s} \Big\|_{\mathcal{B}(L_2)} \, \le \,
\mathrm{c}_{n,\gamma} \hskip1pt 2^{- \gamma s/2} \summ_k \alpha_k
\, = \, \mathrm{c}_{n,\gamma} \hskip1pt 2^{- \gamma s/2}.$$

\section{Calder{\'o}n-Zygmund decomposition} \label{S4}

We now go back to the noncommutative setting and present a
noncommutative form of Calder{\'o}n-Zygmund decomposition. Let us
recall from the Introduction that, for a given semifinite von
Neumann algebra $\M$ equipped with a \emph{n.s.f.} trace $\tau$,
we shall work on the weak-operator closure $\Mn$ of the algebra
$\Mn_B$ of essentially bounded functions $f: \R^n \to \M$. Recall
also the dyadic filtration $(\Mn_k)_{k \in \Z}$ in $\Mn$.

\subsection{Cuculescu revisited}

A difficulty inherent to the noncommutativity is the absence of
maximal functions. It is however possible to obtain noncommutative
maximal weak and strong inequalities. The strong inequalities
follow by recalling that the $L_p$ norm of a maximal function is
an $L_p(\ell_\infty)$ norm. As observed by Pisier \cite{Pis} and
further studied by Junge \cite{J1}, the theory of operator spaces
is the right tool to define noncommutative $L_p(\ell_\infty)$
spaces; see \cite{JX2} for a nice exposition. We shall be
interested on weak maximal inequalities, which already appeared in
Cuculescu's construction above and are simpler to describe.
Indeed, given a sequence $(f_k)_{k \in \Z}$ of positive functions
in $L_1$ and any $\lambda \in \R_+$, we are interested in
describing the noncommutative form of the Lebesgue measure of
$$\Big\{ \sup_{k \in \Z} f_k > \lambda \Big\}.$$ If $f_k \in
L_1(\Mn)_+$ for $k \in \Z$, this is given by
$$\inf \Big\{ \varphi \big( \1_\Mn - q \big) \, \big| \ q \in
\Mn_\pi, \, q f_k q \le \lambda \hskip1pt q \ \ \mbox{for all} \ \
k \in \Z \Big\}.$$ Given a positive dyadic martingale $f = (f_1,
f_2, \ldots)$ in $L_1(\Mn)$ and looking one more time at
Cuculescu's construction, it is apparent that the projection
$q(\lambda)_k$ represents the following set
$$q(\lambda)_k \sim \Big\{ \sup_{1 \le j \le k} f_j \le \lambda
\Big\}.$$ Therefore, we find $$\1_\Mn - \bigwedge_{k \ge 1}
q(\lambda)_k \sim \Big\{\sup_{k \ge 1} f_k > \lambda \Big\}.$$
However, in this paper we shall be interested in the projection
representing the set where $\sup_{k \in \Z} f_k
> \lambda$ since we will work with the full dyadic filtration
$(\Mn_k)_{k \in \Z}$, where $\Mn_k$ stands for $\mathsf{E}_k
(\Mn)$. We shall clarify below why it is not enough to work with
the truncated filtration $(\Mn_k)_{k \ge 1}$. The construction of
the right projection for $\sup_{k \in \Z} f_k > \lambda$ does not
follow automatically from Cuculescu's construction, see
Proposition \ref{HLWMI} below. Moreover, given a general function
$f \in L_1(\Mn)_+$, we are not able at the time of this writing to
construct the right projections $q_\lambda(f,k)$ which represent
the sets $$q_\lambda(f,k) \sim \Big\{ \sup_{j \in \Z, \, j \le k}
f_j \le \lambda \Big\}.$$ Indeed, the weak$^*$ limit procedure
used in the proof of Proposition \ref{HLWMI} below does not
preserve the commutation relation i) of Cuculescu's construction
and we are forced to work in the following dense class of
$L_1(\Mn)$
\begin{equation} \label{Subspace}
\Mn_{c,+} = \Big\{ f: \R^n \to \M \, \big| \ f \in \Mn_+, \
\overrightarrow{\mathrm{supp}} \hskip1pt f \ \ \mathrm{is \
compact} \Big\} \subset L_1(\Mn).
\end{equation}
Here $\overrightarrow{\mathrm{supp}}$ means the support of $f$ as
a vector-valued function in $\R^n$. In other words, we have
$\overrightarrow{\mathrm{supp}} \hskip1pt f = \mathrm{supp}
\hskip1pt \|f\|_\M$. We employ this terminology to distinguish
from $\mathrm{supp} \, f$ (the support of $f$ as an operator in
$\Mn$) defined in Section \ref{S1}. Note that
$\overrightarrow{\mathrm{supp}} \hskip1pt f$ is a measurable
subset of $\R^n$, while $\mathrm{supp} \hskip1pt f$ is a
projection in $\Mn$. In the rest of the paper we shall work with
functions $f$ in $\Mn_{c,+}$. This impose no restriction due to
the density of $\mathrm{span} \hskip1pt \Mn_{c,+}$ in $L_1(\Mn)$.
The following result is an adaptation of Cuculescu's construction
which will be the one to be used in the sequel.

\begin{lemma} \label{Cuculescu}
Let $f \in \Mn_{c,+}$ and $f_k = \mathsf{E}_k (f)$ for $k \in \Z$.
The sequence $(f_k)_{k \in \Z}$ is a $($positive$)$ dyadic
martingale in $L_1(\Mn)$. Given any positive number $\lambda$,
there exists a decreasing sequence $(q_\lambda(f,k))_{k \in \Z}$
of projections in $\mathcal{A}$ satisfying
\begin{itemize}
\item[i)] $q_\lambda(f,k)$ commutes with $q_\lambda(f,k-1)
\hskip1pt f_k \hskip1pt q_\lambda(f,k-1)$ for each $k \in \Z$.

\item[ii)] $q_\lambda(f,k)$ belongs to $\mathcal{A}_k$ for each $k
\in \Z$ and $q_\lambda(f,k) \hskip1pt f_k \hskip1pt q_\lambda(f,k)
\le \lambda \hskip1pt q_\lambda(f,k)$.

\item[iii)] The following estimate holds $$\varphi \Big(
\mathbf{1}_\Mn - \bigwedge_{k \in \Z} q_\lambda(f,k) \Big) \le
\frac{1}{\lambda} \hskip1pt \|f\|_1.$$
\end{itemize}
\end{lemma}

\dem Since $f \in \Mn_{c,+}$ we have for all $Q \in \Q_j$
$$f_Q = \frac{1}{|Q|} \int_Q f(x) \hskip1pt dx \le 2^j \hskip1pt
\|f\|_\Mn \, \big| \overrightarrow{\mathrm{supp}} \hskip1pt f
\big| \hskip1pt \1_\Mn \longrightarrow 0 \quad \mbox{as} \quad j
\to -\infty.$$ In particular, given any $\lambda \in \R_+$, we
have $f_j \le \lambda \hskip1pt \1_\Mn$ for all $j < m_\lambda <
0$ and certain $m_\lambda \in \Z \setminus \N$ with $|m_\lambda|$
large enough. Then we define the desired projections by the
following relations $$q_\lambda(f,k) = \begin{cases} \1_\Mn &
\mbox{if} \ k < m_\lambda, \\
\chi_{(0,\lambda]}(f_k) & \mbox{if} \ k = m_\lambda, \\
\chi_{(0,\lambda]} \big( q_\lambda(f,k-1) \hskip1pt f_k \hskip1pt
q_\lambda(f,k-1) \big) & \mbox{if} \ k > m_\lambda.
\end{cases}$$ To prove iii) we observe that our projections are
exactly the ones obtained when applying Cuculescu's construction
over the truncated filtration $(\Mn_k)_{k \ge m_\lambda}$. Thus we
get $$\varphi \Big( \mathbf{1}_\Mn - \bigwedge_{k \in \Z}
q_\lambda(f,k) \Big) = \varphi \Big( \mathbf{1}_\Mn - \bigwedge_{k
\ge m_\lambda} q(\lambda)_k \Big) \le \frac{1}{\lambda} \hskip1pt
\sup_{k \ge m_\lambda} \|f_k\|_1 = \frac{1}{\lambda} \hskip1pt
\|f\|_1.$$ The rest of the properties of the sequence
$(q_\lambda(f,k))_{k \in \Z}$ are easily verifiable. \fin

\subsection{The maximal function}

We now recall the Hardy-Littlewood weak maximal inequality. In
what follows it will be quite useful to have another expression
for the $q_\lambda(f,k)$'s constructed in Lemma \ref{Cuculescu}.
It is not difficult to check that $$q_\lambda(f,k) = \sum_{Q \in
\Q_k} \xi_\lambda(f,Q) 1_Q$$ for $k \in \Z$, with
$\xi_\lambda(f,Q)$ projections in $\M$ defined by
$$\xi_\lambda(f,Q) = \begin{cases} \1_\M & \mbox{if} \ k <
m_\lambda, \\ \chi_{(0,\lambda]}(f_Q) & \mbox{if} \ k = m_\lambda,
\\ \chi_{(0,\lambda]} \big( \xi_\lambda(f,\widehat{Q}) f_Q
\xi_\lambda(f,\widehat{Q}) \big) & \mbox{if} \ k > m_\lambda.
\end{cases}$$
As for Cuculescu's construction, we have
\begin{itemize}
\item $\xi_\lambda(f,Q) \in \M_\pi$.

\item $\xi_\lambda(f,Q) \le \xi_\lambda(f,\widehat{Q})$.

\item $\xi_\lambda(f,Q)$ commutes with $\xi_\lambda(f,\widehat{Q})
f_Q \hskip1pt \xi_\lambda(f,\widehat{Q})$.

\item $\xi_\lambda(f,Q) f_Q \hskip1pt \xi_\lambda(f,Q) \le \lambda
\hskip1pt \xi_\lambda(f,Q)$.
\end{itemize}

The noncommutative weak type $(1,1)$ inequality for the
Hardy-Littlewood dyadic maximal function \cite{Me} follows as a
consequence of this. We give a proof including some details
(reported by Quanhua Xu to the author) not appearing in \cite{Me}.

\begin{proposition} \label{HLWMI}
If $(f,\lambda) \in L_1(\Mn) \times \R_+$, there exists
$q_\lambda(f) \in \Mn_\pi$ with $$\sup_{k \in \Z} \big\|
q_\lambda(f) f_k q_\lambda(f) \big\|_{\Mn} \le 16 \hskip1pt
\lambda \qquad \mbox{and} \qquad \varphi \big( \1_\Mn -
q_\lambda(f) \big) \le \frac8\lambda \, \|f\|_1.$$
\end{proposition}

\dem Let us fix an integer $m \in \Z \setminus \N$. Assume $f \in
L_1(\Mn)_+$ and consider the sequence $(q(\lambda)_{m,k})_{k \ge
m}$ provided by Cuculescu's construction applied over the
filtration $(\Mn_k)_{k \ge m}$. Define $$q_m(\lambda) =
\bigwedge_{k \ge m} q(\lambda)_{m,k} \quad \mbox{for each} \quad m
\in \Z \setminus \N.$$ Let us look at the family
$(q_m(\lambda))_{m \in \Z \setminus \N}$. By the weak* compactness
of the unit ball $\mathsf{B}_{\Mn}$ and the positivity of our
family, there must exists a cluster point $a \in
\mathsf{B}_{\Mn_+}$. In particular, we may find a subsequence with
$q_{m_j}(\lambda) \to a$ as $j \to \infty$ (note that $m_j \to -
\infty$ as $j \to \infty$) in the weak$^*$ topology. Then we set
$q_\lambda(f) = \chi_{[1/2,1]}(a)$ and define positive operators
$\delta(a)$ and $\beta(a)$ bounded by $2 \1_\Mn$ and determined by
$$q_\lambda(f) = a \delta(a) = \delta(a) a,$$ $$\1_\Mn - q_\lambda(f) =
\chi_{(1/2,1]}(\1_\Mn - a) = (\1_\Mn - a) \beta(a) = \beta(a)
(\1_\Mn - a).$$ In order to prove the first inequality stated
above, we note that
$$\big\| q_\lambda(f) f_k \hskip1pt q_\lambda(f) \big\|_{\Mn} =
\sup_{\|b\|_{L_1(\Mn)} \le 1} \varphi \big( q_\lambda(f) f_k
\hskip1pt q_\lambda(f) b \big).$$ However, we have
\begin{eqnarray*}
\varphi \big( q_\lambda(f) f_k \hskip1pt q_\lambda(f) b \big) & =
& \varphi \big( a f_k a \hskip1pt \delta(a) b \hskip1pt \delta(a)
\big)
\\ [3pt] & = & \lim_{j \to \infty} \varphi \big( q_{m_j}(\lambda)
f_k \hskip1pt q_{m_j}(\lambda) \delta(a) b \hskip1pt \delta(a)
\big) \\ & \le & \lim_{j \to \infty} \big\| q_{m_j}(\lambda) f_k
\hskip1pt q_{m_j}(\lambda) \big\|_{\infty} \|\delta(a) b \hskip1pt
\delta(a)\|_1.
\end{eqnarray*}
Therefore we conclude $$\varphi \big( q_\lambda(f) f_k \hskip1pt
q_\lambda(f) b \big) \le \|b\|_1 \|\delta(a)\|_\infty^2 \lim_{j
\to \infty} \big\| q(\lambda)_{m_j,k} f_k \hskip1pt
q(\lambda)_{m_j,k} \big\|_\infty \le 4 \lambda.$$ This proves the
first inequality, as for the second
\begin{eqnarray*}
\varphi \big( \1_\Mn - q_\lambda(f) \big) & = & \varphi \big(
(\1_\Mn - a) \beta(a) \big) \, \le \, 2 \hskip1pt \varphi (\1_\Mn
- a) \\ & = & 2 \lim_{j \to \infty} \varphi \big( \1_{\Mn} -
q_{m_j}(\lambda) \big) \, \le \, \frac2\lambda \, \|f\|_1.
\end{eqnarray*}
Finally, for a general $f \in L_1(\Mn)$ we decompose $$f = (f_1 -
f_2) + i (f_3 - f_4)$$ with $f_j \in L_1(\Mn)_+$ and define
$$q_\lambda(f) = \bigwedge_{1 \le j \le 4} q_\lambda(f_j).$$ Then, the
estimate follows easily with constants $16 \hskip1pt \lambda$ and
$8/\lambda$ respectively. \fin

\subsection{The good and bad parts}

If $f \in L_1$ positive and $\lambda \in \R_+$, define
$$M_d f(x) = \sup_{x \in Q \in \Q} \frac{1}{|Q|} \int_Q f(y) \, dy \quad
\mbox{and} \quad \mathrm{E}_\lambda = \Big\{ x \in \R^n \, \big|
\, M_d f(x) > \lambda \Big\}.$$ Writing $\mathrm{E}_\lambda =
\bigcup_j Q_j$ as a disjoint union of maximal dyadic cubes with
$f_{Q} \le \lambda < f_{Q_j}$ for all dyadic $Q \supset Q_j,$ we
may decompose $f = g+b$ where the good and bad parts are given by
$$g = f 1_{\mathrm{E}_{\lambda}^c} + \summ_j f_{Q_j} 1_{Q_j} \quad
\mbox{and} \quad b = \summ_j (f - f_{Q_j}) 1_{Q_j}.$$ If $b_j = (f
- f_{Q_j}) 1_{Q_j}$, we have
\begin{itemize}
\item[i)] $\|g\|_1 \le \|f\|_1$ and $\|g\|_\infty \le 2^n
\lambda$.

\item[ii)] $\mbox{supp} \, b_j \subset Q_j$, $\int_{Q_j} b_j = 0$
and $\sum_j \|b_j\|_1 \le 2 \|f\|_1$.
\end{itemize}

Before proceeding with the noncommutative Calder{\'o}n-Zygmund
decomposition, we simplify our notation for the projections
$\xi_\lambda(f,Q)$ and $q_\lambda(f,k)$. Namely, $(f,\lambda)$
will remain fixed in $\Mn_{c,+} \times \R_+$, see
\eqref{Subspace}. These choices lead us to set
$$\big( \xi_Q, q_k, q \big) = \Big( \xi_\lambda(f,Q),
q_\lambda(f,k), \bigwedge_{k \in \Z} q_\lambda(f,k) \Big).$$
Moreover, we shall write $(p_k)_{k \in \Z}$ for the projections
\begin{equation} \label{pks}
p_k = q_{k-1} - q_k = \sum_{Q \in \Q_k} (\xi_{\widehat{Q}} -
\xi_Q) 1_Q = \sum_{Q \in \Q_k} \pi_Q 1_Q.
\end{equation}
The terminology $\pi_Q = \xi_{\widehat{Q}} - \xi_Q$ will be
frequently used below. Recall that the $p_k$'s are pairwise
disjoint and (according to our new terminology) we have $q_j =
\1_\Mn$ for all $j < m_\lambda$. In particular, we find $$\summ_k
p_k = \1_\Mn - q.$$ If we write $p_\infty$ for $q$ and
$\widehat{\Z}$ stands for $\Z \cup \{\infty\}$, our noncommutative
analogue for the Calder{\'o}n-Zygmund decomposition can be stated as
follows. If $f \in \Mn_{c,+}$ and $\lambda \in \R_+$, we consider
the decomposition $f=g+b$ with
\begin{equation} \label{allterms}
g = \sum_{i,j \in \widehat{\Z}} p_i f_{i \vee j} p_j \quad
\mbox{and} \quad b = \sum_{i,j \in \widehat{\Z}} p_i (f - f_{i
\vee j}) p_j,
\end{equation}
where $i \vee j = \max (i,j)$. Note that $i \vee j = \infty$
whenever $i$ or $j$ is infinite. In particular, since $f=f_\infty$
by definition, the extended sum defining $b$ is just an ordinary
sum over $\Z \times \Z$. Note also that our expressions are
natural generalizations of the classical good and bad parts stated
in the classical decomposition. Indeed, recalling the
orthogonality of the $p_k$'s, all the off-diagonal terms vanish in
the commutative setting and we find something like
\begin{equation} \label{diagonalterms}
g_d = qfq + \summ_k p_k f_k p_k \quad \mbox{and} \quad b_d =
\summ_k p_k (f - f_k) p_k.
\end{equation}
In this form, and recalling that for $\M = \C$ we have $$q \sim
\R^n \setminus \mathrm{E}_\lambda \quad \mbox{and} \quad p_k \sim
\Big\{ Q_j \subset \mathrm{E}_\lambda \, \big| \ Q_j \in \Q_k
\Big\},$$ it is not difficult to see that we recover the classical
decomposition.

\begin{remark} \label{Relabelling}
\emph{In the following we shall use the square-diagram in Figure V
below to think of our decomposition. Namely, we first observe that
for any $f \in \Mn_{c,+}$ and for any $\lambda \in \R_+$ there
will be an $m_\lambda \in \Z$ such that $f_j \le \lambda \hskip1pt
\1_\Mn$ for all $j < m_\lambda$, see the proof of Lemma
\ref{Cuculescu} above. In particular, since $f$ and $\lambda$ will
remain fixed, by a simple relabelling we may assume with no loss
of generality that $m_\lambda = 1$. This will simplify very much
the notation, since now we have $p_k = 0$ for all non-positive
$k$. Therefore, the terms $p_i f_{i \vee j} p_j$ and $p_i (f -
f_{i \vee j}) p_j$ in our decomposition may be located in the
$(i,j)$-th position of an $\infty \times \infty$ matrix where the
\lq last\rq${}$ row and column are devoted to the projection $q =
p_\infty$.}
\end{remark}

\noindent
\begin{picture}(360,100)(-180,-50)
    \put(-101.5,-31.5){\path(0,0)(0,63)(63,63)(63,0)(0,0)}
    \put(-92.5,22.5){\path(0,0)(0,9)}
    \put(-92.5,22.5){\path(0,0)(-9,0)}
    \put(-47.5,-22.5){\path(0,0)(9,0)}
    \put(-47.5,-22.5){\path(0,0)(0,-9)}
    \multiput(-92.5,22.5)(9,-9){5}{\path(0,0)(9,0)(9,-9)(0,-9)(0,0)}

    \put(-101.5,34.5){$p_1 p_2 \ldots \hskip23pt q$}
    \put(-111.5,24.5){$p_1$}
    \put(-111.5,15.5){$p_2$}
    \put(-108.5,1.5){$\vdots$}
    \put(-108.5,-28.5){$q$}

    \put(-97,-43.5){$f=g_d+b_d$}

    \put(39.5,-31.5){\path(0,0)(0,63)(63,63)(63,0)(0,0)}
    \multiput(48.5,-31.5)(9,0){6}{\line(0,1){63}}
    \multiput(39.5,22.5)(0,-9){6}{\line(1,0){63}}

    \put(39.5,34.5){$p_1 p_2 \ldots \hskip23pt q$}
    \put(29.5,24.5){$p_1$}
    \put(29.5,15.5){$p_2$}
    \put(32.5,1.5){$\vdots$}
    \put(32.5,-28.5){$q$}

    \put(51,-43.5){$f = g + b$}
\end{picture}

\null

\vskip-25pt

\null
\begin{center}
\textsc{Figure V} \\ Commutative and noncommutative decompositions
\end{center}

\section{Weak type estimates for diagonal terms}
\label{S5}

In this section we start with the proof of Theorem A. Before that,
a couple of remarks are in order. First, according to the
classical theory it is clearly no restriction to assume that
$q=2$. In particular, since $L_2(\Mn)$ is a Hilbert space valued
$L_2$ space, boundedness in $L_2(\Mn)$ will hold. Second, we may
assume the function $f \in L_1(\Mn)$ belongs to $\Mn_{c,+}$.
Indeed, this follows by decomposing $f$ as a linear combination
$(f_1 - f_2) + i (f_3 - f_4)$ of positive functions $f_j \in
L_1(\Mn)_+$ and using the quasi-triangle inequality on
$L_{1,\infty}(\Mn)$ stated in Section \ref{S1}. Then we
approximate each $f_j \in L_1(\Mn)_+$ by functions in $\Mn_{c,+}$.
Third, since $f \ge 0$ by assumption, we may break it for any
fixed $\lambda \in \R_+$ following our Calder{\'o}n-Zygmund
decomposition. In this section we prove our main result for the
\emph{diagonal} terms in \eqref{diagonalterms}. According to the
quasi-triangle inequality for $L_{1,\infty}(\Mn)$, this will
reduce the problem to estimate the off-diagonal terms.

\subsection{Classical estimates}
\label{diagonalcase}

The standard estimates i) and ii) satisfied by the good and bad
parts of Calder{\'o}n-Zygmund decomposition are satisfied by the
diagonal terms \eqref{diagonalterms}. Indeed, since $f$ is
positive so is $g_d$ and $$\|g_d\|_1 = \varphi(qfq) + \sum_{k \ge
1} \varphi \big( p_k f_k p_k \big) = \varphi \big( fq + f
(\1_{\Mn}-q) \big) \, = \, \|f\|_1.$$ On the other hand, by
orthogonality we have
$$\|g_d\|_\infty \, = \, \max \Big\{ \|qfq\|_\infty, \, \sup_{k
\ge 1} \big\| p_k f_k p_k \big\|_\infty\Big\}.$$ To estimate the
first term, take $a \in L_1(\Mn)$ of norm $1$ with
$$\|qfq\|_\infty \le \varphi(qfq \hskip1pt a) + \delta.$$ Since
$f_k \to f$ as $k \to \infty$ in the weak* topology, we deduce
that
\begin{eqnarray*}
\|qfq\|_\infty & \le & \varphi(qfq \hskip1pt a) + \delta \\ [3pt]
& = & \lim_{k \to \infty} \varphi(q f_k q \hskip1pt a) + \delta \\
& \le & \lim_{k \to \infty} \|q f_k q\|_{\infty} \|a\|_1 + \delta
\, \le \, \lambda + \delta,
\end{eqnarray*}
where the last inequality follows from $q f_k q = q q_k f_k q_k q
\le \lambda q$. Therefore, taking $\delta \to 0^+$ we deduce that
$\|qfq\|_\infty \le \lambda$. Let us now estimate the second
terms. To that aim, we observe that $$f_k = \sum_{Q \in \Q_k}
\frac{1}{|Q|} \int_Q f(y) \, dy \, 1_Q = 2^n \sum_{Q \in \Q_k}
\frac{1}{|\widehat{Q}|} \int_Q f(y) \, dy \, 1_Q \le 2^n
f_{k-1}.$$ Therefore, we obtain
\begin{equation} \label{FullFiltrationNeeded}
\big\| p_k f_k p_k \big\|_{\infty} \le 2^n \big\| q_{k-1} f_{k-1}
q_{k-1} \big\|_{\infty} \le 2^n \lambda.
\end{equation}
This completes the proof of our assertions for $g_d$. Let us now
prove the assertions for $b_d$. If we take $b_{d,k}$ to be $p_k
(f-f_k) p_k$, it is clear that $b_{d,k} = b_{d,k}^*$ and also that
$\mbox{supp} \, b_{d,k} \le p_k$. Moreover, recalling that
$$b_{d,k} = \sum_{Q \in \Q_k} (\xi_{\widehat{Q}} - \xi_Q) (f-f_Q)
(\xi_{\widehat{Q}} - \xi_Q) \, 1_Q,$$ the following identity holds
for any $Q_0 \in \Q_k$ $$\int_{Q_0} b_{d,k}(y) \, dy =
(\xi_{\widehat{Q}_0} - \xi_{Q_0}) \Big( \int_{Q_0} f(y) -
f_{Q_0}(y) \, dy \Big) (\xi_{\widehat{Q}_0} - \xi_{Q_0}) = 0.$$
Finally, we observe that
$$\sum_{k \ge 1} \|b_{d,k}\|_1 \le \sum_{k \ge 1} \varphi \big(
p_k (f+f_k) p_k \big) = 2 \varphi \big( f (\1_{\Mn}-q) \big) \le 2
\|f\|_1.$$ This completes the proof of our assertions for the
function $b_d$. As we shall see in the following section, the
estimates for the off-diagonal terms require more involved
arguments which do not appear in the classical (scalar-valued)
theory.

\begin{remark}
\emph{It is important to note that the \emph{doubling estimate}
\eqref{FullFiltrationNeeded} is crucial for our further analysis
and also that such inequality is the one which imposes to work
with the full filtration $(\Mn_k)_{k \in \Z}$ instead with the
truncated one $(\Mn_k)_{k \ge 1}$. Indeed, if we truncate at $k
\ge 1$ (not at $k \ge m_\lambda$ as we have done), then condition
\eqref{FullFiltrationNeeded} fails in general for $k=1$. This is
another difference with the approach in \cite{PR}, where the
doubling condition above was not needed.}
\end{remark}

\subsection{An $\R^n$-dilated projection}

As above, given a positive function $f \in L_1$, let
$\mathrm{E}_\lambda$ be the set in $\R^n$ where the dyadic
Hardy-Littlewood maximal function $M_df$ is greater than
$\lambda$. If we decompose $\mathrm{E}_\lambda = \bigcup_j Q_j$ as
a disjoint union of maximal dyadic cubes, let us write $9
\hskip1pt \mathrm{E}_\lambda$ for the dilation $$9 \hskip1pt
\mathrm{E}_\lambda = \bigcup_j 9 Q_j.$$ As we pointed out in the
Introduction, this is a key set to give a weak type estimate for
the bad part in Calder{\'o}n-Zygmund decomposition. On the other hand,
we know from Cuculescu's construction that $\1_\Mn-q$ represents
the noncommutative analog of $\mathrm{E}_\lambda$, so that the
noncommutative analog of $9 \hskip1pt \mathrm{E}_\lambda$ should
look like \lq$9 (\1_\Mn - q)$\rq${}$ in the sense that we dilate
on $\R^n$ but not on $\M$. In the following result we construct
the right noncommutative analog of $9 \hskip1pt
\mathrm{E}_\lambda$.

\begin{lemma} \label{keylem}
There exists $\zeta \in \Mn_\pi$ such that
\begin{itemize}
\item[i)] $\lambda \hskip1pt \varphi \big( \1_\Mn - \zeta \big)
\le 9^n \hskip1pt \|f\|_1$.

\item[ii)] If $Q_0 \in \Q$ and $x \in 9 \hskip1pt Q_0$, then
$$\zeta(x) \le \1_\M - \xi_{\widehat{Q}_0} + \xi_{Q_0}.$$ In
particular, in this case we immediately find $\zeta(x) \le
\xi_{Q_0}$.
\end{itemize}
\end{lemma}

\dem Given $k \in \Z_+$, we define
$$\psi_k = \sum_{s=1}^k \sum_{Q \in \Q_s} (\xi_{\widehat{Q}} -
\xi_Q) 1_{9 Q} \quad \mbox{and} \quad \zeta_k = \1_\Mn -
\mbox{supp} \, \psi_k.$$ Since we have $\xi_{Q} \le
\xi_{\widehat{Q}}$ for all dyadic cube $Q$, it turns out that
$(\psi_k)_{k \ge 1}$ is an increasing sequence of positive
operators. However, enlarging $Q$ by its concentric father $9
\hskip1pt Q$ generates overlapping and the $\psi_k$'s are not
projections. This forces us to consider the associated support
projections and define $\zeta_1, \zeta_2, \ldots$ as above. The
sequence of projections $(\zeta_k)_{k \ge 1}$ is clearly
decreasing and we may define $$\zeta = \bigwedge_{k \ge 1}
\zeta_k.$$ Now we are ready to prove the first estimate
\begin{eqnarray*}
\lambda \hskip1pt \varphi \big( \1_\Mn - \zeta \big) & = & \lambda
\lim_{k \to \infty} \varphi \big( \1_{\Mn} - \zeta_k \big) \, \le
\, \lambda \, \sum_{s=1}^\infty \sum_{Q \in \Q_s} \varphi \big(
(\xi_{\widehat{Q}} - \xi_Q) 1_{9Q} \big) \\ & = & 9^n \lambda
\sum_{s=1}^\infty \sum_{Q \in \Q_s} \varphi \big(
(\xi_{\widehat{Q}} - \xi_Q) 1_{Q}\big) \, = \, 9^n \lambda \,
\varphi (\1_{\Mn} - q) \, \le \, 9^n \hskip1pt \|f\|_1.
\end{eqnarray*}
Now fix $Q_0 \in \Q$, say $Q_0 \in \Q_{k_0}$ for some $k_0 \in
\Z$. If $k_0 \le 0$, the assertion is trivial since we know from
Remark \ref{Relabelling} that $\xi_{Q_0} = \xi_{\widehat{Q}_0} =
\1_\M$. Thus, we assume that $k_0 \ge 1$. Then we have
\begin{eqnarray*}
(\xi_{\widehat{Q}_0} - \xi_{Q_0}) 1_{9 Q_0} \le \psi_{k_0} &
\Rightarrow & \zeta_{k_0} \le \1_{\Mn} - (\xi_{\widehat{Q}_0} -
\xi_{Q_0}) 1_{9 Q_0}
\\ & \Rightarrow & \zeta(x) \le \zeta_{k_0}(x) \le \1_\M -
\xi_{\widehat{Q}_0} + \xi_{Q_0}
\end{eqnarray*}
for any $x \in 9 Q_0$. It remains to prove that in fact $\zeta(x)
\le \xi_{Q_0}$. Let us write $Q_j$ for the $j$-th dyadic
antecessor of $Q_0$. In other words, $Q_1$ is the dyadic father of
$Q_0$, $Q_2$ is the dyadic father of $Q_1$ and so on until $Q_{k_0
-1} \in \Q_1$. Since the family $Q_0, Q_1, \ldots$ is increasing,
the same happens for their concentric fathers and we find $$x \in
\bigcap_{j=0}^{k_0-1} 9 \hskip1pt Q_j.$$ In particular, applying
the estimate proved so far $$\zeta(x) \le \bigwedge_{j=0}^{k_0-1}
\big( \1_\M - \xi_{\widehat{Q}_j} + \xi_{Q_j} \big) = \xi_{Q_0}.$$
The last identity easily follows from
$$\xi_{\widehat{Q}_{k_0-1}} = \1_\M.$$ Indeed, we have agreed in
Remark \ref{Relabelling} to assume $q_k = \1_\Mn$ for all $k \le
0$. \fin

\subsection{Chebychev's inequalities}
\label{Chev}

By Paragraph \ref{diagonalcase}, we have
$$\|g_d\|_2^2 = \varphi \big( g_d^{\frac12} g_d^{\null} \hskip1pt
g_d^{\frac12} \big) \le \|g_d\|_1 \|g_d\|_\infty \le 2^n \lambda
\hskip1pt \|f\|_1.$$ In particular, the estimate below follows
from Chebychev's inequality
$$\lambda \hskip1pt \varphi \Big\{ |Tg_d| > \lambda \Big\} \le
\frac{1}{\lambda} \hskip1pt \|T g_d\|_2^2 \lesssim
\frac{1}{\lambda} \hskip1pt \|g_d\|_2^2 \le 2^n \hskip1pt
\|f\|_1.$$ As it is to be expected, here we have used our
assumption on the $L_2$-boundedness of $T$. Now we are interested
on a similar estimate with $b_d$ in place of $g_d$. Using the
projection $\zeta$ introduced in Lemma \ref{keylem}, we may
consider the following decomposition $$T b_d = (\1_\Mn - \zeta) T
(b_d) (\1_\Mn - \zeta) + \zeta \hskip1pt T (b_d) (\1_\Mn - \zeta)
+ (\1_\Mn - \zeta) T (b_d) \zeta + \zeta \hskip1pt T (b_d)
\zeta.$$ In particular, we find
$$\lambda \hskip1pt \varphi \Big\{ \big| T b_d \big| > \lambda
\Big\} \lesssim \lambda \hskip1pt \varphi (\1_\Mn - \zeta) +
\lambda \hskip1pt \varphi \Big\{ \big| \zeta \hskip1pt T (b_d)
\zeta \big| > \lambda \Big\}.$$ Indeed, according to our
decomposition of $T b_d$ and the quasi-triangle inequality on
$L_{1,\infty}(\Mn)$, the estimate above reduces to observe that
the first three terms in such decomposition are left or right
supported by $\1_\Mn - \zeta$. Hence, since the quasi-norm in
$L_{1,\infty}(\Mn)$ is adjoint-invariant \cite{FK}, we easily
deduce it. On the other hand, according to the first estimate in
Lemma \ref{keylem}, it suffices to study the last term above. Let
us analyze the operator $\zeta \hskip1pt T (b_d) \zeta$. In what
follows we shall freely manipulate infinite sums with no worries
of convergence. This is admissible because we may assume from the
beginning (by a simple approximation argument) that $f \in \Mn_n$
for some finite $n \ge 1$. In particular, we could even think that
all our sums are in fact finite sums. We may write $$\zeta
\hskip1pt T (b_d) \zeta = \sum_{k \ge 1} \zeta \hskip1pt T
(b_{k,d}) \zeta$$ with $b_{k,d} = p_k (f - f_k) p_k$ for all $k
\ge 1$. Then, Chebychev's inequality gives
$$\lambda \hskip1pt \varphi \Big\{ \big| \zeta \hskip1pt T (b_d)
\zeta \big| > \lambda \Big\} \le \sum_{k = 1}^\infty \big\| \zeta
\hskip1pt T (b_{d,k}) \zeta \big\|_1.$$ According to Lemma
\ref{keylem} and using $\xi_Q \pi_Q = \pi_Q \xi_Q = 0$ (recall the
definition of $\pi_Q$ from \eqref{pks} above), we have $\zeta(x)
b_{d,k}(y) \zeta(x) = 0$ whenever $x$ lies in the concentric
father $9Q$ of the cube $Q \in \Q_k$ for which $y \in Q$. In other
words, we know that $x$ lives far away from the singularity of the
kernel $k$ and
\begin{eqnarray*}
\Big[ \zeta \hskip1pt T(b_{d,k}) \zeta \Big](x) & = & \int_{\R^n}
k(x,y) \big( \zeta(x) b_{d,k}(y) \zeta(x) \big) \, dy \\ [4pt] & =
& \sum_{Q \in \Q_k} \Big( \int_{Q} k(x,y) \big( \zeta(x)
b_{d,k}(y) \zeta(x) \big) \, dy \Big) \, 1_{(9Q)^c}(x) \\ [1pt] &
= & \zeta(x) \Big( \sum_{Q \in \Q_k} \Big( \int_{Q} k(x,y)
b_{d,k}(y) \, dy \Big) \, 1_{(9Q)^c}(x) \Big) \zeta(x).
\end{eqnarray*}
Now we use the mean-zero condition of $b_{d,k}$ from Paragraph
\ref{diagonalcase}
$$\Big[ \zeta \hskip1pt T (b_{d,k}) \zeta \Big](x) = \zeta(x)
\Big( \sum_{Q \in \Q_k} \Big( \int_{Q} \big( k(x,y) - k(x,c_Q)
\big) b_{d,k}(y) \, dy \Big) \, 1_{(9Q)^c}(x) \Big) \zeta(x),$$
where $c_Q$ is the center of $Q$. Then we use the Lipschitz
$\gamma$-smoothness to obtain
\begin{eqnarray*}
\lefteqn{\sum_{k=1}^\infty \big\| \zeta \hskip1pt T (b_{d,k})
\zeta \big\|_1} \\ & \le & \sum_{k = 1}^\infty \sum_{Q \in \Q_k}
\int_Q \Big\| \big( k(\hskip2pt \cdot,y) - k(\hskip2pt \cdot,c_Q)
\big) b_{d,k}(y) \hskip1pt 1_{(9Q)^c}
(\hskip1pt \cdot \hskip1pt) \Big\|_1 \, dy \\
& = & \sum_{k = 1}^\infty \sum_{Q \in \Q_k} \int_Q \tau \Big[
\Big( \int_{(9Q)^c} \big| k(x,y) - k(x,c_Q) \big| \, dx \Big)
|b_{d,k}(y)| \Big] \, dy \\
& \lesssim & \sum_{k = 1}^\infty \sum_{Q \in \Q_k} \int_Q \Big(
\int_{(9Q)^c} \frac{|y-c_Q|^\gamma}{|x-c_Q|^{n+\gamma}} \, dx
\Big) \tau |b_{d,k}(y)| \, dy \\ & \lesssim & \sum_{k=1}^\infty
\sum_{Q \in \Q_k} \int_Q \tau |b_{d,k}(y)| \, dy \, = \,
\sum_{k=1}^\infty \int_{\R^n} \tau |b_{d,k}(y)| \, dy \, = \,
\sum_{k=1}^\infty \|b_{d,k}\|_1 \, \le \, 2 \|f\|_1,
\end{eqnarray*}
where the last inequality follows once more from Paragraph
\ref{diagonalcase}. This completes the argument for the diagonal
part. Indeed, for any fixed $\lambda \in \R_+$ we have seen that
the diagonal parts of $g$ and $b$ (which depend on the chosen
$\lambda$) satisfy
\begin{equation} \label{weaktineqdiag}
\lambda \hskip1pt \varphi \Big\{ |Tg_d|
> \lambda \Big\} + \lambda \hskip1pt \varphi \Big\{ |Tb_d|
> \lambda \Big\} \le \mathrm{c}_n \, \|f\|_1.
\end{equation}

\section{Weak type estimates for off-diagonal terms}
\label{S6}

Given $\lambda \in \R_+$, we have broken $f$ with our
Calder{\'o}n-Zygmund decomposition for such $\lambda$. In the last
section, we have estimated the diagonal terms $g_d$ and $b_d$. Let
us now consider the off-diagonal terms $g_\mathit{off}$ and
$b_\mathit{off}$ determined by $g = g_d + g_\mathit{off}$ and $b =
b_d + b_\mathit{off}$. As we pointed out in the Introduction, it
is paradoxical that the bad part behaves (when dealing with
off-diagonal terms) better than the good one!

\subsection{An expression for $g_\mathit{off}$}
\label{MDE}

We have
\begin{eqnarray*}
g_{\mathit{off}} & = & \sum_{\begin{subarray}{c} i \neq j
\\ i,j \in \widehat{\Z} \end{subarray}} p_i f_{i \vee j} p_j \\
& = & q f (\1_\Mn-q) + (\1_\Mn - q) f q + \sum_{s=1}^\infty
\sum_{k=1}^\infty p_k f_{k+s} p_{k+s} + p_{k+s} f_{k+s} p_k.
\end{eqnarray*}
Here we have restricted the sum $\sum_{k \in \Z}$ to $\sum_{k \ge
1}$ according to Remark \ref{Relabelling}. Applying property i) of
Cuculescu's construction, we know that the projection $q_j$
commutes with $q_{j-1} f_j q_{j-1}$ for all $j \ge 1$. Taking $i
\wedge j = \min(i,j)$, this immediately gives that $p_i f_{i
\wedge j} p_j = 0$ for $i \neq j$. Indeed, we have
$$\begin{array}{rclcl} p_i f_{i \wedge j} p_j & = & p_i \hskip1pt
q_{i-1} f_i \hskip1pt q_{i-1} p_j & = & 0 \qquad \mbox{if} \qquad
i < j, \\ p_i f_{i \wedge j} p_j & = & p_i q_{j-1} f_j q_{j-1} p_j
& = & 0 \qquad \mbox{if} \qquad i > j.
\end{array}$$
Using this property and inverting the order of summation, we
deduce
\begin{eqnarray*}
\lefteqn{\sum_{s=1}^\infty \sum_{k=1}^\infty p_k f_{k+s} p_{k+s} +
p_{k+s} f_{k+s} p_k} \\ & = & \sum_{s,k=1}^\infty p_k
(f_{k+s}-f_k) p_{k+s} + p_{k+s} (f_{k+s}-f_k) p_k \\ & = &
\sum_{s,k=1}^\infty \sum_{j=1}^s p_k df_{k+j} p_{k+s} + p_{k+s}
df_{k+j} p_k = \sum_{j,k=1}^\infty \sum_{s=j}^\infty p_k df_{k+j}
p_{k+s} + p_{k+s} df_{k+j} p_k.
\end{eqnarray*}
Recall that we may use Fubini theorem since, as we observed in
Paragraph \ref{Chev}, we may even assume that all our sums are
finite sums. Now we can sum in $s$ and apply the commutation
property above to obtain
\begin{eqnarray*}
\lefteqn{\sum_{s=1}^\infty \sum_{k=1}^\infty p_k f_{k+s} p_{k+s} +
p_{k+s} f_{k+s} p_k} \\ & = & \sum_{j,k=1}^\infty p_k df_{k+j}
(q_{k+j-1}-q) + (q_{k+j-1}-q) df_{k+j} p_k \\ & = &
\sum_{j,k=1}^\infty p_k df_{k+j} q_{k+j-1} + q_{k+j-1} df_{k+j}
p_k - \sum_{k=1}^\infty p_k (f-f_k) q + q (f-f_k) p_k
\\ & = & \sum_{j,k=1}^\infty p_k df_{k+j} q_{k+j-1} + q_{k+j-1}
df_{k+j} p_k - \sum_{k=1}^\infty p_k f q + q f p_k \\ & = &
\sum_{j,k=1}^\infty p_k df_{k+j} q_{k+j-1} + q_{k+j-1} df_{k+j}
p_k - (\1_\Mn - q) f q + q f (\1_\Mn - q).
\end{eqnarray*}
Indeed, we have used $$p_k f_k q = p_k q_{k-1} f_k q_{k-1} q = 0 =
q q_{k-1} f_k q_{k-1} p_k = q f_k p_k.$$ Combined the identities
obtained so far, we get
$$g_{\mathit{off}} = \sum_{s=1}^\infty \sum_{k=1}^\infty p_k df_{k+s}
q_{k+s-1} + q_{k+s-1} df_{k+s} p_k = \sum_{s=1}^\infty
\sum_{k=1}^\infty g_{k,s}.$$ We shall use through out this
expression for $g_{\mathit{off}}$ in terms of the functions
$g_{k,s}$.

\subsection{Noncommutative pseudo-localization}

Now we formulate and prove the noncommutative extension of our
pseudo-localization principle. We need a weak notion of support
from \cite{PR} which is quite useful when dealing with weak type
inequalities. For a non-necessarily self-adjoint $f \in \Mn$, the
\emph{two-sided null projection} of $f$ is the greatest projection
$q$ in $\Mn_\pi$ satisfying $qfq = 0$. Then we define the
\emph{weak support projection} of $f$ as
$$\mathrm{supp}^* f = \1_\Mn - q.$$ It is clear that
$\mathrm{supp}^* f = \mathrm{supp} f$ when $\Mn$ is abelian.
Moreover, this notion is weaker than the usual support projection
in the sense that we have $\mathrm{supp}^* f \le \mathrm{supp} f$
for any self-adjoint $f \in \Mn$ and $\mathrm{supp}^* f$ is a
subprojection of both the left and right supports in the
non-self-adjoint case.

\begin{remark} \label{charweaksupp}
\emph{Below we shall use the following characterization of the
weak support projection. The projection $\mathrm{supp}^* f$ is the
smallest projection $p$ in $\Mn_\pi$ satisfying the identity
$$f = p f + f p - p f p.$$ Indeed, let $q$ be the
two-sided null projection of $f$ and let $p = \1_\Mn - q$. Then we
have $(\1_\Mn - p) f (\1_\Mn - p) = 0$ by definition. In other
words, $f = p f + f p - p f p$ and $p$ is the smallest projection
with this property because $q$ is the greatest projection
satisfying the identity $qfq=0$.}
\end{remark}

The following constitutes a noncommutative analog of the
pseudo-localization principle that we have stated in the
Introduction. The terminology has been chosen to fit with that of
the noncommutative Calder{\'o}n-Zygmund decomposition. This will make
the exposition more transparent.

\begin{theorem} \label{pseudolocal}
Let us fix a positive integer $s$. Given a function $f \in
L_2(\Mn)$ and any integer $k$, let us consider any projection
$q_k$ in $\Mn_\pi \cap \Mn_k$ satisfying that $\1_\Mn - q_k$
contains $\mathrm{supp}^* df_{k+s}$ as a subprojection. If we
write $$q_k = \sum_{Q \in \Q_k} \xi_Q 1_Q$$ with $\xi_Q \in
\M_\pi$, we may further consider the projection
$$\zeta_{f,s} = \bigwedge_{k \in \Z} \Big( \1_\Mn -
\bigvee_{Q \in \Q_k} (\1_\M - \xi_Q) \hskip1pt 1_{9Q} \Big).$$
Then we have the following localization estimate in $L_2(\Mn)$
$$\Big( \int_{\R^n} \tau \Big( \left| \big[ \zeta_{f,s} \, T \! f
\, \zeta_{f,s} \big] (x) \right|^2 \Big) \, dx \Big)^{\frac12} \le
\mathrm{c}_{n,\gamma} s \hskip1pt 2^{- \gamma s/4} \Big(
\int_{\R^n} \tau \big( |f(x)|^2 \big) \, dx \Big)^{\frac12},$$ for
any $L_2$-normalized Calder{\'o}n-Zygmund operator with Lipschitz
parameter $\gamma$.
\end{theorem}

\dem We shall reduce this result to its commutative counterpart.
More precisely to the shifted form of the $T1$ theorem proved
above. According to Remark \ref{charweaksupp} and the shift
condition $\mathrm{supp}^* df_{k+s} \prec \1_\Mn - q_k$, we have
$$df_{k+s} = q_k^\perp df_{k+s} + df_{k+s} q_k^\perp - q_k^\perp
df_{k+s} q_k^\perp$$ where we write $q_k^\perp = \1_\Mn - q_k$ for
convenience. On the other hand, let $$\zeta_k = \1_\Mn -
\bigvee_{Q \in \Q_k} (\1_\M - \xi_Q) \hskip1pt 1_{9Q},$$ so that
$\zeta_{f,s} = \bigwedge_k \zeta_k$. Following Lemma, \ref{keylem}
it is easily seen that $\1_\Mn - \zeta_k$ represents the
$\R^n$-dilated projection associated to $\1_\Mn - q_k$ with a
factor $9$. Let $\mathcal{L}_a$ and $\mathcal{R}_a$ denote the
left and right multiplication maps by the operator $a$. Let also
$\mathcal{LR}_a$ stand for $\mathcal{L}_a + \mathcal{R}_a -
\mathcal{L}_a \mathcal{R}_a$ Then our considerations so far and
the fact that $\mathcal{L}_{\zeta_k}, \mathcal{R}_{\zeta_k}$ and
$\mathcal{LR}_{q_k^{\perp}}$ commute with $\mathsf{E}_j$ for $j
\ge k$ give
\begin{eqnarray*}
\lefteqn{\zeta_{f,s} \, T \! f \, \zeta_{f,s}} \\ & = &
\mathcal{L}_{\zeta_{f,s}} \mathcal{R}_{\zeta_{f,s}} \Big( \summ_k
\mathsf{E}_k T \Delta_{k+s} \mathcal{LR}_{q_k^\perp} + \summ_k (id
- \mathsf{E}_k) \mathcal{L}_{\zeta_k} \mathcal{R}_{\zeta_k} T
\mathcal{LR}_{q_k^\perp} \Delta_{k+s} \Big) (f).
\end{eqnarray*}
Now we claim that $$\mathcal{L}_{\zeta_k} \mathcal{R}_{\zeta_k} T
\mathcal{LR}_{q_k^\perp} = \mathcal{L}_{\zeta_k}
\mathcal{R}_{\zeta_k} T_{4 \cdot 2^{-k}}
\mathcal{LR}_{q_k^\perp}.$$ Indeed, this clearly reduces to see
$$\mathcal{L}_{\zeta_k} T
\mathcal{L}_{q_k^\perp} = \mathcal{L}_{\zeta_k} T_{4 \cdot 2^{-k}}
\mathcal{L}_{q_k^\perp} \qquad \mbox{and} \qquad
\mathcal{R}_{\zeta_k} T \mathcal{R}_{q_k^\perp} =
\mathcal{R}_{\zeta_k} T_{4 \cdot 2^{-k}}
\mathcal{R}_{q_k^\perp}.$$ By symmetry, we just prove the first
identity $$\mathcal{L}_{\zeta_k} T \mathcal{L}_{q_k^\perp} f (x)
\, = \sum_{Q \in \Q_k} \zeta_k(x) (\1_\M - \xi_Q) \int_Q k(x,y)
f(y) \, dy.$$ Assume that $x \in 9Q$ for some $Q \in \Q_k$, then
it easily follows from the definition of the projection $\zeta_k$
that $\zeta_k(x) \le \xi_Q$. Note that this is simpler than the
argument in Lemma \ref{keylem} because we do not need to prove
here a property like i) there. In particular, we deduce from the
expression above that for each $y \in Q$ we must have $x \in \R^n
\setminus 9Q$. This implies that $|x-y| \ge 4 \cdot 2^{-k}$ as
desired. Finally, since the operators $\mathcal{L}$ and
$\mathcal{R}$ were created from properties of $f$ and
$\zeta_{f,s}$, we can eliminate them and obtain
$$\zeta_{f,s} \, T \! f \, \zeta_{f,s} \, = \, \mathcal{L}_{\zeta_{f,s}}
\mathcal{R}_{\zeta_{f,s}} \Big( \summ_k \mathsf{E}_k \hskip1pt T
\Delta_{k+s} + \summ_k (id - \mathsf{E}_k) \hskip1pt T_{4 \cdot
2^{-k}} \Delta_{k+s} \Big) (f).$$ Assume that $T^*1 = 0$.
According to our shifted form of the $T1$ theorem, we know that
the operator inside the brackets has norm in $\mathcal{B}(L_2)$
controlled by $\mathrm{c}_{n,\gamma} s \hskip1pt 2^{- \gamma
s/4}$. In particular, the same happens when we tensor with the
identity on $L_2(\M)$, which is the case. This proves the
assertion for convolution-type operators. When $T^*1$ is a
non-zero element of $\mathrm{BMO}$, we may follow verbatim the
paraproduct argument given above since $\mathcal{LR}_{q_k^\perp}$
commutes with $\mathsf{E}_k$ and $\zeta_{f,s} q_k^\perp =
q_k^\perp \zeta_{f,s} = 0$. \fin

\begin{remark}
\emph{It is apparent that $\1_\Mn - q_k$ represents in the
noncommutative setting the set $\Omega_k$ in the commutative
formulation. Moreover, $\zeta_{f,s}$ and $\zeta_k$ represent $\R^n
\setminus \Sigma_{f,s}$ and $\R_n \setminus 9 \hskip1pt \Omega_k$
respectively. The only significant difference is that in the
commutative statement we take $\Omega_k$ to be the \emph{smallest}
$\mathcal{R}_k$-set containing $\mathrm{supp} \hskip2pt df_{k+s}$.
This is done to optimize the corresponding localization estimate.
Indeed, the smaller are the $\Omega_k$'s the larger is $1_{\R^n
\setminus \Sigma_{f,s}} Tf$. However, it is in general false that
the smaller are the $\1_\Mn - q_k$'s the larger is $\zeta_{f,s} \,
T \! \, \zeta_{f,s}$. That is why we consider \emph{any} sequence
of $q_k$'s satisfying the shift condition.}
\end{remark}

\subsection{Estimation of $T g_\mathit{off}$}

Our aim is to estimate $$\lambda \hskip1pt \varphi \Big\{ \big| T
g_{\mathit{off}} \big| > \lambda \Big\}.$$ As usual, we decompose
the term $T g_{\mathit{off}}$ in the following way
$$(\1_\Mn - \zeta) T (g_\mathit{off}) (\1_\Mn - \zeta) + \zeta
\hskip1pt T (g_\mathit{off}) (\1_\Mn - \zeta) + (\1_\Mn - \zeta) T
(g_\mathit{off}) \zeta + \zeta \hskip1pt T (g_\mathit{off})
\zeta,$$ where $\zeta$ denotes the projection constructed in Lemma
\ref{keylem}. According to this lemma and the argument in
Paragraph \ref{Chev}, we are reduced to estimate the last term
above. This will be done in several steps.

\subsubsection{Orthogonality}

It is not difficult to check that the terms $g_{k,s}$ in Paragraph
\ref{MDE} are pairwise orthogonal. It follows from the
trace-invariance of conditional expectations and the mutual
orthogonality of the $p_k$'s. We first prove the following
implication $$\varphi (g_{k,s}^{\null} g_{k',s'}^*) \neq 0 \
\Rightarrow \ k+s = k'+ s'.$$ Indeed, if we assume w.l.o.g. that
$k+s > k'+s'$, we get $$\varphi (g_{k,s}^{\null} g_{k',s'}^*) =
\varphi \big( \mathsf{E}_{k+s-1}(g_{k,s}^{\null} g_{k',s'}^*)
\big) = \varphi \big( \mathsf{E}_{k+s-1}(g_{k,s}^{\null})
g_{k',s'}^* \big) = 0.$$ Now, assume that $k \neq k'$ and $k+s =
k'+s'$. By the orthogonality of the $p_k$'s
\begin{eqnarray*}
\varphi (g_{k,s}^{\null} g_{k',s'}^*) & = & \varphi \big( p_k
df_{k+s} q_{k+s-1} p_{k'} df_{k+s} q_{k+s-1} \big) \\ & + &
\varphi \big( q_{k+s-1} df_{k+s} p_k q_{k+s-1} df_{k+s} p_{k'}
\big) \\ & = & \varphi \big( p_{k'} df_{k+s} q_{k+s-1} p_k
df_{k+s} q_{k+s-1} \big) \\ & + & \varphi \big( q_{k+s-1} df_{k+s}
p_k q_{k+s-1} df_{k+s} p_{k'} \big) \, = \, 0
\end{eqnarray*}
since $p_k q_{k+s-1} = q_{k+s-1} p_k = 0$. This means that
$\varphi (g_{k,s}^{\null} g_{k',s'}^*) = 0$ unless $k = k'$ and
$k+s = k'+s'$ or, equivalently, $(k,s) = (k',s')$. Therefore, the
$g_{k,s}$'s are pairwise orthogonal and
$$\|g_{\mathit{off}}\|_2^2 = \sum_{s=1}^\infty \sum_{k=1}^\infty
\|g_{k,s}\|_2^2.$$

\subsubsection{An $\ell_\infty(\ell_2)$ estimate}

Following the classical argument in Calder{\'o}n-Zygmund decomposition
or our estimate for the diagonal terms in Paragraph
\ref{diagonalcase}, it would suffice to prove that
$\|g_{\mathit{off}}\|_2^2 \lesssim \lambda \hskip1pt \|f\|_1$.
According to the pairwise orthogonality of the $g_{k,s}$'s, that
is to say $$\sum_{s=1}^\infty \sum_{k=1}^\infty \|g_{k,s}\|_2^2
\lesssim \lambda \hskip1pt \|f\|_1.$$ However, we just have the
weaker inequality
\begin{equation} \label{L2estimate}
\sup_{s \ge 1} \, \sum_{k=1}^\infty \|g_{k,s}\|_2^2 \lesssim
\lambda \hskip1pt \|f\|_1.
\end{equation}
Let us prove this estimate before going on with the proof
\begin{eqnarray*}
\|g_{k,s}\|_2^2 & = & 2 \hskip1pt \varphi \big( p_k df_{k+s}
q_{k+s-1} df_{k+s} p_k \big) \\ & = & 2 \hskip1pt \varphi \big(
p_k f_{k+s} q_{k+s-1} f_{k+s} p_k \big) \, - \, 2 \hskip1pt
\varphi \big( p_k f_{k+s} q_{k+s-1} f_{k+s-1} p_k \big) \\ & - & 2
\hskip1pt \varphi \big( p_k f_{k+s-1} q_{k+s-1} f_{k+s} p_k \big)
\, + \, 2 \hskip1pt \varphi \big( p_k f_{k+s-1} q_{k+s-1}
f_{k+s-1} p_k \big).
\end{eqnarray*}
By Cuculescu's construction ii) and $f_j \le 2^n f_{j-1}$ (see
Paragraph \ref{diagonalcase}), we find
$$\begin{array}{rclcl} \displaystyle \big\| f_{k+s}^{\frac12}
\hskip1pt q_{k+s-1}^{\null} f_{k+s}^{\frac12} \big\|_\infty & = &
\big\| q_{k+s-1} f_{k+s} \hskip1pt q_{k+s-1} \big\|_\infty &
\lesssim & \lambda, \\ [3pt] \displaystyle \big\|
f_{k+s-1}^{\frac12} q_{k+s-1}^{\null} f_{k+s-1}^{\frac12}
\big\|_\infty & = & \big\| q_{k+s-1} f_{k+s-1} q_{k+s-1}
\big\|_\infty & \le & \lambda.
\end{array}$$
The crossed terms require H\"older's inequality
\begin{eqnarray*}
\varphi \big( p_k f_{k+s} q_{k+s-1} f_{k+s-1} p_k \big) & \le &
\varphi \big( p_k f_{k+s} q_{k+s-1} f_{k+s} p_k \big)^{\frac12}
\\ & \times & \varphi \big( p_k f_{k+s-1} q_{k+s-1} f_{k+s-1}
p_k \big)^{\frac12} \\ & \lesssim & \lambda \hskip2pt \varphi
\big( p_k f_{k+s} p_k \big)^{\frac12} \varphi \big( p_k f_{k+s-1}
p_k \big)^{\frac12} \, = \, \lambda \hskip1pt \varphi \big( p_k f
p_k \big),
\end{eqnarray*}
where the last identity uses the trace-invariance of the
conditional expectations $\mathsf{E}_{k+s}$ and
$\mathsf{E}_{k+s-1}$ respectively. The same estimate holds for the
remaining crossed term. This proves that
$$\sup_{s \ge 1} \, \sum_{k=1}^\infty \|g_{k,s}\|_2^2 \lesssim
\lambda \, \sup_{s \ge 1} \, \sum_{k=1}^\infty \hskip1pt \varphi
\big( p_k f p_k \big) \le \lambda \hskip1pt \|f\|_1.$$

\subsubsection{The use of pseudo-localization}

Consider the function $$g_{(\!s)} = \summ_k g_{k,s}.$$ It is
straightforward to see that $d {g_{(\! s)}}_{k+s} = g_{k,s}$. In
particular, we have
\begin{equation} \label{ncshift}
\mathrm{supp}^* \, d {g_{(\! s)}}_{k+s} \le p_k = q_{k-1} - q_k
\le \1_\Mn - q_k.
\end{equation}
According to the terminology of Theorem \ref{pseudolocal}, we
consider the projection $$\zeta_{g_{(\! s)},s} = \bigwedge_{k \ge
1} \Big( \1_\Mn - \bigvee_{Q \in \Q_k} (\1_\M - \xi_Q) 1_{9Q}
\Big).$$ Notice that we are just taking $k \ge 1$ and not $k \in
\Z$ as in Theorem \ref{pseudolocal}. This is justified by the fact
that the $q_k$'s are now given by Cuculescu's construction applied
to our $f \in \mathcal{A}_{c,+}$ and our assumption in Remark
\ref{Relabelling} implies that $\1_\M - \xi_Q = 0$ for all $Q \in
\Q_k$ with $k < 1$. Now, if we compare this projection with the
one provided by Lemma \ref{keylem} $$\zeta = \bigwedge_{k \ge 1}
\Big( \1_\Mn - \bigvee_{\begin{subarray}{c} 1 \le j \le k \hskip3pt \\
Q \in \Q_j \end{subarray}} (\xi_{\widehat{Q}} - \xi_Q) 1_{9Q}
\Big),$$ it becomes apparent that $\zeta \le \zeta_{g_{(\!
s)},s}$. On the other hand, Chebychev's inequality gives $$\lambda
\hskip1pt \varphi \Big\{ \big| \zeta \hskip1pt T(g_{\mathit{off}})
\zeta \big| > \lambda \Big\} = \lambda \hskip1pt \varphi \Big\{
\Big| \sum_{s=1}^\infty \zeta \hskip1pt T(g_{(\! s)}) \zeta \Big|
> \lambda \Big\} \le \frac{1}{\lambda} \, \Big[
\sum_{s=1}^\infty \big\| \zeta \hskip1pt T(g_{(\! s)}) \zeta
\big\|_2 \Big]^2.$$ This automatically implies $$\lambda \hskip1pt
\varphi \Big\{ \big| \zeta \hskip1pt T(g_{\mathit{off}}) \zeta
\big| > \lambda \Big\} \le \frac{1}{\lambda} \, \Big[
\sum_{s=1}^\infty \big\| \zeta_{g_{(\! s)},s} \hskip1pt T(g_{(\!
s)}) \zeta_{g_{(\! s)},s} \big\|_2 \Big]^2.$$ Now, combining
\eqref{L2estimate} and \eqref{ncshift}, we may use
pseudo-localization and deduce
\begin{eqnarray*}
\lambda \hskip1pt \varphi \Big\{ \big| \zeta \hskip1pt
T(g_{\mathit{off}}) \zeta \big| > \lambda \Big\} & \le &
\frac{\mathrm{c}_{n,\gamma}^2}{\lambda} \, \Big[ \sum_{s=1}^\infty
s \hskip1pt 2^{- \gamma s/4} \|g_{(\! s)}\|_2 \Big]^2 \\ & = &
\frac{\mathrm{c}_{n,\gamma}^2}{\lambda} \, \Big[ \sum_{s=1}^\infty
s \hskip1pt 2^{- \gamma s/4} \Big( \summ_k \|g_{k,s}\|_2^2
\Big)^{\frac12} \Big]^2 \ \le \ \mathrm{c}_{n,\gamma} \hskip1pt
\|f\|_1.
\end{eqnarray*}
This completes the argument for the off-diagonal terms of $g$.

\subsection{Estimation of $T b_\mathit{off}$}
\label{boff}

As above, it suffices to estimate $$\zeta \hskip1pt T
(b_{\mathit{off}}) \zeta = \sum_{s=1}^\infty \sum_{k=1}^\infty
\zeta \hskip1pt T \Big( p_k (f - f_{k+s}) p_{k+s} + p_{k+s} (f -
f_{k+s}) p_k \Big) \zeta.$$ In the sequel we use the following
notation. For any dyadic cube $Q \in \Q_{k+s}$, we shall write
$Q_k$ to denote the $s$-th antecessor of $Q$. That is, $Q_k$ is
the only dyadic cube in $\Q_k$ containing $Q$. If we set $$b_{k,s}
= p_k (f - f_{k+s}) p_{k+s} + p_{k+s} (f - f_{k+s}) p_k,$$ the
identity below follows from $\xi_{Q_k} \pi_{Q_k} = \pi_{Q_k}
\xi_{Q_k} = 0$ and Lemma \ref{keylem}
\begin{eqnarray*}
\lefteqn{\hskip-10pt \zeta \hskip1pt T (b_{k,s}) \zeta (x)} \\
& = & \int_{\R^n} k(x,y) \hskip1pt \big(
\zeta(x) \hskip1pt b_{k,s}(y) \zeta(x) \big) \, dy \\
[4pt] & = & \zeta(x) \hskip1pt \Big( \sum_{Q \in \Q_{k+s}} \int_Q
k(x,y) \hskip1pt \pi_{Q_k} (f(y) - f_Q) \pi_Q \, dy \
1_{(9Q_k)^c}(x) \Big) \zeta(x) \\ & + & \zeta(x) \hskip1pt \Big(
\sum_{Q \in \Q_{k+s}} \int_{Q} k(x,y) \hskip1pt \pi_Q (f(y) - f_Q)
\pi_{Q_k} \, dy \ 1_{(9Q_k)^c}(x) \Big) \zeta(x) \\ & = & \zeta(x)
\hskip1pt \Big( \sum_{Q \in \Q_{k+s}} \int_{Q} k(x,y) \hskip1pt
b_{k,s}(y) \, dy \ 1_{(9Q_k)^c}(x) \Big) \zeta(x).
\end{eqnarray*}

Before going on with the proof, let us explain a bit our next
argument. Our terms $b_{s,k}$ are located in the $(s+1)$-th upper
and lower diagonals and we want to compare their \emph{size} with
that of the main diagonal. To do so we write each $b_{k,s}$,
located in the entries $(k,k+s)$ and $(k+s,k)$, as a linear
combination of four diagonal boxes in a standard way. However,
this procedure generates overlapping and we are forced to consider
only those integers $k$ congruent to a fixed $1 \le j \le s+1$ at
a time. The figure below will serve as a model ($s=2$) for our
forthcoming estimates.

\noindent
\begin{picture}(360,150)(-180,-100)
    \multiput(-151.35,25.5)(7,-7){7}{$\bullet$}
    \multiput(-165.35,11.5)(7,-7){7}{$\bullet$}

    \put(-166.5,-31.5){\path(0,0)(0,63)(63,63)(63,0)(0,0)}
    \put(-159.5,24.5){\path(0,0)(0,7)}
    \put(-159.5,24.5){\path(0,0)(-7,0)}
    \put(-110.5,-24.5){\path(0,0)(7,0)}
    \put(-110.5,-24.5){\path(0,0)(0,-7)}
    \multiput(-159.5,24.5)(7,-7){7}{\path(0,0)(7,0)(7,-7)(0,-7)(0,0)}
    \put(-152.5,24.5){\path(0,0)(0,7)}
    \put(-110.5,-17.5){\path(0,0)(7,0)}
    \multiput(-152.5,24.5)(7,-7){6}{\path(0,0)(7,0)(7,-7)}
    \put(-159.5,17.5){\path(0,0)(-7,0)}
    \put(-117.5,-24.5){\path(0,0)(0,-7)}
    \multiput(-159.5,17.5)(7,-7){6}{\path(0,0)(0,-7)(7,-7)}
    \put(-145.5,24.5){\path(0,0)(0,7)}
    \put(-110.5,-10.5){\path(0,0)(7,0)}
    \multiput(-145.5,24.5)(7,-7){5}{\path(0,0)(7,0)(7,-7)}
    \put(-159.5,10.5){\path(0,0)(-7,0)}
    \put(-124.5,-24.5){\path(0,0)(0,-7)}
    \multiput(-159.5,10.5)(7,-7){5}{\path(0,0)(0,-7)(7,-7)}
    \put(-166,-43.5){$3$-th \emph{diagonals}}

    \put(-93,-2.5){=}

    \multiput(-61.35,25.5)(21,-21){3}{$\bullet$}
    \multiput(-75.35,11.5)(21,-21){3}{$\bullet$}

    \put(-76.5,-31.5){\path(0,0)(0,63)(63,63)(63,0)(0,0)}
    \put(-69.5,24.5){\path(0,0)(0,7)}
    \put(-69.5,24.5){\path(0,0)(-7,0)}
    \put(-20.5,-24.5){\path(0,0)(7,0)}
    \put(-20.5,-24.5){\path(0,0)(0,-7)}
    \multiput(-69.5,24.5)(7,-7){7}{\path(0,0)(7,0)(7,-7)(0,-7)(0,0)}
    \put(-62.5,24.5){\path(0,0)(0,7)}
    \put(-20.5,-17.5){\path(0,0)(7,0)}
    \multiput(-62.5,24.5)(7,-7){6}{\path(0,0)(7,0)(7,-7)}
    \put(-69.5,17.5){\path(0,0)(-7,0)}
    \put(-27.5,-24.5){\path(0,0)(0,-7)}
    \multiput(-69.5,17.5)(7,-7){6}{\path(0,0)(0,-7)(7,-7)}
    \put(-55.5,24.5){\path(0,0)(0,7)}
    \put(-20.5,-10.5){\path(0,0)(7,0)}
    \multiput(-55.5,24.5)(7,-7){5}{\path(0,0)(7,0)(7,-7)}
    \put(-69.5,10.5){\path(0,0)(-7,0)}
    \put(-34.5,-24.5){\path(0,0)(0,-7)}
    \multiput(-69.5,10.5)(7,-7){5}{\path(0,0)(0,-7)(7,-7)}
    \put(-72,-43.5){$k \equiv 1 \ \mathrm{mod} \, 3$}

    \put(-3,-2.5){+}

    \multiput(35.65,18.5)(21,-21){2}{$\bullet$}
    \multiput(21.65,4.5)(21,-21){2}{$\bullet$}

    \put(13.5,-31.5){\path(0,0)(0,63)(63,63)(63,0)(0,0)}
    \put(20.5,24.5){\path(0,0)(0,7)}
    \put(20.5,24.5){\path(0,0)(-7,0)}
    \put(69.5,-24.5){\path(0,0)(7,0)}
    \put(69.5,-24.5){\path(0,0)(0,-7)}
    \multiput(20.5,24.5)(7,-7){7}{\path(0,0)(7,0)(7,-7)(0,-7)(0,0)}
    \put(27.5,24.5){\path(0,0)(0,7)}
    \put(69.5,-17.5){\path(0,0)(7,0)}
    \multiput(27.5,24.5)(7,-7){6}{\path(0,0)(7,0)(7,-7)}
    \put(20.5,17.5){\path(0,0)(-7,0)}
    \put(62.5,-24.5){\path(0,0)(0,-7)}
    \multiput(20.5,17.5)(7,-7){6}{\path(0,0)(0,-7)(7,-7)}
    \put(34.5,24.5){\path(0,0)(0,7)}
    \put(69.5,-10.5){\path(0,0)(7,0)}
    \multiput(34.5,24.5)(7,-7){5}{\path(0,0)(7,0)(7,-7)}
    \put(20.5,10.5){\path(0,0)(-7,0)}
    \put(55.5,-24.5){\path(0,0)(0,-7)}
    \multiput(20.5,10.5)(7,-7){5}{\path(0,0)(0,-7)(7,-7)}
    \put(18.5,-43.5){$k \equiv 2 \ \mathrm{mod} \, 3$}

    \put(87,-2.5){+}

    \multiput(132.65,11.5)(21,-21){2}{$\bullet$}
    \multiput(118.65,-2.5)(21,-21){2}{$\bullet$}

    \put(103.5,-31.5){\path(0,0)(0,63)(63,63)(63,0)(0,0)}
    \put(110.5,24.5){\path(0,0)(0,7)}
    \put(110.5,24.5){\path(0,0)(-7,0)}
    \put(159.5,-24.5){\path(0,0)(7,0)}
    \put(159.5,-24.5){\path(0,0)(0,-7)}
    \multiput(110.5,24.5)(7,-7){7}{\path(0,0)(7,0)(7,-7)(0,-7)(0,0)}
    \put(117.5,24.5){\path(0,0)(0,7)}
    \put(159.5,-17.5){\path(0,0)(7,0)}
    \multiput(117.5,24.5)(7,-7){6}{\path(0,0)(7,0)(7,-7)}
    \put(110.5,17.5){\path(0,0)(-7,0)}
    \put(152.5,-24.5){\path(0,0)(0,-7)}
    \multiput(110.5,17.5)(7,-7){6}{\path(0,0)(0,-7)(7,-7)}
    \put(124.5,24.5){\path(0,0)(0,7)}
    \put(159.5,-10.5){\path(0,0)(7,0)}
    \multiput(124.5,24.5)(7,-7){5}{\path(0,0)(7,0)(7,-7)}
    \put(110.5,10.5){\path(0,0)(-7,0)}
    \put(145.5,-24.5){\path(0,0)(0,-7)}
    \multiput(110.5,10.5)(7,-7){5}{\path(0,0)(0,-7)(7,-7)}
    \put(109,-43.5){$k \equiv 3 \ \mathrm{mod} \, 3$}

    \put(-94.5,-77){\path(0,0)(21,0)}
    \put(-94.5,-84){\path(0,0)(21,0)}
    \put(-87.5,-70){\path(0,0)(0,-21)}
    \put(-80.5,-70){\path(0,0)(0,-21)}

    \put(-79.45,-76.2){$\bullet$}
    \put(-93.2,-89.9){$\bullet$}

    \put(-67,-83){=}

    \put(-52.5,-77){\path(0,0)(21,0)}
    \put(-52.5,-84){\path(0,0)(21,0)}
    \put(-45.5,-70){\path(0,0)(0,-21)}
    \put(-38.5,-70){\path(0,0)(0,-21)}

    \put(-51.2,-76.2){$\bullet$}
    \put(-51.2,-89.9){$\bullet$}
    \put(-37.45,-76.2){$\bullet$}
    \put(-37.45,-89.9){$\bullet$}
    \put(-44.2,-76.2){$\bullet$}
    \put(-51.2,-83.2){$\bullet$}
    \put(-44.2,-83.2){$\bullet$}
    \put(-44.2,-90.2){$\bullet$}
    \put(-37.45,-83.2){$\bullet$}

    \put(-23.5,-83){--}

    \put(-10.5,-77){\path(0,0)(21,0)}
    \put(-10.5,-84){\path(0,0)(21,0)}
    \put(-3.5,-70){\path(0,0)(0,-21)}
    \put(3.5,-70){\path(0,0)(0,-21)}

    \put(-9.2,-76.2){$\bullet$}
    \put(-2.2,-76.2){$\bullet$}
    \put(-9.2,-83.2){$\bullet$}
    \put(-2.2,-83.2){$\bullet$}

    \put(18.5,-83){--}

    \put(31.5,-77){\path(0,0)(21,0)}
    \put(31.5,-84){\path(0,0)(21,0)}
    \put(38.5,-70){\path(0,0)(0,-21)}
    \put(45.5,-70){\path(0,0)(0,-21)}

    \put(39.8,-89.9){$\bullet$}
    \put(39.8,-82.9){$\bullet$}
    \put(46.55,-89.9){$\bullet$}
    \put(46.55,-82.9){$\bullet$}

    \put(59,-83){+}

    \put(73.5,-77){\path(0,0)(21,0)}
    \put(73.5,-84){\path(0,0)(21,0)}
    \put(80.5,-70){\path(0,0)(0,-21)}
    \put(87.5,-70){\path(0,0)(0,-21)}

    \put(81.8,-82.9){$\bullet$}

    \Thicklines

    \put(-76.5,31.5){\path(0,0)(0,-21)(21,-21)(21,0)(0,0)}
    \put(-55.5,10.5){\path(0,0)(0,-21)(21,-21)(21,0)(0,0)}
    \put(-34.5,-10.5){\path(0,0)(0,-21)(21,-21)(21,0)(0,0)}

    \put(20.5,24.5){\path(0,0)(0,-21)(21,-21)(21,0)(0,0)}
    \put(41.5,3.5){\path(0,0)(0,-21)(21,-21)(21,0)(0,0)}
    \put(62.5,-17.5){\path(0,0)(14,0)}
    \put(62.5,-17.5){\path(0,0)(0,-14)}

    \put(117.5,17.5){\path(0,0)(0,-21)(21,-21)(21,0)(0,0)}
    \put(138.5,-3.5){\path(0,0)(0,-21)(21,-21)(21,0)(0,0)}
    \put(159.5,-24.5){\path(0,0)(7,0)}
    \put(159.5,-24.5){\path(0,0)(0,-7)}

    \put(-94.5,-70){\path(0,0)(0,-21)(21,-21)(21,0)(0,0)}
    \put(-52.5,-70){\path(0,0)(0,-21)(21,-21)(21,0)(0,0)}
    \put(-10.5,-70){\path(0,0)(0,-21)(21,-21)(21,0)(0,0)}
    \put(31.5,-70){\path(0,0)(0,-21)(21,-21)(21,0)(0,0)}
    \put(73.5,-70){\path(0,0)(0,-21)(21,-21)(21,0)(0,0)}
\end{picture}

\begin{center}
\textsc{Figure VI} \\ Decomposition into disjoint diagonal boxes
for $s=2$
\end{center}

\noindent According to Chebychev's inequality we obtain
\begin{eqnarray*}
\lefteqn{\hskip-30pt \lambda \hskip1pt \varphi \Big\{ \Big|
\sum_{s=1}^\infty \sum_{k=1}^\infty \zeta \hskip1pt T (b_{k,s})
\zeta \Big| > \lambda \Big\}} \\ & \le & \sum_{s=1}^\infty \Big\|
\sum_{k=1}^\infty
\zeta \hskip1pt T (b_{k,s}) \zeta \Big\|_1 \\
& \le & \sum_{s=1}^\infty \sum_{k=1}^\infty \sum_{Q \in \Q_{k+s}}
\Big\| \int_{Q} k(\hskip2pt \cdot ,y) \hskip1pt b_{k,s}(y) \, dy \
1_{(9Q_k)^c}(\hskip1pt \cdot \hskip1pt) \Big\|_1 \\ & \le &
\sum_{s=1}^\infty \sum_{j=0}^s \sum_{\begin{subarray}{c} k \equiv
j \\ \!\!\!\!\!\!\! \mod s+1 \end{subarray}} \sum_{Q \in \Q_{k+s}}
\Big\| \int_{Q} k(\hskip2pt \cdot,y) \hskip1pt b_{k,s}(y) dy \,
1_{(9Q_k)^c}(\hskip1pt \cdot \hskip1pt) \Big\|_1.
\end{eqnarray*}
We now use the decomposition
\begin{eqnarray*}
b_{k,s} & = & \Big( \sum_{r=0}^s p_{k+r} \Big) (f - f_{k+s}) \Big(
\sum_{r=0}^s p_{k+r} \Big) \\ & - & \Big( \sum_{r=0}^{s-1} p_{k+r}
\Big) (f - f_{k+s}) \Big( \sum_{r=0}^{s-1} p_{k+r} \Big) \\ & - &
\Big( \sum_{r=1}^s p_{k+r} \Big) (f - f_{k+s}) \Big( \sum_{r=1}^s
p_{k+r} \Big) \\ & + & \Big( \sum_{r=1}^{s-1} p_{k+r} \Big) (f -
f_{k+s}) \Big( \sum_{r=1}^{s-1} p_{k+r} \Big) \, = \, b_{k,s}^1 -
b_{k,s}^2 - b_{k,s}^3 + b_{k,s}^4,
\end{eqnarray*}
of $b_{k,s}$ as a linear combination of four \emph{diagonal}
terms. Let us recall that the four projections $\sum_r p_{k+r}$
above (with $0 \preceq r \preceq s$ and $\preceq$ meaning either
$<$ or $\le$) belong to $\Mn_{k+s}$. In particular, since
$\mathsf{E}_{k+s}(f - f_{k+s}) = 0$, the following identity holds
for any $Q \in \Q_{k+s}$ and any $1 \le i \le 4$ $$\int_Q
b_{k,s}^i(y) \, dy = 0.$$ Therefore, we find
\begin{eqnarray*}
\lefteqn{\lambda \hskip1pt \varphi \Big\{ \Big| \sum_{s=1}^\infty
\sum_{k=1}^\infty \zeta \hskip1pt T (b_{k,s}) \zeta \Big| >
\lambda \Big\}} \\ & \le & \sum_{s=1}^\infty
\sum_{i=1}^4 \sum_{j=0}^s \sum_{\begin{subarray}{c} k \equiv j \\
\!\!\!\!\!\!\! \mod s+1 \end{subarray}} \sum_{Q \in \Q_{k+s}}
\int_{Q} \Big\| \big( k(\hskip2pt \cdot,y) - k(\hskip2pt
\cdot,c_Q) \big) \hskip1pt b^i_{k,s}(y) \, 1_{(9Q_k)^c}(\hskip1pt
\cdot \hskip1pt) \Big\|_1 \, dy.
\end{eqnarray*}
However, by Lipschitz $\gamma$-smoothness we have
\begin{eqnarray*}
\lefteqn{\hskip-40pt \int_{Q} \Big\| \big( k(\hskip2pt \cdot,y) -
k(\hskip2pt \cdot,c_Q) \big) \hskip1pt b^i_{k,s}(y) \,
1_{(9Q_k)^c}(\hskip1pt
\cdot \hskip1pt) \Big\|_1 \, dy} \\
& = & \int_Q \tau \Big[ \Big( \int_{(9Q_k)^c} \big| k(x,y) -
k(x,c_Q) \big| \, dx \Big) | b^i_{k,s}(y)|  \Big] \, dy \\ &
\lesssim & \int_Q \Big( \int_{(9Q_k)^c} \frac{|y - c_Q|^\gamma}{|x
- c_Q|^{n+\gamma}} \, dx \Big) \hskip1pt \tau | b^i_{k,s}(y)| \,
dy
\\ & \lesssim & \ell(Q)^\gamma / \ell(Q_k)^\gamma \, \ \varphi
\Big[ \Big( \sum_{r=0}^s p_{k+r}
\Big) (f + f_{k+s}) \Big( \sum_{r=0}^s p_{k+r} \Big) 1_Q \Big] \\
& \lesssim & 2^{-\gamma s} \varphi \Big[ \Big( \sum_{r=0}^s
p_{k+r} \Big) f \Big( \sum_{r=0}^s p_{k+r} \Big) 1_Q \Big].
\end{eqnarray*}
Finally, summing over $(s,i,j,k,Q)$ we get
\begin{eqnarray*}
\lefteqn{\hskip-20pt \lambda \hskip1pt \varphi \Big\{ \Big|
\sum_{s=1}^\infty \sum_{k=1}^\infty \zeta \hskip1pt T (b_{k,s})
\zeta \Big| > \lambda \Big\}} \\ & \le & \sum_{s=1}^\infty
\sum_{i=1}^4
\sum_{j=0}^s \sum_{\begin{subarray}{c} k \equiv j \\
\!\!\!\!\!\!\! \mod s+1 \end{subarray}} 2^{-\gamma s} \varphi
\Big[ \Big( \sum_{r=0}^s p_{k+r} \Big) f \Big( \sum_{r=0}^s
p_{k+r} \Big) \Big] \\ & \le & \Big( \sum_{s=1}^\infty
\sum_{i=1}^4 \sum_{j=0}^s 2^{-\gamma s} \Big) \, \|f\|_1 \, = \, 4
\Big( \sum_{s=1}^\infty \frac{s+1}{2^{\gamma s}} \Big) \, \|f\|_1
= 4 \hskip1pt \mathrm{c}_\gamma \, \|f\|_1.
\end{eqnarray*}
This completes the argument for the off-diagonal terms of $b$.

\subsection{Conclusion}

Combining the results obtained so far in Sections \ref{S5} and
\ref{S6}, we obtain the weak type inequality announced in Theorem
A. The strong $L_p$ estimates follow for $1 < p < 2$ from the real
interpolation method, see e.g. \cite{PX2} for more information on
the real interpolation of noncommutative $L_p$ spaces. In the case
$2 < p < \infty$, our estimates follow from duality since our
size/smoothness conditions on the kernel are symmetric in $x$ and
$y$.

\begin{remark} \label{TypeIII}
\emph{Recent results in noncommutative harmonic analysis
\cite{J2,JP1,JP2,JX1} show the relevance of non-semifinite von
Neumann algebras in the theory. The definition of the
corresponding $L_p$ spaces (so called Haagerup $L_p$ spaces) is
more involved, see \cite{H,T1}. A well-known reduction argument
due to Haagerup \cite{H2} allows us to extend our strong $L_p$
estimates in Theorems A and B to functions $f: \R^n \to \M$ with
$\M$ a type III von Neumann algebra $\M$. Indeed, if $\sigma$
denotes the one-parameter unimodular group associated to
$(\Mn,\varphi)$, we take the crossed product $\mathcal{R} = \Mn
\rtimes_\sigma \mathrm{G}$ with the group $\mathrm{G} = \bigcup_{n
\in \N} 2^{-n} \Z$. According to \cite{H2}, $\mathcal{R}$ is the
closure of a union of finite von Neumann algebras $\bigcup_{k \ge
1} \Mn_k$ directed by inclusion. We know that our result holds on
$L_p(\Mn_k)$ for $1 < p < \infty$ and with constants independent
of $k$. Therefore, the same will hold on $L_p(\mathcal{R})$. Then,
using that $L_p(\Mn)$ is a (complemented) subspace of
$L_p(\mathcal{R})$, the assertion follows.}
\end{remark}

\begin{remark} \label{RemWeakHypo} \emph{According to the
classical theory, it seems that some hypotheses of Theorem A could
be weakened. For instance, the size condition on the kernel is not
needed for scalar-valued functions. Moreover, it is well-known
that the classical theory only uses Lipschitz smoothness on the
second variable to produce weak type $(1,1)$ estimates. Going even
further, it is unclear whether or not we can use weaker smoothness
conditions, like H\"ormander type conditions. Nevertheless, all
these apparently extra assumptions become quite natural if we
notice that all of them where used to produce our
pseudo-localization principle, a key point in the whole argument.
Under this point of view, we have just imposed the natural
hypotheses which appear around the $T1$ theorem. This leads us to
pose the following problem.}

\noindent \emph{\textbf{Problem.} Can we weaken the hypotheses on
the kernel as pointed above?}
\end{remark}

\begin{remark}
\emph{We believe that our methods should generalize if we replace
$\R^n$ by any other space of homogeneous type. In other words, a
metric space equipped with a non-negative Borel measure which is
doubling with respect to the given metric. More general notions
can be found in \cite{CW,MS1,MS2}. Of course, following recent
results by Nazarov/Treil/Volberg and Tolsa, it is also possible to
study extensions of non-doubling Calder{\'o}n-Zygmund theory in our
setting. It is not so clear that the methods of this paper can be
easily adapted to this case.}
\end{remark}

\section{Operator-valued kernels}
\label{S7}

We now consider Calder{\'o}n-Zygmund operators associated to
operator-valued kernels $k: \R^{2n} \setminus \Delta \to \M$
satisfying the canonical size/smoothness conditions. In other
words, we replace the absolute value by the $\M$-norm, see the
Introduction for details. We begin by constructing certain bad
kernels which show that there is no hope to extend Theorem A in
full generality to this context. Then we obtain positive results
assuming some extra hypotheses.

\subsection{Negative results}

The origin of the counterexample we are constructing goes back to
a lack (well-known to experts in the field) of noncommutative
martingale transforms
$$\summ_k^{\null} df_k \mapsto \summ_k \xi_{k-1} df_k.$$ Indeed, the
boundedness of this operator on $L_p$ might fail when the
predictable sequence of $\xi_k$'s is operator-valued. Here is a
simple example. Let $\Mn$ be the algebra of $m \times m$ matrices
equipped with the standard trace $\mathrm{tr}$ and consider the
filtration $\Mn_1, \Mn_2, \ldots, \Mn_m$, where $\Mn_s$ denotes
the subalgebra spanned by the matrix units $e_{ij}$ with $1 \le
i,j \le s$ and the matrix units $e_{kk}$ with $k > s$.
\begin{itemize}
\item If $1 < p < 2$, we take $f = \sum_{k=2}^m e_{1k}$ and $\xi_k
= e_{k1}$, so that $$\Big\| \summ_k^{\null} \xi_{k-1} df_k
\Big\|_p = (m-1)^{1/p} \gg \sqrt{m-1} = \Big\| \summ_k df_k
\Big\|_p.$$

\item If $2 < p < \infty$, we take $f = \sum_{k=2}^m e_{k-1,k}$
and $\xi_k = e_{1k}$, so that $$\Big\| \summ_k^{\null} df_k
\Big\|_p = (m-1)^{1/p} \ll \sqrt{m-1} = \Big\| \summ_k \xi_{k-1}
df_k \Big\|_p.$$
\end{itemize}

Letting $m \to \infty$, we see that $L_p$ boundedness might fail
for any $p \neq 2$ even having $L_2$ boundedness. Our aim is to
prove that the same phenomenon happens in the context of singular
integrals with operator-valued kernels. The examples above show us
the right way to proceed. Namely, we shall construct a similar
operator using Littlewood-Paley type arguments. Note that a dyadic
martingale approach is also possible here, but this would give
rise to certain operators having non-smooth kernels and we want to
show that smoothness does not help in this particular case.

Let $\mathcal{S}_\R$ be the Schwarz class in $\R$ and consider a
non-negative function $\psi$ in $\mathcal{S}_\R$ bounded above by
$1$, supported in $1 \le |\xi| \le 2$ and identically $1$ in $5/4
\le |\xi| \le 7/4$. Define
$$\psi_k(\xi) = \psi(2^{-k} \xi).$$ Let $\Psi$ denote the inverse
Fourier transform of $\psi$, so that $\widehat{\Psi} = \psi$. If
we construct the functions $\Psi_k(x) = 2^k \hskip1pt \Psi(2^k
x)$, we have $\widehat{\Psi}_k = \psi_k$ and we may define the
following convolution-type operators
\begin{eqnarray*}
T_1 f (x) & = & \sum_{k \ge 1} e_{k1} \hskip1pt \Psi_k \! * \! f,
\\ T_2 f (x) & = & \sum_{k \ge 1} e_{1k} \hskip1pt \Psi_k \! * \!
f.
\end{eqnarray*}
In this case we are taking $\M = \mathcal{B}(\ell_2)$ and both
$T_1$ and $T_2$ become contractive operators in $L_2(\Mn)$.
Indeed, let $\mathcal{F}_\Mn = \mathcal{F}_\R \otimes
id_{L_2(\M)}$ denote the Fourier transform on $L_2(\Mn)$.
According to Plancherel's theorem, $\mathcal{F_A}$ is an isometry
and the following inequality holds
\begin{eqnarray*}
\|T_1f\|_2 & = & \Big\| \sum_{k \ge 1} e_{k1} \widehat{\Psi_k * f}
\Big\|_2 \ = \ \Big( \sum_{k \ge 1} \big\| e_{k1} \psi_k
\widehat{f} \hskip1pt \big\|_2^2 \Big)^{\frac12} \\ & \le & \Big(
\sum_{k \ge 1} \, \int_{2^k \le |\xi| \le 2^{k+1}}
|\widehat{f}(\xi)|^2 \, d\xi \Big)^{\frac12} \ \le \ \|f\|_2.
\end{eqnarray*}
The same argument works for $T_2$. Now we show that the kernels of
$T_1$ and $T_2$ also satisfy the expected size and smoothness
conditions. These are convolution-type kernels given by
$$k_1(x,y) = \sum_{k \ge 1} e_{k1} \Psi_k(x-y) \quad \mbox{and}
\quad k_2(x,y) = \sum_{k \ge 1} e_{1k} \Psi_k(x-y).$$ We clearly
have $$\begin{array}{rclcl} \|k_1(x,y)\|_\M & = & \displaystyle
\Big\| \sum_{k \ge 1} e_{k1} \Psi_k(x-y) \Big\|_\M & = &
\displaystyle \Big( \sum_{k \ge 1} |\Psi_k(x-y)|^2
\Big)^{\frac12}, \\ \|k_2(x,y)\|_\M & = & \displaystyle \Big\|
\sum_{k \ge 1} e_{1k} \Psi_k(x-y) \Big\|_\M & = & \displaystyle
\Big( \sum_{k \ge 1} |\Psi_k(x-y)|^2 \Big)^{\frac12}.
\end{array}$$ Therefore, for the size condition it suffices to see
that
\begin{equation} \label{size2}
\Big( \sum_{k \in \Z} |\Psi_k(x)|^2 \Big)^{\frac12} \lesssim
\frac{1}{|x|}.
\end{equation}
Similarly, using the mean value theorem in the usual way, the
condition
\begin{equation} \label{smooth2}
\Big( \sum_{k \in \Z} |\Psi_k'(x)|^2 \Big)^{\frac12} \lesssim
\frac{1}{|x|^2}
\end{equation}
implies Lipschitz smoothness for any $0 < \gamma \le 1$. The proof
of \eqref{size2} and \eqref{smooth2} is standard. Namely, since
$\Psi$ and $\Psi'$ belong to the Schwarz class $\mathcal{S}_\R$,
there exist absolute constants $\mathrm{c}_1$ and $\mathrm{c}_2$
such that $$|\Psi(x)| \le \mathrm{c}_1 \min \Big\{ 1,
\frac{1}{|x|^2} \Big\} \quad \mbox{and} \quad |\Psi'(x)| \le
\mathrm{c}_2 \min \Big\{ 1, \frac{1}{|x|^3} \Big\}.$$ If $2^{-j}
\le |x| < 2^{-j+1}$, we find the estimate
$$\Big( \sum_{k \in \Z} |\Psi_k(x)|^2 \Big)^{\frac12} \le \Big(
\mathrm{c}_1 \sum_{k \le j} 2^{2k} + \mathrm{c}_1 |x|^{-4} \sum_{k
> j} 2^{-2k} \Big)^{\frac12} \lesssim \Big( 2^{2j} +
\frac{1}{2^{2j}|x|^4} \Big)^\frac12 \lesssim \frac{1}{|x|}.$$
Similarly, using that $\Psi_k'(x) = 2^{2k} \Psi'(2^kx)$, we have
$$\Big( \sum_{k \in \Z} |\Psi_k'(x)|^2 \Big)^{\frac12} \le \Big(
\mathrm{c}_2 \sum_{k \le j} 2^{4k} + \mathrm{c}_2 |x|^{-6} \sum_{k
> j} 2^{-2k} \Big)^{\frac12} \lesssim \Big( 2^{4j} +
\frac{1}{2^{2j}|x|^6} \Big)^\frac12 \lesssim \frac{1}{|x|^2}.$$
Thus, $T_1$ and $T_2$ are bounded on $L_2(\Mn)$ with
operator-valued kernels satisfying the standard size and
smoothness conditions. Now we shall see how the boundedness on
$L_p(\Mn)$ fails for $p \neq 2$. By definition, we know that
\begin{itemize}
\item $\psi_k$ is supported by $2^k \le |\xi| \le 2^{k+1}$.

\item $\psi_k$ is identically $1$ in $5 \cdot 2^k/4 \le |\xi| \le
7 \cdot 2^k/4$.
\end{itemize}
If $\mathcal{I}_0 = [5/4,7/4]$ and $\mathcal{I}_k = \mathcal{I}_0
+ \frac32 (2^k-1)$, it is easily seen that
\begin{equation} \label{absorption}
\psi_k 1_{\mathcal{I}_k} = 1_{\mathcal{I}_k}
\end{equation}
for all nonnegative integer $k$. Now we are ready to show the
behavior of $T_1$ and $T_2$ on $L_p$. Indeed, let us fix an
integer $m \ge 1$ and let $g_k$ be the inverse Fourier transform
of $1_{\mathcal{I}_k}$ for $1 \le k \le m$. Then we set
$$f_1 = \sum_{k=1}^{m} e_{1k} \hskip1pt g_k \quad
\mbox{and} \quad f_2 = \sum_{k=1}^{m} e_{kk} \hskip1pt g_k.$$ By
\eqref{absorption} we have $\Psi_j * g_k = \delta_{jk} g_k$ for $1
\le k \le m$. Moreover, $$\widehat{g}_j(\xi) = \widehat{g}_k \Big(
\xi + \frac32 (2^k - 2^j) \Big) \Rightarrow |g_j(x)| = |g_k(x)|.$$
These observations allow us to obtain the following identities
$$\begin{array}{rclclcl} \displaystyle
\frac{\|T_1f_1\|_p}{\|f_1\|_p} & = & \frac{\displaystyle \Big\|
\sum_{k=1}^m e_{kk} g_k \Big\|_p}{\displaystyle \Big\|
\sum_{k=1}^m e_{1k} g_k \Big\|_p} & = & \frac{\displaystyle \Big\|
\Big( \sum_{k=1}^m |g_k|^p \Big)^{\frac1p} \Big\|_p}{\displaystyle
\Big\| \Big( \sum_{k=1}^m |g_k|^2 \Big)^{\frac12} \Big\|_p} & = &
m^{\frac1p - \frac12} \, \|g_1\|_p, \\ [29pt] \displaystyle
\frac{\|T_2f_2\|_p}{\|f_2\|_p} & = & \frac{\displaystyle \Big\|
\sum_{k=1}^m e_{1k} g_k \Big\|_p}{\displaystyle \Big\|
\sum_{k=1}^m e_{kk} g_k \Big\|_p} & = & \frac{\displaystyle \Big\|
\Big( \sum_{k=1}^m |g_k|^2 \Big)^{\frac12} \Big\|_p}{\displaystyle
\Big\| \Big( \sum_{k=1}^m |g_k|^p \Big)^{\frac1p} \Big\|_p} & = &
m^{\frac12 - \frac1p} \, \|g_1\|_p.
\end{array}$$
Therefore, letting $m \to \infty$ we see that $T_1$ and $T_2$ are
not bounded on $L_p(\Mn)$ for $1 < p < 2$ and $2 < p < \infty$
respectively. Since we have seen that both are bounded on
$L_2(\Mn)$ and are equipped with \emph{good} kernels, we deduce
that Theorem A does not hold for $T_1$ and $T_2$. This is a
consequence of the matrix units we have included in the kernels of
our operators.

\begin{remark}
\emph{We refer to \cite{Me2} and \cite{NPTV} for a study of
paraproducts associated to operator-valued kernels. There it is
shown that certain classical estimates also fail when dealing with
noncommuting operator-valued kernels. The results in \cite{NPTV}
give new light to Carleson embedding theorem.}
\end{remark}

\subsection{The $L_\infty \to \mathrm{BMO}$ boundedness}

In what follows we shall work under the hypotheses of Theorem B.
In other words, with Calder{\'o}n-Zygmund operators which are
$\M$-bimodule maps bounded on $L_q(\Mn)$ and are associated to
operator-valued kernels satisfying the standard size/smoothness
conditions, see the Introduction for further details. Let us
define the noncommutative form of dyadic $\mathrm{BMO}$ associated
to our von Neumann algebra $\Mn$. According to \cite{Me,PX1}, we
may define the space $\mathrm{BMO}_{\! \Mn}$ as the closure of
functions $f$ in $L_{1,\mathrm{loc}}(\R^n; \M)$ with
$$\|f\|_{\mathrm{BMO}_{\! \Mn}} \, = \, \max \Big\{
\|f\|_{\mathrm{BMO}_{\! \Mn}^r}, \|f\|_{\mathrm{BMO}_{\! \Mn}^c}
\Big\} \, < \, \infty,$$ where the row and column $\mathrm{BMO}$
norms are given by
\begin{eqnarray*}
\|f\|_{\mathrm{BMO}_{\! \Mn}^r} & = & \sup_{Q \in \Q} \Big\| \Big(
\frac{1}{|Q|} \int_Q \big( f(x) - f_Q \big) \big( f(x) - f_Q
\big)^* \, dx \Big)^{\frac12} \Big\|_{\M}, \\
\|f\|_{\mathrm{BMO}_{\! \Mn}^c} & = & \sup_{Q \in \Q} \Big\| \Big(
\frac{1}{|Q|} \int_Q \big(f(x) - f_Q \big)^* \big(f(x) - f_Q \big)
\, dx \Big)^{\frac12} \Big\|_{\M}.
\end{eqnarray*}
In order to extend our pseudo-localization result to the framework
of Theorem B, we shall need to work with the identity $\1_\Mn$ and
show that $T^* \1_\Mn$ belongs to the noncommutative form of
$\mathrm{BMO}$. In fact, the (still unpublished) result below due
to Tao Mei \cite{TaoMei} gives much more.

\begin{theorem} \label{TaoMei}
If $T$ is as above, then $$\|T \! f\|_{\mathrm{BMO}_{\! \Mn}} \
\le \ \mathrm{c}_{n,\gamma} \hskip1pt \|f\|_\Mn.$$
\end{theorem}

Mei's argument for Theorem \ref{TaoMei} is short and nice for
$q=2$. The case $q \neq 2$ requires the noncommutative analog of
John-Nirenberg theorem obtained by Junge and Musat in \cite{JM}.

\begin{remark} \label{bestconstante}
\emph{Let us fix an index $q < p < \infty$. By a recent result of
Musat \cite{Mu} adapted to our setting by Mei \cite{Me}, we know
that $L_q(\Mn)$ and $\mathrm{BMO}_{\! \Mn}$ form an interpolation
couple. Moreover, both the real and complex methods give the
isomorphism
$$\big[ \mathrm{BMO}_{\! \Mn}, L_q(\Mn) \big]_{q/p} \simeq
L_p(\Mn)$$ with constant $\mathrm{c}_p \sim p$ for $p$ large. The
proof of the latter assertion was achieved in \cite{JM}, refining
the argument of \cite{Mu}. In particular, the $L_p$ estimates
announced in Theorems A and B automatically follow from Theorem
\ref{TaoMei} combined with Musat's interpolation. Although this
approach might look much simpler, the proof of the necessary
interpolation results from \cite{Mu} and of the noncommutative
John-Nirenberg theorem (used in Mei's argument) are also quite
technical.}
\end{remark}

\begin{remark}
\emph{It also follows from Theorem \ref{TaoMei} that the problem
posed in Remark \ref{RemWeakHypo} is only interesting for weak
type inequalities. Indeed, if we are given a kernel with no size
condition and only satisfying the H\"ormander smoothness condition
in the second variable, then we may obtain the strong $L_p$
estimates provided by Theorems A and B for $1 < p \le 2$. We just
need to apply Mei's argument for Theorem \ref{TaoMei} (which works
under these weaker assumptions) to the adjoint mapping and dualize
backwards. A similar argument holds for H\"ormander smooth kernels
in the first variable and $2 \le p < \infty$.}
\end{remark}

\subsection{Proof of Theorem B}

Before proceeding with the argument, we set some preliminary
results. According to Theorem \ref{TaoMei} and the symmetry of the
conditions on the kernel, we know that $T^* \1_\Mn$ belongs to
$\mathrm{BMO}_{\! \Mn}$. In the following result, we shall write
$\mathcal{H}_1$ for the Hardy space associated to the dyadic
filtration on $\R^n$. That is, the predual of dyadic
$\mathrm{BMO}$, see \cite{GA}.

\begin{lemma} \label{cancellem}
If $T$ is as above and $T^* \1_\Mn = 0$, then
$$\int_{\R^n} T \! f (x) \, dx = 0 \quad \mbox{for any} \quad f
\in \mathcal{H}_1.$$
\end{lemma}

\dem Since $T^* \1_\Mn = 0$ vanishes as an element of
$\mathrm{BMO}_{\! \Mn}$, we will have
\begin{equation} \label{standardform}
\tau \Big( \int_{\R^n} T \! \phi (x) \, dx \Big) = \big\langle T
\! \phi, \1_\Mn \big\rangle = \big\langle \phi, T^* \1_\Mn
\big\rangle = 0
\end{equation}
for any $\phi \in \mathcal{H}_1(\Mn)$, the Hardy space associated
to the dyadic filtration $(\Mn_k)_{k \in \Z}$, see \cite{PX1} for
details and for the noncommutative analogue of Fefferman's duality
theorem $\mathcal{H}_1(\Mn)^* = \mathrm{BMO}_{\! \Mn}$. Given any
projection $q \in \M_\pi$ of finite trace and $f \in
\mathcal{H}_1$, it is clear that $\phi = f q \in
\mathcal{H}_1(\Mn)$. In particular, using $\M$-modularity again
$$\tau \Big( \int_{\R^n} T \! \phi(x) \, dx \Big) = \tau \Big( q
\int_{\R^n} T \! f(x) \, dx \Big) = 0$$ for any such projection
$q$. Clearly, this immediately implies the assertion. \fin

\begin{lemma} \label{NClocaliz}
Let $T$ be as above for $q=2$ and $L_2(\Mn)$-normalized. Then,
given $x_0 \in \R^n$ and $r_1, r_2 > 0$ with $r_2 > 2 \hskip1pt
r_1$, the following estimate holds for any pair $f,g$ of bounded
scalar-valued functions respectively supported by
$\mathsf{B}_{r_1}(x_0)$ and $\mathsf{B}_{r_2}(x_0)$
$$\Big\| \int_{\R^n} T \! f(x) \hskip1pt g(x) \, dx \Big\|_\M
\le \mathrm{c}_n \hskip1pt r_1^n \hskip1pt \log(r_2/r_1) \hskip1pt
\|f\|_\infty \|g\|_\infty.$$
\end{lemma}

\dem We proceed as in the proof of the localization estimate given
in Paragraph \ref{3AR}. Let $\mathsf{B}$ denote the ball
$\mathsf{B}_{3r_1/2}(x_0)$ and consider a smooth function $\rho$
identically $1$ on $\mathsf{B}$ and $0$ outside
$\mathsf{B}_{2r_1}(x_0)$. Taking $\eta = 1-\rho$, we may decompose
$$\Big\| \int_{\R^n} T \! f(x) \hskip1pt g(x) \, dx \Big\|_\M
= \Big\| \int_{\R^n} T \! f(x) \hskip1pt \rho g(x) \, dx \Big\|_\M
+ \Big\| \int_{\R^n} T \! f(x) \hskip1pt \eta g(x) \, dx
\Big\|_\M.$$ For the first term we adapt the commutative argument
using the convexity of the function $a \mapsto |a|^2$. Indeed, if
$\M$ embeds isometrically in $\mathcal{B(H)}$, it suffices to see
that $a \mapsto \langle a^*a \hskip1pt h, h \rangle_\mathcal{H}$
is a convex function for any $h \in \mathcal{H}$. However, this
follows from the identity $\langle a^*a \hskip1pt h, h
\rangle_\mathcal{H} = \| a h \|_\mathcal{H}^2$. As an immediate
consequence of this, we find the inequality
$$\Big| \frac{1}{|\mathsf{B}_{2r_1}(x_0)|}
\int_{\mathsf{B}_{2r_1}(x_0)} T \! f(x) \hskip1pt \rho g(x) \, dx
\Big|^2 \le \frac{1}{|\mathsf{B}_{2r_1}(x_0)|}
\int_{\mathsf{B}_{2r_1}(x_0)} \big| T \! f(x) \hskip1pt \rho g(x)
\big|^2 \, dx.$$ This combined with $\M$-modularity gives
\begin{eqnarray*}
\lefteqn{\hskip-15pt \Big\| \int_{\R^n} T \! f(x) \hskip1pt \rho
g(x) \, dx \Big\|_\M} \\ & = & |\mathsf{B}_{2r_1}(x_0)| \hskip1pt
\Big\| \frac{1}{|\mathsf{B}_{2r_1}(x_0)|}
\int_{\mathsf{B}_{2r_1}(x_0)} T \! f(x) \hskip1pt \rho g(x) \, dx
\Big\|_\M \\ & \le & |\mathsf{B}_{2r_1}(x_0)| \hskip1pt \Big\|
\frac{1}{|\mathsf{B}_{2r_1}(x_0)|} \int_{\mathsf{B}_{2r_1}(x_0)}
\big| T \! f(x) \hskip1pt \rho g(x) \big|^2 \, dx \Big\|_\M^\frac12 \\
& = & \mathrm{c}_n \hskip1pt r_1^{n/2} \sup_{\|a\|_{L_2(\M)} \le
1} \Big( \int_{\R^n} \tau \Big[ 1_{\mathsf{B}_{2r_1}(x_0)}(x)
\big| T \! f(x) \hskip1pt \rho g(x) a \big|^2 \Big] \, dx
\Big)^\frac12 \\ & \le & \mathrm{c}_n \hskip1pt r_1^{n/2}
\sup_{\|a\|_{L_2(\M)} \le 1} \Big( \int_{\R^n} \tau \Big[ \big| T
\! (fa)(x) \big|^2 \Big] \, dx \Big)^\frac12 \|g\|_\infty \\ & \le
& \mathrm{c}_n \hskip1pt r_1^{n/2} \hskip1pt \sup_{\|a\|_{L_2(\M)}
\le 1} \|f\|_2 \hskip1pt \|a\|_2 \, \|g\|_\infty \ \le \
\mathrm{c}_n \hskip1pt r_1^n \hskip1pt \|f\|_\infty \hskip1pt
\|g\|_\infty,
\end{eqnarray*}
since $\mathrm{supp} f \subset \mathsf{B}_{r_1}(x_0)$. On the
other hand, the second term equals
$$\Big\| \int_{\R^n} T \! f(x) \hskip1pt \eta g(x) \, dx
\Big\|_\M = \Big\| \int_{\mathsf{B}_{r_2}(x_0) \setminus
\mathsf{B}} \Big( \int_{\mathsf{B}_{r_1}(x_0)} k(x,y) \hskip1pt
f(y) \, dy \Big) \hskip2pt \eta g(x) \, dx \Big\|_\M.$$ This term
is estimated exactly in the same way as in Paragraph \ref{3AR}.
\fin

\skeB As in the proof of Theorem A, we first observe that there is
no restriction by assuming that $q=2$. Indeed, according to
Theorem \ref{TaoMei} and Remark \ref{bestconstante}, it is easily
seen that boundedness on $L_q(\Mn)$ is equivalent to boundedness
on $L_2(\Mn)$. Moreover, we may assume that $f \in \Mn_{c,+}$ and
decompose it for fixed $\lambda \in \R_+$ applying the
noncommutative Calder{\'o}n-Zygmund decomposition. This gives rise to
$f = g + b$. The diagonal parts are estimated in the same way.
Indeed, since we have $\|g_d\|_2^2 \le 2^n \lambda \|f\|_1$, the
$L_2$-boundedness of $T$ suffices for the good part. On the other
hand, we use Lemma \ref{keylem} for the bad part $b_d$ in the
usual way. This reduces the problem to estimate $\zeta \hskip1pt
T(b_d) \hskip1pt \zeta$. By $\M$-bimodularity, we can proceed
verbatim with the argument given for this term in the proof of
Theorem A. Moreover, exactly the same reasoning leads to control
the off-diagonal part $b_{\mathit{off}}$. It remains to estimate
the term associated to $g_{\mathit{off}}$. By Lemma \ref{keylem}
one more time, it suffices to study the quantity
$$\lambda \hskip1pt \varphi \Big\{ \big| \zeta \hskip1pt
T(g_{\mathit{off}}) \hskip1pt \zeta \big| > \lambda \Big\}.$$ As
in the proof of Theorem A, we write $g_{\mathit{off}} = \sum_{k,s}
g_{k,s}$ as a sum of martingale differences and use
pseudo-localization. To justify our use of pseudo-localization we
follow the argument in Theorem \ref{pseudolocal} using
$\M$-bimodularity. This reduces the problem to study the validity
of the \emph{paraproduct argument} and of the \emph{shifted form
of the $T1$ theorem} for our new class of Calder{\'o}n-Zygmund
operators.

The paraproduct argument is simple. Indeed, since $T$ is
$\M$-bimodular, the same holds for $T^*$ so that $T^* \1_\Mn$
becomes an element of $\mathrm{BMO}_{\! \mathcal{Z_A}}$ where
$\mathcal{Z_A}$ denotes the center of $\Mn$. According to
\cite{TaoMei}, the dyadic paraproduct $\Pi_\xi$ associated to the
term $\xi = T^* \1_\Mn$ defines a bounded map on $L_2(\Mn)$.
Moreover, since it is clear that $\Pi_\xi$ is $\M$-bimodular, this
allows us to consider the usual decomposition $T = T_0 +
\Pi_{\xi}^*$. Now following the argument in Paragraph \ref{PA},
with the characteristic functions $1_{\R^n \setminus
\Sigma_{f,s}}$ and $1_{\Omega_k}$ replaced by the corresponding
projections provided by Theorem \ref{pseudolocal}, we see that the
estimate of the paraproduct also reduces here to the shifted $T1$
theorem. At this point we make crucial use of the fact that $\xi =
T^* \1_\Mn$ is commuting, so that the same holds for
$\Delta_j(\xi)$ for all $j \in \Z$.

Let us now sketch the main (slight) differences that appear when
reproving the shifted $T1$ theorem for operator-valued kernels.
Lemma \ref{cancellem} will play the role of the cancellation
condition \eqref{T*1=0}. On the other hand, we also have at our
disposal the three auxiliary results (suitably modified) in
Paragraph \ref{3AR}. Namely, Cotlar lemma as it was stated there
will be used below with the only difference that we apply it over
the Hilbert space $\mathcal{H} = L_2(\Mn)$ instead of the
classical $L_2$. Regarding Schur lemma, it is evident how to adapt
it to the present setting. We just need to replace the Schur
integrals by
$$\mathcal{S}_1(x) = \int_{\R^n} \|k(x,y)\|_\M \, dy \quad
\mbox{and} \quad \mathcal{S}_2(y) = \int_{\R^n} \|k(x,y)\|_\M \,
dx.$$ We leave the reader to complete the straightforward
modifications in the original argument. Finally, Lemma
\ref{NClocaliz} given above is the counterpart in our context of
the localization estimate that we use several times in the proof
of Theorem A. Once these tools are settled, the proof follows
verbatim just replacing the absolute value $| \cdot |$ by the norm
$\| \cdot \|_\M$ when corresponds. Maybe it is also worthy of
mention that the two instances in the proof of the shifted $T1$
theorem where the Lebesgue differentiation theorem is mentioned,
we should apply its noncommutative analog from \cite{Me}. This
completes the proof. \fin

\begin{remark}
\emph{After Theorem B, it is also natural to wonder about a
vector-valued noncommutative Calder{\'o}n-Zygmund theory. Let us be
more precise, if the von Neumann algebra $\M$ is hyperfinite,
Pisier's theory \cite{Pis} allows us to consider the spaces
$L_p(\Mn;\mathrm{X})$ with values in the operator space
$\mathrm{X}$. Here it is important to recall that we must impose
on $\mathrm{X}$ an operator space structure since a Banach space
structure is not rich enough. Then, we can consider vector-valued
noncommutative singular integrals and study for which operator
spaces we obtain weak type $(1,1)$ and/or strong type $(p,p)$
inequalities. Of course, this is closely related to the geometry
of the operator space in question and in particular to the notion
of $\mathrm{UMD}_p$ operator spaces, also defined by Pisier. In
this context a great variety of problems come into scene, like the
independence of the $\mathrm{UMD}_p$ condition with respect to $p$
(see \cite{Mu2} for some advances) or the operator space analog of
Burkholder's geometric characterization of the $\mathrm{UMD}$
property in terms of $\zeta$-convexity \cite{Bu}.}
\end{remark}

\begin{remark}
\emph{Another related problem is the existence of $T1$ type
theorems for Calder{\'o}n-Zygmund operators associated to
operator-valued kernels. Here we should mention the closely
related works of Hyt\"onen \cite{Hy} and Hyt\"onen/Weis
\cite{HW2}. Namely given two Banach spaces $\mathrm{X}$ and
$\mathrm{Y}$, they consider $\mathrm{X}$-valued functions and
operator valued kernels taking values in $\mathcal{B}(\mathrm{X},
\mathrm{Y})$. Note that in our setting both $\mathrm{X}$ and
$\mathrm{Y}$ coincide with $L_2(\mathcal{M})$. The only drawback
of their approach in our setting is that, in the context of
general Banach spaces, they need to impose
$\mathcal{R}$-boundedness conditions on the kernel and it is
presumable that no such stronger conditions should be necessary in
our case.}
\end{remark}

\section*{Appendix A. On pseudo-localization}

\subsection*{\textnormal{A.1.} Applicability.}

We begin by analyzing how the pseudo-localization principle is
applied to a given $L_2$-function. At first sight, it is only
applicable to functions $f$ in $L_2$ satisfying that
$\mathsf{E}_m(f) = f_m = 0$ for some integer $m$. Indeed,
according to the statement of the pseudo-localization principle we
have
$$\mathrm{supp} \hskip1pt f \subset \bigcup_{k \in \Z}
\mathrm{supp} \hskip2pt df_{k+s} \subset \bigcup_{k \in \Z}
\Omega_k.$$ Given $\varepsilon = (\varepsilon_1, \varepsilon_2,
\ldots, \varepsilon_n)$ with $\varepsilon_j = \pm 1$ for $1 \le j
\le n$, let $$\R^n_{(\varepsilon)} = \Big\{ x \in \R^n \, \big| \
\mathrm{sgn} \hskip1pt x_j = \varepsilon_j \ \mbox{for} \ 1 \le j
\le n \Big\}$$ be the $n$-dimensional quadrant associated to
$\varepsilon$ and define $f_{(\varepsilon)}$ to be the restriction
of $f$ to such quadrant. If $f_m \neq 0$ for all $m \in \Z$, the
same will happen to $f_{(\varepsilon)}$ for some index
$\varepsilon$. Assume (with no loss of generality) that
$\varepsilon = (1,1,\ldots,1)$ or, in other words, that $f$ itself
is supported by the first quadrant. Let $\Lambda_f$ be the set of
negative $k$'s satisfying
$$\mathrm{supp} \hskip2pt df_{k+s} \neq \emptyset.$$ Our
hypothesis $f_m \neq 0$ for all $m \in \Z$ implies that
$\Lambda_f$ has infinitely many elements. According to the shift
condition, we know that $\Omega_k \neq \emptyset$ for each $k \in
\Lambda_f$ and therefore contains at least a cube in $\Q_k$, since
$\Omega_k$ is an $\mathcal{R}_k$-set. In fact, for $k$ small
enough the $\Q_k$-cube in the first quadrant closest to the origin
will be large enough to intersect the support of $f$. A moment of
thought gives rise to the conclusion that $\Omega_k$ contains such
cube for infinitely many negative $k$'s and
$$\Sigma_{f,s} = \bigcup_{k \in \Z} 9 \hskip1pt \Omega_k = \R^n.$$
Therefore, our result does not provide any information in this
case.

It is convenient to explain how to apply our result for an
arbitrary function $f$ in $L_2$ not satisfying the condition $f_m
= 0$. By homogeneity, we may assume that $\|f\|_2 = 1$. On the
other hand, if $\mathrm{supp} \hskip1pt f$ is not compact we
approximate $f$ by a compactly supported function $f_0$ such that
$\|f-f_0\|_2 \le \mathrm{c}_{n,\gamma} \hskip1pt s \hskip1pt 2^{-
\gamma s/4}$. This clearly reduces our problem to find the set
$\Sigma_{f,s}$ around the support of $f_0$. Next we decompose $f_0
= \sum_{1 \le j \le 2^n} f_j$, with $f_j$ being the restriction of
$f_0$ to the $j$-th quadrant and work independently with each of
these functions. In other words our localization problem reduces
to study functions $f$ in $L_2$ with compact support contained in
the first $n$-dimensional quadrant. Let $f$ be such a function and
take $Q$ to be the smallest dyadic cube containing the support of
$f$. We have $Q \in \Q_m$ for some integer $m$. Then we find $f_m
= \lambda 1_Q$ with $\lambda = \frac{1}{|Q|} \int_{\R^n} f(x)
\hskip1pt dx$ and thus we decompose
$$f = \big( f - \lambda 2^{- \gamma s/2} 1_{Q_s}
\big) + \lambda 2^{- \gamma s/2} 1_{Q_s} = f^1 + f^2$$ where $Q_s$
is a cube satisfying:
\begin{itemize}
\item $Q_s$ contains $Q$.

\item $Q_s$ is contained in a dyadic antecessor of $Q$.

\item The Lebesgue measure of $Q_s$ is $|Q_s| = 2^{\gamma s /2}
|Q|$.
\end{itemize}
It is clear that we have
$$\Big( \int_{\R^n} |f^2(x)|^2 dx \Big)^{\frac12} = \lambda
2^{- \gamma s/2} 2^{\gamma s/4} \sqrt{|Q|} = 2^{- \gamma s/4}
\|f_m\|_2 \le 2^{- \gamma s/4}.$$ Therefore, $f^2$ is small enough
for our aims. On the other hand, let $\widehat{Q}_s$ be the dyadic
$Q$-antecessor of generation $m-j_0$ with $j_0$ being the smallest
positive integer such that $j_0 \ge \gamma s /2 n$. In other
words, this cube is the smallest dyadic $Q$-antecessor containing
$Q_s$. If we set $m_0 = m-j_0$, we clearly have $f_{m_0}^1 = 0$.
When $k \le m_0-s$ we have $df^1_{k+s} = 0$ and $\mathrm{supp}
\hskip2pt df^1_{k+s} = \emptyset$ so that there is no set to
control. When $k+s > m_0$ we use $$\mathrm{supp} \hskip2pt
df^1_{k+s} \subset \widehat{Q}_s.$$ Hence, we may choose
$\Omega_k$ to be the smallest $\mathcal{R}_k$-set containing
$\widehat{Q}_s$. In the worst case $k = m_0-s+1$ we are forced to
take $\Omega_k$ as the $(s-1)$-th dyadic antecessor of
$\widehat{Q}_s$. That is, the $(j_0 + s -1)$-dyadic antecessor
$\widehat{Q}(j_0+s-1)$ of $Q$. This gives rise to the set
$$\Sigma_{f,s} = \bigcup_{k \in \Z} 9 \hskip1pt \Omega_k = 9
\hskip1pt \widehat{Q}(j_0+s-1) \subset 9 \cdot 2^{j_0+s} \hskip2pt
Q \sim 9 \cdot 2^{(1+ \frac{\gamma}{2n})s} \hskip2pt \mathrm{supp}
\hskip1pt f$$ and completes the argument for arbitrary $L_2$
functions. To conclude, we should mention that the dependance on
the $n$-dimensional quadrants, due to the geometry imposed by the
standard dyadic filtration, is fictitious. Indeed, we can always
translate the dyadic filtration, so that the role of the origin is
played by another point which leaves the support of $f$ in the
\emph{new} first quadrant.

\begin{Aremark}
\emph{Given a function $f$ in $L_2$ and a parameter $\delta \in
\R_+$, we have analyzed so far how to find appropriate sets
$\Sigma_{f,\delta}$ satisfying the localization estimate which
motivated our pseudo-localization principle $$\Big( \int_{\R^n
\setminus \Sigma_{f \! ,\delta}} |T \! f(x)|^2 \, dx
\Big)^{\frac12} \le \delta \Big( \int_{\R^n} |f(x)|^2 \, dx
\Big)^{\frac12}.$$ Reciprocally, given a set $\Sigma$ in $\R^n$
and $\delta \in \R_+$, it is quite simple to find functions
$f_{\Sigma, \delta}$ satisfying such estimate on $\R^n \setminus
\Sigma$. Indeed, let $s \ge 1$ be the smallest possible integer
satisfying $\mathrm{c}_{n, \gamma} s 2^{- \gamma s/4} \le \delta$
and write $\Sigma = \bigcup_{k \in \Z} 9 \Omega_k$ as a disjoint
union of $9$-dilations of maximal $\mathcal{R}_k$-sets. In this
case, any function of the form
$$f_{\Sigma,\delta} = \sum_{k \in \Z} 1_{\Omega_k} dg_{k+s}$$ with
$g \in L_2$ satisfies the hypotheses of our pseudo-localization
principle with $\Sigma$ as the final localization set. Indeed, we
have $d(f_{\Sigma,\delta})_{k+s} = 1_{\Omega_k} dg_{k+s}$ because
$1_{\Omega_k}$ is $(k+s)$-predictable and we deduce that
$f_{\Sigma,\delta}$ satisfies the shift condition.}
\end{Aremark}

\subsection*{\textnormal{A.2.}
Decreasing rate of singular integrals in the $L_2$ metric.}

As an immediate consequence of the pseudo-localization principle,
we can give a \emph{lower} estimate of how fast decreases a
singular integral far away from a set $\Sigma_f$ associated to
$f$. To be more specific, the following result holds.

\begin{Acorollary}
Let $f$ be in $L_2$ and define $$\Sigma_f = \bigcup_{k \in
\Z}^{\null} 9 \hskip1pt \Gamma_k \qquad \mbox{with} \qquad
\Gamma_k = \mathrm{supp} \hskip2pt df_k \in \mathcal{R}_k.$$ Then,
the following holds for any $\xi > 4$ $$\Big( \int_{\R^n \setminus
\hskip1pt \xi \Sigma_f} |T \! f(x)|^2 \, dx \Big)^{\frac12} \ \le
\ \mathrm{c}_{n,\gamma} \hskip1pt \xi^{-\gamma/4} \log \xi
\hskip1pt \Big( \int_{\R^n} |f(x)|^2 dx \Big)^{\frac12}$$ and any
$L_2$-normalized Calder{\'o}n-Zygmund operator with Lipschitz
parameter $\gamma$.
\end{Acorollary}

\dem Let $\Omega_k$ be the smallest $\mathcal{R}_k$-set containing
$\Gamma_{k+s}$. In the worst case, $\Gamma_{k+s}$ can be written
as a union $\bigcup_\alpha Q_\alpha$ of $\Q_{k+s}$-cubes. Taking
$\widehat{Q}_{\alpha}(s)$ to be the $s$-th dyadic antecessor of
$Q_\alpha$, we observe that $$\Omega_k \subset \bigcup_\alpha
\widehat{Q}_{\alpha}(s) \subset 2^{s+1} \Gamma_{k+s}.$$ Then we
construct
$$\Sigma_{f,s} = \bigcup_{k \in \Z} 9 \hskip1pt \Omega_k
\subset 2^{s+1} \Sigma_f,$$ and the theorem above automatically
gives $$\Big( \int_{\R^n \setminus \hskip1pt 2^{s+1} \Sigma_f} |T
\! f(x)|^2 \, dx \Big)^{\frac12} \ \le \ \mathrm{c}_{n,\gamma}
\hskip1pt s \hskip1pt 2^{- \gamma s/4} \hskip1pt \Big( \int_{\R^n}
|f(x)|^2 dx \Big)^{\frac12}.$$ Since this holds for every positive
integer $s$, the assertion follows. \fin

\begin{Aremark}
\emph{All the considerations in Paragraph A.1 apply to this
result.}
\end{Aremark}

\begin{Aremark}
\emph{This result might be quite far from being optimal, see
below.}
\end{Aremark}

\subsection*{\textnormal{A.3.} Atomic pseudo-localization
in $L_1$.}

\renewcommand{\theequation}{A.1}
\addtocounter{equation}{-1}

Maybe the oldest localization result was already implicit in the
Calder{\'o}n-Zygmund decomposition. Indeed, let $b$ denote the bad
part of $f$ associated to a fixed $\lambda > 0$ and let
$\Sigma_\lambda$ be the level set where the dyadic
Hardy-Littlewood maximal function $M_d f$ is bigger than
$\lambda$. Note that $b$ is supported by $\Sigma_\lambda$. Then,
we have
$$\int_{\R^n \setminus 2 \Sigma_\lambda} |T b(x)| \, dx
\le \mathrm{c}_n \summ_{j} \|b_j\|_1 \le \mathrm{c}_n \|f\|_1,$$
where the $b_j$'s are the atoms in which we decompose $b$. In
fact, this reduces to a well-known localization result for dyadic
atoms in $L_1$. Namely, let $a$ denote an atom supported by a
dyadic cube $Q_a$. Then, the mean-zero of $a$ gives the following
estimate for any $\xi > 2$
\begin{eqnarray} \label{atomiclocal}
\hskip30pt \int_{\R^n \setminus \hskip1pt \xi \hskip1pt Q_a}
|Ta(x)| \, dx \!\!\! & = & \!\!\! \int_{\R^n \setminus \hskip1pt
\xi \hskip1pt Q_a} \Big| \int_{\R^n} \big[ k(x,y) - k(x,c_{Q_a})
\big] \hskip1pt a(y) \, dy \Big| \, dx \\ \nonumber \!\!\! & \le &
\!\!\! \int_{\R^n \setminus \hskip1pt \xi \hskip1pt Q_a}
\int_{\R^n} \frac{|y-c_{Q_a}|^\gamma}{|x-y|^{n+\gamma}} \hskip1pt
|a(y)| \, dy \, dx \le \mathrm{c}_n \hskip1pt \xi^{-\gamma}
\hskip1pt \|a\|_1.
\end{eqnarray}
Note that the only condition on $T$ that we use is the
$\gamma$-Lipschitz smoothness on the second variable, not even an
a priori boundedness condition. Under these mild assumptions, we
may generalize \eqref{atomiclocal} in the language of our
pseudo-localization principle. Namely, the following result (maybe
known to experts) holds.

\begin{Atheorem} \label{THA}
Let us fix a positive integer $s$. Given a function $f$ in $L_1$
and any integer $k$, we define $\Omega_k$ to be the smallest
$\mathcal{R}_k$-set containing the support of $df_{k+s}$ and
consider the set
$$\Sigma_{f,s} = \bigcup_{k \in \Z} 3 \hskip1pt \Omega_k.$$
Then, we have for any Calder{\'o}n-Zygmund operator as above
$$\int_{\R^n \setminus \hskip1pt \Sigma_{f \! ,s}} |T \! f(x)| \,
dx \ \le \ \mathrm{c}_n \hskip1pt 2^{- \gamma s} \int_{\R^n}
|f(x)| dx.$$
\end{Atheorem}

\dem We may clearly assume that $f_m = 0$ for some integer $m$.
Namely, otherwise we can argue as in the previous paragraph to
deduce that $\Sigma_{f,s} = \R^n$ and the assertion is vacuous.
Define inductively
\begin{eqnarray*}
\mathsf{A}_1 & = & \mathrm{supp} \hskip2pt df_{m+1}, \\
\mathsf{A}_j & = & \mathrm{supp} \hskip2pt df_{m+j} \setminus
\Big( \bigcup_{w < j} \mathsf{A}_w \Big).
\end{eqnarray*}
Use that $\mathrm{supp} f \subset \bigcup_j \mathsf{A}_j$ and
pairwise disjointness of $\mathsf{A}_j$'s to obtain
\begin{eqnarray*}
\int_{\R^n \setminus \hskip1pt \Sigma_{f \! ,s}} |T \! f(x)| \, dx
& \le & \summ_j \, \sum_{\begin{subarray}{c} Q \in \Q_{m+j}
\\ Q \subset \mathsf{A}_j \end{subarray}} \, \int_{\R^n \setminus
\hskip1pt \Sigma_{f \! ,s}} |T (f 1_Q)(x)| \, dx \\ & = & \summ_j
\, \sum_{\begin{subarray}{c} Q \in \Q_{m+j}
\\ Q \subset \mathsf{A}_j \end{subarray}} \, \int_{\R^n \setminus
\hskip1pt \Sigma_{f \! ,s}} \Big| T \Big(1_Q \sum_{k=m+j}^{\infty}
df_k \Big)(x) \Big| \, dx.
\end{eqnarray*}
Let $\widehat{Q}_s$ be the $s$-th dyadic antecessor of $Q$. Since
$$Q \subset \mathsf{A}_j \subset \mathrm{supp} \hskip2pt
df_{m+j} \subset \Omega_{m+j-s}$$ and $Q \in \Q_{m+j}$, we deduce
$2^s Q \subset 3 \widehat{Q}_s \subset \Sigma_{f,s} \Rightarrow
\R^n \setminus \Sigma_{f,s} \subset \R^n \setminus 2^s Q$ and
$$\int_{\R^n \setminus \hskip1pt \Sigma_{f \! ,s}} |T \! f(x)| \,
dx \ \le \ \summ_j \, \sum_{\begin{subarray}{c} Q \in \Q_{m+j}
\\ Q \subset \mathsf{A}_j \end{subarray}} \, \int_{\R^n \setminus
\hskip1pt 2^s Q} \Big| T \Big(1_Q \sum_{k=m+j}^{\infty} df_k
\Big)(x) \Big| \, dx.$$ On the other hand, by \eqref{atomiclocal}
\begin{eqnarray*}
\int_{\R^n \setminus \hskip1pt \Sigma_{f \! ,s}} |T \! f(x)| \, dx
& \le & \mathrm{c}_n 2^{- \gamma s} \summ_j \,
\sum_{\begin{subarray}{c} Q \in \Q_{m+j} \\ Q \subset \mathsf{A}_j
\end{subarray}} \, \Big\| 1_Q \sum_{k=m+j}^{\infty} df_k \Big\|_1
\\ & = & \mathrm{c}_n 2^{- \gamma s} \summ_j \Big\|
1_{\mathsf{A}_j} \, \sum_{k=m+j}^{\infty} df_k \Big\|_1 \ = \
\mathrm{c}_n 2^{- \gamma s} \summ_j \|1_{\mathsf{A}_j} f\|_1.
\end{eqnarray*}
Using once more the pairwise disjointness of the $\mathsf{A}_j$'s
we deduce the assertion. \fin

\begin{Aremark} \emph{Here we should notice that the condition $f_m
= 0$ can not be removed as we did in the $L_2$ case and its
applicability is limited to this atomic setting. On the other
hand, if we try to use the argument of Theorem \ref{THA} for
$p=2$, we will find a nice illustration of why the ideas around
almost orthogonality that we have used in the paper come into
play. In the $L_1$ framework, almost orthogonality is replaced by
the triangle inequality.}
\end{Aremark}

\subsection*{\textnormal{A.4.} Other forms of pseudo-localization.}

\renewcommand{\theequation}{A.2}
\addtocounter{equation}{-1}

Once we have obtained results in $L_1$ and $L_2$, it is quite
natural to wonder about $L_p$ pseudo-localization for other values
of $p$. If we only deal with atoms, it easily seen that
\eqref{atomiclocal} generalizes to any $p > 1$ in the following
way
\begin{equation} \label{patomiclocal}
\Big( \int_{\R^n \setminus \hskip1pt \xi \hskip1pt Q_a} |Ta(x)|^p
\, dx \Big)^{\frac1p} \ \le \ \mathrm{c}_n \hskip1pt
\xi^{-(\gamma+ n/p')} \hskip1pt \Big( \int_{\R^n} |a(x)|^p \, dx
\Big)^{\frac1p}.
\end{equation}
This gives rise to two interesting problems:
\begin{itemize}
\item[i)] In Theorem \ref{THA} we showed that \eqref{atomiclocal}
generalizes to more general functions in $L_1$, those satisfying
$f_m =0$ for some integer $m$. On the other hand, as we have seen
in Paragraph A.1, the condition $f_m=0$ is not a serious
restriction for $p=2$, or any $p > 1$. Therefore, inequality
\eqref{patomiclocal} suggests that our pseudo-localization
principle might hold with $s \hskip1pt 2^{- \gamma s/4}$ replaced
by the better constant $2^{- (\gamma + n/2) s}$. However, this
result and its natural $L_p$ generalization are out of the scope
of this paper.

\vskip7pt

\item[ii)] Although the constant that we have obtained in our
pseudo-localization principle on $L_2$ might be far from being
optimal, it still makes a lot of sense to wonder whether or not
the corresponding interpolated inequality holds for $1 < p < 2$.
Below we give some guidelines which might lead to such a result.
We have not checked details, since the necessary estimates might
be quite technical, as those in the proof for $p=2$. All our ideas
below can be thought as problems for the interested reader.
\end{itemize}
\renewcommand{\theequation}{A.3}
\addtocounter{equation}{-1} The \emph{interpolated inequality}
that comes to mind is $$\Big( \int_{\R^n \setminus \Sigma_{f,s}}
|T \! f(x)|^p \, dx \Big)^{\frac1p} \, \le \,
\mathrm{c}_{n,\gamma} \frac{s \hskip1pt 2^{- \gamma s/4}}{(s 2^{3
\gamma s/4})^{\frac2p-1}} \, \Big( \int_{\R^n} |f(x)|^p \, dx
\Big)^{\frac1p}.$$ However, by the presence of $\Sigma_{f,s}$, a
direct interpolation argument does not apply and we need a more
elaborated approach. Namely, following the proof of our result in
$L_2$ verbatim, it suffices to find suitable upper bounds for
$\Phi_s$ and $\Psi_s$ in $\mathcal{B}(L_p)$. Here we might use
Rubio de Francia's idea of extrapolation and content ourselves
with a rough estimate (i.e. independent of $s$) for the norm of
these operators from $L_1$ to $L_{1,\infty}$. Of course, by real
interpolation this would give rise to the weaker inequality
\begin{equation} \label{AA22}
\Big( \int_{\R^n \setminus \Sigma_{f,s}} |T \! f(x)|^p \, dx
\Big)^{\frac1p} \, \le \, \mathrm{c}_{n,\gamma} \frac{(s \hskip1pt
2^{- \gamma s/4})^{2 - \frac2p}}{p-1} \, \Big( \int_{\R^n}
|f(x)|^p \, dx \Big)^{\frac1p}.
\end{equation}
However, this would be good enough for many applications. The
Calder{\'o}n-Zygmund method will be applicable to both $\Phi_s$ and
$\Psi_s$ if we know that their kernels satisfy a suitable
smoothness estimate. The lack of regularity of $\mathsf{E}_k$ and
$\Delta_{k+s}$ appears again as the main difficulty to overcome.
In this case, it is natural to wonder if the H\"ormander condition
$$\int_{|x| > 2 |y|} |k(x,y) - k(x,0)| \, dx \le
\mathrm{c}_{n,\gamma},$$ holds for the kernels of $\Phi_s$ and
$\Psi_s$. We believe this should be true. Anyway, a more in depth
application of Rubio's extrapolation method (which we have not
pursued so far) might be quite interesting here.

\begin{Aremark}
\emph{According to the classical theory \cite{St}, it is maybe
more natural to replace (in the shifted form of the $T1$ theorem)
the dyadic martingale differences $\Delta_{k+s}$ by a
Littlewood-Paley decomposition and the conditional expectations
$\mathsf{E}_k$ by their partial sums. This result will be surely
easier to prove since there is no lack of regularity as in the
dyadic martingale setting. This alternative approach to the
shifted $T1$ theorem might give rise to some sort of
pseudo-localization result in terms of Littlewood-Paley
decompositions. Although this is not helpful in the noncommutative
setting (by our dependance on Cuculescu's construction), it makes
the problem on the smoothness of the kernels of $\Phi_s$ and
$\Psi_s$ more accessible.}
\end{Aremark}

\begin{Aremark}
\emph{If the argument sketched above for inequality \eqref{AA22}
works, another natural question is whether results for $p
> 2$ can be deduced by duality. On one hand, the operator $\Phi_s$ behaves well
with respect to duality. In fact, the analysis of $\sum_k
\Delta_{k+s} T \mathsf{E}_k$ just requires (in analogy with the
$T1$ theorem) to assume first that we have $T1 = 0$. As pointed
above, this kind of cancellation conditions are only necessary for
$\Phi_s$, since the presence of the terms $id-\mathsf{E}_k$ in
$\Psi_s$ produce suitable cancellations. However, this is exactly
why the adjoint $$\Psi_s^* = \summ_k \Delta_{k+s} T_{4 \cdot
2^{-k}}^* (id - \mathsf{E}_k)$$ does not behave as expected. This
leaves open the problem for $p > 2$.}
\end{Aremark}

\section*{Appendix B. On Calder{\'o}n-Zygmund
decomposition}

\subsection*{\textnormal{B.1.} Weighted inequalities}

Given a positive function $f$ in $L_1$ and $\lambda \in \R_+$, let
us consider the Calder{\'o}n-Zygmund decomposition $f = g+b$
associated to $\lambda$. As pointed out and well-known, the most
significant inequalities satisfied by these functions are
$$\int_{\R^n} |g(x)|^2 \, dx \le 2^n \lambda \hskip1pt \|f\|_1
\quad \mbox{and} \quad \summ_j \int_{\R^n} |b_j(x)| \, dx \le 2
\hskip1pt \|f\|_1,$$ where the $b_j$'s are the atoms in which $b$
is decomposed. We already saw in Section \ref{S5} that these
inequalities remain true for the diagonal terms of the
noncommutative Calder{\'o}n-Zygmund decomposition. However, we do not
have at our disposal (see Paragraph B.2 below) such inequalities
for the off-diagonal terms. As we have explained in the
Introduction, our way to solve this lack has been to prove the
off-diagonal estimates
\begin{itemize}
\item $\big\|\zeta \hskip1pt T \big( \summ_k b_{k,s} \big)
\hskip1pt \zeta \big\|_1 \lesssim \alpha_s \|f\|_1$,

\vskip5pt

\item $\big\|\zeta \hskip1pt T \big( \summ_k g_{k,s} \big)
\hskip1pt \zeta \big\|_2^2 \lesssim \beta_s \hskip1pt \lambda
\hskip1pt \|f\|_1,$
\end{itemize}
for some fast decreasing sequences $\alpha_s, \beta_s$. The proof
of these estimates has exploited the properties of the projection
$\zeta$ in conjunction with our localization results. We have
therefore hidden the actual inequalities satisfied by the
off-diagonal terms which are independent of the behavior of $\zeta
\hskip1pt T(\cdot) \hskip1pt \zeta$. Namely, we have
\begin{itemize}
\item[a)] Considering the \emph{atoms} $$b_{k,s} = p_k (f-f_{k+s})
p_{k+s} + p_{k+s} (f-f_{k+s}) p_k$$ in $b_{\mathit{off}} =
\sum_{k,s} b_{k,s}$, we have for any positive sequence
$(\alpha_s)_{s \ge 1}$
$$\summ_s \summ_k \alpha_s \hskip1pt \|b_{k,s}\|_1 \ \lesssim \ \Big(
\summ_s s \hskip1pt \alpha_s \Big) \, \|f\|_1.$$

\item[b)] Considering the \emph{martingale differences}
$$g_{k,s} = p_k df_{k+s} q_{k+s-1} + q_{k+s-1} df_{k+s} p_k$$ in
$g_{\mathit{off}} = \sum_{k,s} g_{k,s}$, we have for any positive
sequence $(\beta_s)_{s \ge 1}$ $$\Big\| \summ_s \summ_k \beta_s
\hskip1pt g_{k,s} \Big\|_2^2 = \summ_s \summ_k \beta_s^2 \hskip1pt
\|g_{k,s}\|_2^2 \ \lesssim \ \Big( \summ_s \beta_s^2 \Big) \,
\lambda \, \|f\|_1.$$
\end{itemize}

As the careful reader might have noticed, the proof of these
estimates is implicit in our proof of Theorem A. It is still to be
determined whether these estimates for the weights $\alpha_s$ and
$\beta_s$ are sharp. On the other hand, it is also possible to
study weighted $L_p$ estimates for the off-diagonal terms of the
good part and $p>1$. We have not pursued any of these lines.

\subsection*{\textnormal{B.2.} On the lack of a classical $L_2$
estimate}

\renewcommand{\theequation}{B.1}
\addtocounter{equation}{-1}

The pseudo-localization approach of this paper has been motivated
by the lack of the key estimate $\|g\|_2^2 \lesssim \lambda
\hskip1pt \|f\|_1$ in the noncommutative setting. Although we have
not disproved such inequality so far, we end this paper by giving
some evidences that it must fail. Recalling from Section \ref{S5}
that the diagonal terms of $g$ satisfy the estimate $\|g_d\|_2^2
\lesssim \lambda \hskip1pt \|f\|_1$, it suffices disprove the
inequality
$$\|g_{\mathit{off}}\|_2^2 \lesssim \lambda \hskip1pt \|f\|_1.$$
By the original expression for $g_{\mathit{off}}$, we have
$$g_{\mathit{off}} = \sum_i \sum_{j < i} p_i f_i p_j + \sum_j
\sum_{i < j} p_i f_j p_j = \sum_k p_k f_k (\1_\Mn - q_{k-1}) +
(\1_\Mn - q_{k-1}) f_k p_k.$$ By orthogonality of the $p_k$'s and
the tracial property, it is easily seen that
\begin{eqnarray*}
\frac{1}{\lambda} \hskip1pt \|g_{\mathit{off}}\|_2^2 & = &
\frac{2}{\lambda} \hskip1pt \varphi \Big( \summ_k (\1_\Mn -
q_{k-1}) f_k p_k f_k (\1_\Mn - q_{k-1}) \Big) \\ & = & 2 \summ_k
\varphi \Big( \frac{f_k p_k f_k}{\lambda} \Big) + 2 \summ_k
\varphi \Big( \frac{q_{k-1} f_k p_k f_k q_{k-1}}{\lambda} \Big) \
= \ \mathsf{A} + \mathsf{B}.
\end{eqnarray*}
By the tracial property $$\mathsf{B} = \frac{2}{\lambda} \sum_k
\varphi(p_k f_k q_{k-1} f_k p_k).$$ Moreover, we also have
$$\|f_k^{\frac12} q_{k-1} f_k^{\frac12}\|_\infty =
\|q_{k-1} f_k q_{k-1}\|_\infty \le 2^n \hskip1pt \|q_{k-1} f_{k-1}
q_{k-1}\|_\infty \le 2^n \hskip1pt \lambda.$$ Thus, we find the
inequality $$\mathsf{B} \le \mathrm{c}_n \hskip1pt \sum_k
\varphi(p_k f_k) \le \mathrm{c}_n \hskip1pt \|f\|_1$$ and our
problem reduces to disprove
\begin{equation} \label{countexaim}
\Big\| \summ_k f_k p_k \Big\|_2^2 \lesssim \lambda \hskip1pt
\|f\|_1.
\end{equation}

As in the argument given in Section \ref{S7} to find a bad-behaved
noncommuting kernel, our motivation comes from a matrix
construction. Namely, let $\mathcal{A}$ be the algebra of $2m
\times 2m$ matrices equipped with the standard trace $\mathrm{tr}$
and consider the filtration $\mathcal{A}_1, \mathcal{A}_2, \ldots,
\mathcal{A}_{2m}$, where $\mathcal{A}_s$ denotes the subalgebra
spanned by the matrix units $e_{ij}$ with $1 \le i,j \le s$ and
the matrix units $e_{kk}$ with $k > s$. Let us set $\lambda = 1$
and define $$f = \sum_{i,j=1}^{2m} e_{ij}.$$ It is easily checked
that
\begin{eqnarray*}
q_1 & = & \chi_{(0,1]}(f_1) \ = \ \1_\mathcal{A}, \\ q_2 & = &
\chi_{(0,1]}(q_1 f_2 q_1) \ = \ \mbox{$\sum_{k > 2}$} \hskip1pt
e_{kk}, \\ q_3 & = & \chi_{(0,1]}(q_2 f_3 q_2) \ = \
\chi_{(0,1]}(q_2) = q_2, \\ q_4 & = & \chi_{(0,1]}(q_3 f_4 q_3) \
= \ \mbox{$\sum_{k>4}$} \hskip1pt e_{kk}, \\ q_5 & = &
\chi_{(0,1]}(q_4 f_5 q_4) \ = \ \chi_{(0,1]}(q_4) = q_4, \\ q_6 &
= & \ldots
\end{eqnarray*}
Hence, $p_{2k-1} = 0$ and $p_{2k} = e_{2k-1,2k-1} + e_{2k,2k}$.
This gives $$\sum_{k=1}^{2m} f_k p_k = \sum_{k=1}^m f_{2k} p_{2k}
= \sum_{k=1}^m \sum_{j=1}^{2k} e_{j,2k-1} + e_{j,2k}.$$ We have
$\lambda =1$ and it is clear that $\|f\|_1 = 2m$, while the $L_2$
norm is $$\Big\| \sum_{k=1}^{2m} f_k p_k \Big\|_2^2 = \sum_{k=1}^m
4k = 2 m(m+1).$$ Therefore, if we let $m \to \infty$ we see that
\eqref{countexaim} fails in this particular setting.

\noindent \textbf{Problem.} Adapt the construction given above to
the usual von Neumann algebra $\Mn$ of operator-valued functions
$f: \R^n \to \M$, equipped with the standard dyadic filtration.
This would disprove \eqref{countexaim} and thereby the inequality
$\|g\|_2^2 \lesssim \lambda \hskip1pt \|f\|_1$ in the
noncommutative setting.

\bibliographystyle{amsplain}

\

\

\hfill \noindent \textbf{Javier Parcet} \\
\null \hfill \textsc{Departamento de Matem{\'a}ticas} \\
\null \hfill \textsc{Instituto de Matem}{\scriptsize
{\'A}}\textsc{ticas y
F}{\scriptsize {\'I}}\textsc{sica Fundamental} \\
\null \hfill \textsc{Consejo
Superior de Investigaciones Cient}{\scriptsize {\'I}}\textsc{ficas} \\
\null \hfill Depto de Matem{\'a}ticas, Univ. Aut{\'o}noma de Madrid,
28049, Spain \\ \null
\hfill\texttt{javier.parcet@uam.es}

\end{document}